\documentclass[12pt, leqno]{article}
\usepackage{amstex}
\usepackage{amssymb}
\usepackage{epsfig}

\newtheorem{thm}{Theorem}[section]
\newtheorem{df}[thm]{Definition}
\newtheorem{prop}[thm]{Proposition}
\newtheorem{cor}[thm]{Corollary}
\newtheorem{lem}[thm]{Lemma}
\newtheorem{ex}[thm]{Example}

\newcommand{\qed}{\hfill\ensuremath{\square}}
\newcommand{\tr}{\text{tr}}

\newcommand{\wed}{\wedge}

\newcommand{\ag}{(\alpha,a)}
\newcommand{\bh}{(\beta,b)}
\newcommand{\GG}{\Gamma\times \mathcal{A}}
\newcommand{\A}{\mathcal{A}}

\newcommand{\ld}{\lambda}
\newcommand{\ch}{\text{ch}}

\newcommand{\m}{\mu}
\newcommand{\n}{\nu}
\newcommand{\pr}{\partial}

\newcommand{\CN}{\mathbb{C}}
\newcommand{\larw}{\longrightarrow}
\newcommand{\arw}{\rightarrow}
\newcommand{\daw}{\downarrow}

\newcommand{\G}{\Gamma}

\numberwithin{equation}{section}

\begin{document}
\title{\textbf{Graded Lie Superalgebras, Supertrace Formula, 
and Orbit Lie Superalgebras}}
\author{
%Victor G. Kac \begin{thanks}
%{Supported in part by NSF grant DMS-9622870.}
%\end{thanks},
Seok-Jin Kang\begin{thanks}
{Supported in part by Basic Science Research Institute Program, 
Ministry of Education of Korea, BSRI-98-1414,
and GARC-KOSEF at Seoul National University.}
\end{thanks}
and  
Jae-Hoon Kwon$^{\diamond}$ \cr 
{ } \cr
%{$^{*}$ Department of Mathematics} \cr
%{Massachusetts Institute of Technology} \cr
%{Cambridge, MA 02139, U. S. A.} \cr
%{ } \cr
{$^{*}$Department of Mathematics} \cr
{Seoul National University}\cr
{Seoul 151-742, Korea} \cr
{ } \cr
{$^{\diamond}$Department of Mathematics} \cr
{Korea Air Force Academy} \cr
{Chong-Ju 363-849, Korea}  }
\date{\today}
\maketitle 

\begin{abstract}

Let $\Gamma$ be a countable abelian semigroup and $\A$ be a countable 
abelian group satisfying a certain finiteness condition. Suppose that a 
group $G$ acts on a $(\Gamma \times \A)$-graded Lie superalgebra
${\frak L}=\bigoplus_{\ag \in \GG} {\frak L}_{\ag}$ by Lie superalgebra
automorphisms preserving the $(\GG)$-gradation. In this paper, we 
show that the Euler-Poincar\'e principle yields the generalized 
denominator identity for ${\frak L}$ and derive a closed form formula
for the supertraces $\text {str}(g|{\frak L}_{\ag})$ for all $g\in G$,
$\ag \in \GG$. We discuss the applications of our supertrace formula
to various classes of infinite dimensional Lie superalgebras 
such as free Lie superalgebras and generalized Kac-Moody superalgebras. 
In particular, we determine the decomposition of free Lie superalgebras
into a direct sum of irreducible $GL(n) \times GL(k)$-modules,
and the supertraces of the Monstrous Lie superalgebras with group
actions. Finally, we prove that the generalized characters of Verma
modules and the irreducible highest weight modules 
over a generalized Kac-Moody superalgebra ${\frak g}$ 
corresponding to the Dynkin
diagram automorphism $\sigma$ 
are the same as the usual characters of 
Verma modules and irreducible highest weight modules over 
the orbit Lie superalgebra $\breve{\frak g}={\frak g}(\sigma)$
determined by $\sigma$. 
\end{abstract}

\baselineskip=16pt 

\vskip 2cm 

\section{Introduction}

In the past three decades, the theory of 
infinite dimensional Lie algebras and their representations 
has been the focus of extensive research activities 
due to its rich and significant applications to many areas of 
mathematics and mathematical physics. 
The most well-known examples are  {\it Kac-Moody algebras}
(\cite{K68}, \cite{Mo68}) and {\it generalized Kac-Moody algebras}
(\cite{B88}), and a lot of exciting new discoveries have been made
using the language of these infinite dimensional Lie algebras. 
For example, the {\it Macdonald identities} (\cite{Mac72}) were shown to be
equivalent to the {\it denominator identities} of 
affine Kac-Moody algebras (\cite{K74}), 
and the {\it Moonshine Conjecture} was proved by investigating 
the structure of a special kind of generalized Kac-Moody algebra
called the {\it Monster Lie algebra} (\cite{B92}). 

\vskip 2mm

On the other hand, since 1970's, the Lie superalgebras and their 
representations have emerged naturally as the fundamental algebraic 
structure behind several areas of mathematical physics. In \cite{K77},
V. G. Kac gave a comprehensive presentation of the mathematical 
theory of Lie superalgebras, and obtained an important classification 
theorem for finite dimensional simple Lie superalgebras over 
algebraically closed fields of characteristic zero (see also \cite{Sch}). 
The theories of Kac-Moody algebras and generalized Kac-Moody algebras
have been extended to those of {\it Kac-Moody superalgebras} (\cite{K78})
and {\it generalized Kac-Moody superalgebras} (\cite{Mi}, \cite{Ray}) 
following the outline given in \cite{K90}. In \cite{KW1} and \cite{KW2},
V. G. Kac and M. Wakimoto developed the representation theory of affine
Kac-Moody superalgebras, and demonstrated some interesting applications 
of affine Kac-Moody superalgebras to number theory. 

\vskip 2mm

In \cite{Ka94}, S.-J. Kang derived a closed form {\it root multiplicity 
formula} for generalized Kac-Moody algebras, and by applying his
formula to the Monster Lie algebra, he obtained some interesting
recursive relations among the coefficients of the elliptic 
modular function $j$. The main idea of \cite{Ka94} 
is that the Euler-Poincar\'e principle for the Lie algebra homology
can be interpreted as the denominator identity for Kac-Moody
algebras or generalized Kac-Moody algebras, and once we are given
the denominator identity, we can derive a root multiplicity formula.
This idea has been extended to the general 
infinite dimensional graded Lie algebras to yield a closed form 
{\it dimension formula} for the homogeneous subspaces (\cite{KaK96}), 
and a lot of interesting applications of our formula to various classes of 
infinite dimensional Lie algebras  have been investigated 
(see also \cite{Ka95} and \cite{KaK95}). 
In \cite{KK}, V. G. Kac and S.-J. Kang  
exploited this idea further to derive a closed form
{\it trace formula} for graded Lie algebras 
with group actions. In particular, by applying their trace formula
to the Monster simple group acting on the Monster Lie algebra,
they obtained a family of interesting recursive relations among the
coefficients of the {\it Thompson series}. 

\vskip 2mm

The Thompson series are
examples of {\it replicable functions} which are normalized $q$-series 
satisfying certain functional equations called the {\it replication 
formulae}. 
In \cite{Ka97}, S.-J. Kang showed that, given the denominator identity
for graded Lie superalgebras, one can derive a closed form 
{\it superdimension formula} for the homogeneous subspaces,
which can be considered as a unifying method of studying the
structure of infinite dimensional Lie superalgebras. 
He discussed a lot of interesting applications of our superdimension
formula to various classes of infinite dimensional Lie superalgebras
such as free Lie superalgebras, generalized Kac-Moody superalgebras,
and Monstrous Lie superalgebras. Furthermore, he characterized the
replicable functions in terms of  certain product identities,
and determined the root multiplicities of the Monstrous Lie superalgebras
associated with replicable functions. 

\vskip 2mm 

In this paper, we combine the main ideas of \cite{KK} and \cite{Ka97}
to derive a closed form {\it supertrace formula} for 
graded Lie superalgebras with group actions. 
More precisely, let $\Gamma$ be a countable (ususally 
infinite) abelian semigroup and $\A$ be a countable (usually finite)
abelian group such that every element $\ag \in \GG$ can be written
as a sum of elements in $\GG$ in only finitely many ways. 
Suppose that a group $G$ acts on a graded Lie superalgebra 
${\frak L}=\bigoplus_{\ag \in \GG} {\frak L}_{\ag}$ by Lie superalgebra
automorphisms preserving the $(\GG)$-gradation. Then, in Section 2, 
we show that the Euler-Poincar\'e principle 
for the Lie superalgebra homology yields a product identity
$$
\begin{aligned}
\prod_{\ag \in \GG} & \exp \left(-\sum_{k=1}^{\infty} 
\frac{1}{k} \text {str}(g^k|{\frak L}_{\ag}) E^{k \ag} \right)\\
& = 1-\sum_{\ag \in \GG} \text {str} (g| H({\frak L})_{\ag}) E^{\ag},
\end{aligned}
$$
which is called the {\it generalized
denominator identity} for $g\in G$.
By taking the logarithm of both sides of the generalized denominator identity
and using M\"obius inversion, we derive a closed form formula for the
{\it supertraces} $\text {str}(g|{\frak L}_{\ag})
=\psi(a) \tr(g|{\frak L}_{\ag})$
for all $g\in G$, $\ag\in \GG$ (Theorem \ref {Thm 2.3}):
$$\text {str}(g|{\frak L}_{\ag})=\sum \Sb d>0 \\ \ag = d (\tau, b) \endSb
\frac{1}{d} \mu(d) W_{g^d} (\tau, b),$$ 
where $\psi(a)$ is the {\it sign map} defined on $\A$ and $W_{g}(\tau,b)$
is the {\it Witt partition function} defined by (\ref{2.16}). 

\vskip 2mm
In Section 3, we present the standard homology theory for Lie superalgebras.
That is, we show that there exists a super-analogue of Koszul's complex 
for Lie superalgebras (Proposition \ref {Prop 3.1}). 
The $k$-th homology module $H_{k}({\frak L}, V)$ of the Lie superalgebra 
${\frak L}$ with coefficients in $V$ is defined by 
$$H_{k}({\frak L}, V)=\text {Tor}^{U}_{k}(V, \mathbb{C})
=H_{k}(V\otimes_{U} M),$$
where $M=(M_{k}, \pr_{k})$ is the super-analogue of Koszul's
complex given by (\ref{3.1}) and (\ref{3.2}), and $U=U({\frak L})$ denotes the
universal enveloping algebra of ${\frak L}$. 
In particular, if ${\frak L}(V)$ is the free Lie superalgebra 
generated by a superspace $V$, then we have $H_{1}({\frak L}(V))=V$
and $H_{k}({\frak L}(V))=0$ for all $k\ge 2$. 
Even though our presentation is a straightforward
generalization of the homology theory of Lie algebras 
(cf. \cite{CE}, \cite{Fu}), we 
decided to include the detailed treatment here, for 
we could not find it in any other literature. 

\vskip 2mm
Section 4 and Section 5 are devoted to the applications of our supertrace
formula to free Lie superalgebras. In Section 4, we  generalize
the classical {\it Witt formula} for free Lie algebras to obtain a closed
form supertrace formula for free Lie superalgebras with group actions. 
As an application, we compare the structure
of the free Lie algebras and free Lie superalgebras generated by the
same group representations (Proposition \ref {Prop 4.2}).
In Section 5, we consider the free Lie superalgebra ${\frak L}(V)$ 
generated by the natural $(k+l)$-dimensional representation $V$ of
the general linear Lie superalgebra $gl(k,l)$, and determine
the decomposition of ${\frak L}(V)$ into a direct sum of irreducible
$GL(k) \times GL(l)$-modules (Proposition \ref {Prop 5.3} and 
Proposition \ref {Prop 5.4}). 
The {\it Hook Schur functions} and the {\it Littlewood-Richardson 
coefficients} are the key ingredients of the decomposition. 

\vskip 2mm
In Section 6 and Section 7, 
we give a conjectural Kostant-type formula for the 
homology of the negative part of 
generalized Kac-Moody superalgebras, and derive 
a closed form supertrace formula for 
generalized Kac-Moody superalgebras with group actions preserving
the root space decomposition (Theorem \ref {Thm 6.5}). 
Moreover, we apply our supertrace
formula to determine the supertraces of the Monstrous Lie
superalgebras with group actions preserving the root space
decomposition (Proposition \ref{Prop 7.1} and Corollary \ref{Cor 7.2}). 
We believe that the product identities characterizing the 
replicable functions (\cite{Ka97}) can be understood most naturally as
the generalized denominator identities of some
Mosntrous Lie superalgebra with a suitable group action. 

\vskip 2mm
In Section 8 and Section 9, 
we prove that the generalized characters of Verma modules
and irreducible highest weight modules over a generalized
Kac-Moody superalgebra ${\frak g}$ corresponding to the Dynkin 
diagram automorphism $\sigma$ are the same as the usual characters 
of Verma modules and irreducible highest weight 
modules over the {\it orbit Lie superalgebra}
$\breve{\frak g}={\frak g}(\sigma)$ determined by $\sigma$
(Theorem \ref{Thm 9.3} and Corollary \ref{Cor 9.4}). 
This fact was first proved in \cite{FSS} for symmetrizable Kac-Moody algebras
under the linking condition in solving the fixed point problem 
in conformal field theory and was generalized 
to generalized Kac-Moody algebras 
without the linking condition in \cite{FRS}.
In this paper, we also show that 
the generalized denominator identity for a generalized
Kac-Moody superalgebra ${\frak g}$ for the Dynkin diagram automorphism 
$\sigma$ is the same as the usual denominator identity for the
orbit Lie superalgebra $\breve {\frak g}={\frak g}(\sigma)$ 
(Theorem \ref{Thm 9.5}). 
Our argument for generalized Kac-Moody superalgebras is similar  
to those in \cite{FSS} and \cite{FRS}. 
We observe that the generalized character
of a Verma module can be written as the product part of a 
generalized denominator identity and derive an explicit formula 
for the traces of $\sigma$ on each invariant subspaces
(which are not the same as root spaces in general) 
of ${\frak g}$ in terms of the dimensions of some orbit
Lie superalgebras. 
Therefore, we obtain a dimension formula for each homogeneous
subspaces of the invariant subalgebra, which is also a generalized 
Kac-Moody superalgebra with a natural grading arising from $\sigma$ 
(Proposition \ref{Prop 9.6}).

\vskip 1cm
\noindent
{\bf Acknowledgments.} The authors would like to express their
sincere gratitude to Professor Victor G. Kac for his interest 
in this work and many valuable discussions. 

\vskip 2cm

\section{Graded Lie Superalgebras and Supertrace Formula}

We recall some basic facts about Lie superalgebras (\cite{K77}, 
\cite{Ka97}, \cite{Sch}). 
Let $\A$ be a countable abelian group and 
suppose we have a bimultiplicative 
map $\theta : \A\times\A \longrightarrow \mathbb{C}^{\times}$
satisfying
\begin{equation}
\begin{aligned}
& \theta(a+b,c)= \theta(a,c)\theta(b,c), \\
& \theta(a,b+c)= \theta(a,b)\theta(a,c), \\
& \theta(a,b)\theta(b,a)=1\text{   for all  }a,b\in \A.
\end{aligned}
\end{equation} 
In particular, we have $\theta(a,a)=\pm 1$ for all $a\in \A$.
The  map $\theta$ satisfying the 
condition (2.1) is called a \textit{coloring map on $\A$}.
Define a map $\psi : \A\longrightarrow \{\pm1\}$ by $\psi(a)=\theta(a,a)$.
Then $\psi$ is a well-defined group homomorphism called the
\textit{sign map on $\A$ with respect to $\theta$}.
Let $\A_0=\{\,a\in \A\,|\,\psi(a)=1\,\}$ and
    $\A_1=\{\,a\in \A\,|\,\psi(a)=-1\,\}$.
Then we have a decomposition $\A=\A_0\cup \A_1$ and
the elements of $\A$ in $\A_0$ (resp. $\A_1$) are called 
\textit{even} (resp. \textit{odd}).

\vskip 2mm
A {\it $\theta$-colored superspace} is a pair
$(V, \theta)$, where $V=\bigoplus_{a\in \A} V_{a}$
is an $\A$-graded vector space and $\theta: \A \times 
\A \rightarrow \mathbb {C}^{\times}$ is a coloring
map on $\A$. 
The elements of $V_{a}$ are called
{\it even} (resp. {\it odd}) if $\psi(a)=1$ (resp. $\psi(a)=-1$).
For each $a\in \A$, we define the {\it superdimension} of $V_{a}$
to be
\begin{equation}
\text {sdim} V_{a}=\psi(a) \text {dim} V_{a}.
\end{equation}

\vskip 2mm
Similarly, we define a {\it $\theta$-colored superalgebra} to be a pair
$(U, \theta)$, where $U=\bigoplus_{a\in \A} U_a$ is 
an $\A$-graded associative algebra (i.e., 
$U_a U_b \subset U_{a+b}$ for all $a,b \in \A$)
and $\theta: \A \times \A \rightarrow 
\text {\bf C}^{\times}$ is a coloring map on $\A$. 
The {\it direct sum} of $\A$-graded superalgebras is defined 
in the usual way, but, for $\theta$-colored superalgebras
$U=\bigoplus_{a\in \A} U_{a}$ and 
$U'=\bigoplus_{a'\in \A} U'_{a}$,
we define the {\it tensor product} of $U$ and $U'$ to be the 
$\theta$-colored superspace 
$U\otimes U'$ with the natural $\A$-gradation
$$(U \otimes U')_{a+a'}=\text {Span} \{ u \otimes u' | \
u \in U_{a}, u'\in U'_{a'} \}$$
and the multiplication given by 
$$(u \otimes u') (v \otimes v')
=\theta(a',b) (uv \otimes u'v')$$
for $u\in U_a$, $v\in U_b$, $u'\in U'_{a'}$, $v'\in U'_{b'}$,
$a, a', b, b' \in \A$.

\vskip 2mm

\begin{df} \label {Def 2.1} \ {\rm 
A {\it $\theta$-colored Lie superalgebra}
is a $\theta$-colored superspace ${\frak L}=
\bigoplus_{a\in \A} {\frak L}_a$ 
together with a bilinear operation $[\ ,\ ]: {\frak L} \times {\frak L}
\rightarrow {\frak L}$ satisfying 
\begin{equation}
\begin{aligned}
& [{\frak L}_{a}, {\frak L}_{b}] \subset {\frak L}_{a+b}, \\
& [x,y]=-\theta(a,b) [y,x], \\
& [x, [y, z]]=[[x, y], z]+\theta(a,b) [y, [x,z]]
\end{aligned}
\end{equation}
for all $x\in {\frak L}_{a}$, $y\in {\frak L}_{b}$, $z\in {\frak L}$ 
and $a,b \in \A$. }
\end{df}

\vskip 2mm
Let ${\frak L}_{0}=\bigoplus_{a\in \A_0} {\frak L}_a$ and
${\frak L}_{1}=\bigoplus_{b\in \A_{1}} {\frak L}_b$. 
Then we have a decomposition 
${\frak L}={\frak L}_{0} \oplus {\frak L}_{1}$, and 
the homogeneous elements of ${\frak L}_{0}$ 
(resp. ${\frak L}_{1}$) are called {\it even} (resp. {\it odd}).

\vskip 2mm
The {\it universal enveloping algebra} of 
a $\theta$-colored Lie superalgebra ${\frak L}=\bigoplus_{a \in
\A} {\frak L}_{a}$ is the pair $(U({\frak L}), \iota)$,
where $U({\frak L})$ is a $\theta$-colored superalgebra and 
$\iota: {\frak L} \rightarrow U({\frak L})$ is a linear mapping satisfying 
$$\iota([x,y])=\iota(x)\iota(y)-\theta(a,b) \iota(y)\iota(x) \ \ 
\text {for} \ x \in {\frak L}_{a}, y \in {\frak L}_{b}$$
such that for any $\theta$-colored superalgebra 
$U=\bigoplus_{a\in \A} U_a$
and a linear mapping $j:{\frak L} \rightarrow U$ satisfying 
$$j([x,y])=j(x)j(y)-\theta(a,b) j(y)j(x) \ \ \text {for} \ 
x\in {\frak L}_{a}, y\in {\frak L}_{b},$$
there exits a unique homomorphism $\psi: U({\frak L}) \rightarrow U$ of 
$\theta$-colored superalgebras satisfying $\psi \circ \iota =j$.

\vskip 2mm
The uniqueness of the universal enveloping 
algebra $U({\frak L})$ can
be proved in the usual way. For the existence, consider the tensor algebra
${\mathcal T}({\frak L})
=\bigoplus_{k=0}^{\infty} {\frak L}^{\otimes k}$ of ${\frak L}$ and 
let ${\mathcal R}$ be the ideal of ${\mathcal T}({\frak L})$ 
generated by the elements of the form
$$[x,y]-x\otimes y + \theta(a,b) y\otimes x \ \ 
(x \in {\frak L}_a, y\in {\frak L}_b, \ a, b \in {\A}).$$
Then the quotient algebra $U({\frak L})={\mathcal T}({\frak L})
\big/ {\mathcal R}$ together with the
natural mapping $\iota:{\frak L} \rightarrow U({\frak L})$ 
is the universal enveloping algebra of ${\frak L}$. 

\vskip 2mm
For the structure of the 
universal enveloping algebra $U({\frak L})$ of ${\frak L}$, 
the following version of Poincar\'e-Birkhoff-Witt theorem 
is well-known.

\begin{thm} \label{Thm 2.1} {\rm (\cite{BMPZ}, \cite{K77}, \cite{Sch})}

Let $\A$ be a countable abelian group with a coloring map 
$\theta : \A \times \A \rightarrow \text {\bf C}^{\times}$
and let ${\frak L}=\bigoplus_{a\in \A} {\frak L}_{a}$ be a $\theta$-colored
Lie superalgebra. 
Suppose $X=\{x_{\alpha}| \ \alpha \in \Lambda \}$ {\rm ({\it resp. 
$Y=\{y_{\beta}| \ \beta \in \Omega \}$})} is a homogeneous basis of the
subspace ${\frak L}_{0}=\bigoplus_{a \in \A_{0}} {\frak L}_{a}$ 
{\rm ({\it resp. ${\frak L}_{1}=\bigoplus_{b \in \A_{1}} {\frak L}_{b}$
})}. 
Then the elements of the form 
\begin{equation}
x_{\alpha_1} x_{\alpha_2} \cdots x_{\alpha_k} \  y_{\beta_1} y_{\beta_2} 
\cdots y_{\beta_l} \ \ \text {with} \ \alpha_1 \le \cdots \le \alpha_k, \ 
\beta_1<\cdots < \beta_l 
\end{equation}
together with 1 form a basis of the universal enveloping algebra $U({\frak L})$
of ${\frak L}$. \qed
\end{thm}

\vskip 2mm
\noindent
{\it Remark.} \ The theory of colored Lie superalgebras can
be reduced to the theory of ordinary $\mathbb {Z}_2$-graded Lie
superalgebras by redefining the superbracket in a suitable way. 
More precisely, let ${\frak L}=({\frak L}, [\ , \ ])$ 
be a $\theta$-colored 
Lie superalgebra with a coloring map $\theta : \A \times \A 
\longrightarrow \mathbb{C}^{\times}$ and 
let $X=\{x_{i} | \ i=1,2, \cdots \}$ be an 
ordered set of generators of $\A$. 
We define the bimultiplicative maps $\delta : \A \times \A 
\longrightarrow \mathbb{C}^{\times}$ 
and $\tau : \A \times \A \longrightarrow \mathbb{C}^{\times}$ by 
$$\delta (a, b)  =(-1)^{\psi(a) \psi(b)} \theta(b,a)$$
and 
$$\tau (a,b)  = \prod_{i<j} \delta (x_{i}, x_{j})^{r_i s_j},$$
where $a=\sum r_{i} x_{i}, \ b=\sum s_i x_i \in \A$ 
\ $(r_{i}, s_{i} \in \mathbb{Z})$. 
We now introduce a new bracket $[\ , \ ]^{\tau} : {\frak L} \times 
{\frak L} \longrightarrow {\frak L}$ defined by 
\begin{equation*}
[x, y]^{\tau} = \tau(a,b) [x,y]
\end{equation*}
for $x\in {\frak L}_{a}$, $y\in {\frak L}_{b}$.
Then ${\frak L}^{\tau}=({\frak L}, [ \ , \ ]^{\tau})$ 
becomes an ordinary $\mathbb{Z}_{2}$-graded 
Lie superalgebra. 
Furthermore, there is a 1-1 correspondence between the set of 
$\A$-graded representations of the $\theta$-colored Lie superalgebra
${\frak L}=({\frak L}, [\ , \ ])$ and the set of $\A$-graded representations 
of the ordinary Lie superalgebra 
${\frak L}^{\tau}=({\frak L}, [\ , \ ]^{\tau})$
(cf. \cite{Pop}, \cite{Sch1}). 

\vskip 2mm

Let $\A$ be a countable (usually finite) abelian group
with a coloring map $\theta$, and let 
$\G$ be a countable (usually infinite) abelian semigroup
such that every element $(\alpha, a) \in \G \times \A$ 
can be written as a sum of elements in $\G \times \A$
in only finitely many ways. 
Consider a ($\GG$)-graded $\theta$-colored superspace 
$V=\bigoplus_{\ag\in\GG}V_{\ag}$ with $\dim V_{\ag} < \infty$
for all $\ag \in \GG$.
Suppose that a group $G$ acts on $V$ preserving the
($\GG$)-gradation.
We define the \textit{generalized character} or
the {\it graded trace} of $V$ for $g\in G$ by
\begin{equation}
\ch_{g}V=\sum_{(\alpha, a) \in \GG}
\tr(g|V_{\ag})e^{\ag},
\end{equation}
where $e^{\ag}$ are the basis elements of the semi-group
algebra $\mathbb{C}[\GG]$
with the multiplication given by 
$e^{\ag}e^{\bh}=e^{(\alpha+\beta,a+b)}$.

\vskip 2mm
On the other hand, we define 
the \textit{supertrace} of the homogeneous subspace 
$V_{\ag}$ for $g\in G$ by
\begin{equation}
\text {str}(g|V_{\ag})=\psi(a)\tr(g|V_{\ag}),
\end{equation}
and introduce another basis elements of $\mathbb{C}[\GG]$ by setting
$E^{\ag}=\psi(a)e^{\ag}$.
Then it is easy to verify that 
$E^{\ag}E^{\bh}=E^{(\alpha+\beta,a+b)}$.
We now define the \textit{generalized supercharacter} or 
the {\it graded supertrace} of $V$ for $g\in G$ by 
\begin{equation}
\text {sch}_{g}V=\sum_{\ag\in\GG}
\text {str}(g|V_{\ag})E^{\ag}.
\end{equation}
Note that $\text {sch}_{g}V$ is obtained from $\ch_{g}V$ by replacing
$\tr(g|V_{\ag})=\psi(a)
\text {str}(g|V_{\ag})$ and 
$e^{\ag}=\psi(a)E^{\ag}$.
Since $\psi(a)^2=1$, we have $\ch_{g}V=\text {sch}_{g}V$.

\vskip 2mm
Let ${\frak L}=\bigoplus_{\ag\in\GG}{\frak L}_{\ag}$
be a $\theta$-colored ($\GG$)-graded Lie superalgebra
with $\dim {\frak L}_{\ag} < \infty$ for all $\ag \in \GG$
and suppose that a group $G$ acts on ${\frak L}$ by
Lie superalgebra automorphisms 
preserving the $\GG$-gradation.
In this work, using the Euler-Poincar\'e principle, 
which is equivalent to the {\it generalized denominator
identity}, we will derive a closed form formula for the supertraces
$\text {str} (g|{\frak L}_{\ag})$ for all $g\in G$,  $\ag \in \GG$.

\vskip 2mm
Let $\mathbb{C}$ be the trivial one dimensional ${\frak L}$-module.
The homology modules $H_{k}({\frak L})=H_{k}({\frak L},\mathbb{C})$ are 
determined from the following {\it standard complex}:
\begin{equation}\label {2.8}
\cdots\rightarrow C_{k}({\frak L}) \stackrel{d_k}{\rightarrow}
C_{k-1}({\frak L})\stackrel{d_{k-1}}\rightarrow\cdots
\rightarrow C_{1}({\frak L})\stackrel{d_1}{\rightarrow} C_{0}({\frak L})
\rightarrow 0, 
\end{equation}     
where $C_{k}({\frak L})$ are defined by
\begin{equation*} 
C_{k}({\frak L})=\bigoplus_{p+q=k}
\Lambda^{p}({\frak L}_{0})\otimes S^{q}({\frak L}_{1})
\end{equation*}
and the differentials $d_k:C_{k}({\frak L})\rightarrow C_{k-1}({\frak L})$ 
are given by 
\begin{equation*}
\begin{aligned}
& d_{k}((x_1\wed\cdots\wed x_p)\otimes(y_1\cdots y_q)) \\
& =\sum_{1\leq s < t\leq p}(-1)^{s+t}
([x_s,x_t]\wed x_1\wed\cdots\wed \widehat{x_s}\wed
\cdots\wed\widehat{x_t}\wed\cdots\wed x_p)
\otimes(y_1\cdots y_q) \\
 &+\sum_{s=1}^{p}\sum_{t=1}^{q}(-1)^{s}
(x_1\wed\cdots\wed\widehat{x_s}\wed\cdots\wed x_p)\otimes
([x_s,y_t]y_1\cdots\widehat{y_t}\cdots y_q) \\
 &-\sum_{1\leq s < t\leq q}
([y_s,y_t]\wed x_1\wed\cdots\wed x_p)
\otimes(y_1\cdots\widehat{y_s}\cdots\widehat{y_t}\cdots y_q)
\end{aligned}
\end{equation*}
for $k\geq 2$, $x_i\in {\frak L}_{0}$, $y_j\in {\frak L}_{1}$ and $d_1=0$ 
 (cf. \cite{CE}, \cite{Fu}).
In Section 3, 
we will give a more detailed treatment of the homology of ${\frak L}$.
Since the spaces $C_k({\frak L})$ and the homology modules 
$H_k({\frak L})$ inherit
the ($\GG$)-gradation from that of ${\frak L}$, the action of $G$ and the
generalized supercharacters of $C_k({\frak L})$ and 
$H_k({\frak L})$ are well-defined.
Hence by the Euler-Poincar\'{e} principle, we obtain
\begin{equation}
\sum_{k=0}^{\infty}(-1)^{k}\text {sch}_{g}C_{k}({\frak L})
=\sum_{k=0}^{\infty}(-1)^{k}\text {sch}_{g}H_{k}({\frak L})
\ \ \text {for all} \ g\in G.
\end{equation}

\vskip 2mm 

Let 
\begin{equation}
\begin{aligned}
& C({\frak L})=\sum_{k=0}^{\infty}(-1)^{k}C_{k}({\frak L})
=\mathbb{C}\ominus L\oplus C_{2}({\frak L})\ominus\cdots , \\
& \Lambda({\frak L}_{0})=\sum_{k=0}^{\infty}(-1)^{k}
\Lambda^{k}({\frak L}_{0})
=\mathbb{C}\ominus {\frak L}_{0} \oplus
\Lambda^{2}({\frak L}_{0})\ominus\cdots , \\
& S({\frak L}_{1})=\sum_{k=0}^{\infty}(-1)^{k}S^{k}({\frak L}_{1})
=\mathbb{C}\ominus {\frak L}_{1} \oplus S^{2}({\frak L}_{1})\ominus\cdots , \\ 
& H({\frak L})=\sum_{k=1}^{\infty}(-1)^{k+1}H_{k}({\frak L})
=H_1({\frak L})\ominus H_2({\frak L})\oplus H_3({\frak L})\ominus\cdots ,
\end{aligned}
\end{equation}
the alternating direct sum of superspaces.
Then it is easy to see that $C({\frak L})
=\Lambda({\frak L}_{0})\otimes S({\frak L}_{1})$.

\vskip 2mm
We recall the very basic theory of symmetric functions. 
Let $x_1,\cdots,x_k$ be indeterminates and for 
each $n\in \mathbb{N}$, define
\begin{equation}
\begin{aligned}
p_n(x_1,\cdots,x_k)&=\sum_{i=1}^{k} x_i^n, \\
h_n(x_1,\cdots,x_k)&=\sum_{i_1\leq\cdots\leq i_n}
                      x_{i_1}\cdots x_{i_n}, \\
e_n(x_1,\cdots,x_k)&=\sum_{i_1 < \cdots < i_n}
                      x_{i_1}\cdots x_{i_n},
\end{aligned} 
\end{equation}
which are called the {\it $n$th power sum}, 
the {\it $n$th completely symmetric function}, 
and the {\it $n$th elementary symmetric function}, respectively.
Note that the elementary symmetric function 
$e_n$ is defined only for $n\leq k$.
To compute the generalized supercharacters of 
$\Lambda({\frak L}_{0})$ and $S({\frak L}_{1})$,
we recall the following fundamental lemma in the theory of symmetric
functions:

\begin{lem} \label {Lem 2.2} {\rm (\cite{Mac79})}
\begin{equation}
\exp \left(\sum_{r\geq 1}\frac{p_r}{r}t^r \right)
=\sum_{n\geq 0}h_n t^n
=\frac{1}{\sum_{n\geq 0}(-1)^n e_n t^n}.
%=(\sum_{n\geq 0}(-1)^n e_n t^n)^{-1}.
\end{equation}
\qed     
\end{lem}

\vskip 2mm
Let $A$ be an $n\times n$ complex matrix with eigenvalues
$x_1,\cdots,x_n$. 
By Lemma \ref {Lem 2.2}, we have 
\begin{equation}
\begin{aligned}
\exp \left(\sum_{r\geq 1}\tr(A^r)\frac{t^r}{r} \right)
&=\sum_{m\geq 0}\tr(A|S^m(\mathbb{C}^n))t^m, \\
\exp \left(-\sum_{r\geq 1}\tr(A^r)\frac{t^r}{r} \right)
&=\sum_{m\geq 0}(-1)^m\,\tr(A|\Lambda^m(\mathbb{C}^n))t^m.
\end{aligned}
\end{equation}

Let $t=(t\ag)_{\ag \in \GG}$ be a sequence 
of nonnegative integers indexed by 
$\GG$ with only finitely many nonzero terms,
and set $|t|=\sum t\ag$. 
Since the $k$th exterior power $\Lambda^k ({\frak L}_{0})$ 
is decomposed as 
\begin{equation*}
\Lambda^k({\frak L}_{0})=
\bigoplus_{|t|=k}\left(\bigotimes_{\ag \in \GG_0}
\Lambda^{t\ag}({\frak L}_{\ag})\right),
\end{equation*} 
we have 
{\allowdisplaybreaks\begin{eqnarray*}
\ch_{g}\Lambda({\frak L}_{0})
&=&\sum_{k=0}^{\infty}(-1)^k \ch_{g}\Lambda^k({\frak L}_{0}) \\
&=&\prod_{\ag \in \GG_0}
\left(\sum_{m=0}^{\dim {\frak L}_{\ag}}(-1)^m
\tr(g|\Lambda^m({\frak L}_{\ag}))e^{m\ag} \right) \\
&=&\prod_{\ag \in \GG_0}
\left(\sum_{m\geq 0}(-1)^m\tr(g|\Lambda^m({\frak L}_{\ag}))e^{m\ag}\right) \\ 
&=&\prod_{\ag \in \GG_0}
\exp\left(-\sum_{r=1}^{\infty}\frac{1}{r}\tr(g^r|{\frak L}_{\ag})e^{r\ag}
\right)
\end{eqnarray*}}
for all $g\in G$. 
Similarly, we have 
{\allowdisplaybreaks\begin{eqnarray*}
\ch_{g}S({\frak L}_{1})&=&\sum_{k=0}^{\infty}(-1)^k \ch_{g}S^k({\frak L}_{1}) \\
&=&\prod_{\ag \in \GG_1}
\left(\sum_{m\geq 0}(-1)^m\tr(g|S^m({\frak L}_{\bh}))e^{m\bh}\right) \\ 
&=&\prod_{\ag \in \GG_1}
\left(\sum_{m\geq 0}\tr(g|S^m({\frak L}_{\bh}))(-e^{\bh})^m \right) \\
&=&\prod_{\ag \in \GG_1}
\exp \left(\sum_{r=1}^{\infty}\frac{(-1)^r}{r}
\tr(g^r|{\frak L}_{\bh})e^{r\bh} \right).
\end{eqnarray*}}
It follows that
\begin{equation}
\begin{aligned}
& \ch_{g}C({\frak L})=\sum_{k=0}^{\infty}(-1)^k \ch_{g}C_k({\frak L})
          =\ch_{g}\Lambda({\frak L}_{0})\cdot \ch_{g}S({\frak L}_{1}) \\
&=\prod_{\ag \in \GG_0}
\exp\left(-\sum_{r=1}^{\infty}\frac{1}{r}\tr(g^r|{\frak L}_{\ag})
e^{r\ag}\right)          \\
& \times \prod_{(\beta, b) \in \GG_1}
\exp\left(\sum_{r=1}^{\infty}\frac{(-1)^r}{r}
\tr(g^r|{\frak L}_{\bh})e^{r\bh} \right).
\end{aligned}
\notag
\end{equation}
Replacing $E^{\ag}=\psi(a)e^{\ag}$ and
$\text {str}(g|{\frak L}_{\ag})=\psi(a)\tr(g|{\frak L}_{\ag})$, 
the above identity yields
\begin{equation*}
\text {sch}_{g}C({\frak L})=
\prod_{\ag \in \GG}
\exp\left(-\sum_{r=1}^{\infty}\frac{1}{r}\text {str}(g^r|{\frak L}_{\ag})
E^{r\ag}\right).
\end{equation*}
Hence, by the Euler-Poincar\'{e} principle,
we obtain the {\it generalized denominator identity}
for the $\theta$-colored graded Lie superalgebra
${\frak L}=\bigoplus_{\ag \in \GG} {\frak L}_{\ag}$:
\begin{equation} \label {2.14}
\prod_{\ag \in \GG}
\exp\left(-\sum_{r=1}^{\infty}\frac{1}{r}\text {str}(g^r|{\frak L}_{\ag})
E^{r\ag}\right) =1-\text {sch}_{g}H({\frak L}).
\end{equation}

\vskip 2mm 
When $g=1$, the identity element of $G$, 
we get the \textit{denominator identity}
for the Lie superalgebra ${\frak L}$ (\cite{Ka97}):
\begin{equation*}
\prod_{\GG}(1-E^{\ag})^{\text{sdim}{\frak L}_{\ag}}
= 1-\text {sch} H({\frak L}).
\end{equation*}

\vskip 2mm
Let $P(\GG)=\{ (\alpha, a) \in \GG | \ \text {str}(g| H({\frak L})_{\ag})
\neq 0 \}$ and 
$\{\,(\alpha_i,a_j)\,|i,j=1,2,3,\cdots \}$
be an enumeration of $P(\GG)$.
For $(\tau,b)\in\GG$, set 
\begin{equation} \label {2.15}
T(\tau,b)=\{\,s=(s_{ij})\,|\,s_{ij}\geq 0, \ \ 
\sum_{i,j} s_{ij}(\alpha_i,a_j)=(\tau,b)\,\},
\end{equation}
the set of all partitions of $(\tau, b)$ into a sum of 
$(\alpha_i, a_j)$'s, and define 
\begin{equation} \label {2.16}
W_{g}(\tau,b)=\sum_{s\in T(\tau,b)}
        \frac{(|s|-1)!}{s!}\prod_{i,j} 
\text {str}(g|H({\frak L})_{(\alpha_i,a_j)})^{s_{ij}}.
\end{equation}
We will call $W_g(\tau, b)$ the {\it Witt partition function} for
$g\in G$. In the next theorem, using the generalized 
denominator identity (2.14), we derive a closed form formula for
the supertraces $\text {str}(g|{\frak L}_{\ag})$ $(\ag \in \GG)$ in terms of
the Witt partition functions.

\vskip 2mm
\begin{thm} \label {Thm 2.3}
For $g\in G$ and $\ag\in\GG$, we have
\begin{equation} \label {2.17}
\text {str}(g|{\frak L}_{\ag})=\sum \Sb d > 0 \\ \ag=d(\tau,b) \endSb
\frac{1}{d}\mu(d)W_{g^d}(\tau,b).
\end{equation}
\end{thm}

\noindent
{\it Proof.} \ 
By the generalized denominator idenity (\ref {2.14}), we have
\begin{equation*}
\begin{aligned}
\exp & \left(\sum_{ \ag \in \GG} \sum_{r=1}^{\infty}  
\frac{1}{r}\text {str}(g^r|{\frak L}_{\ag})E^{r\ag} \right) \\
&= \frac {1}
{1-\sum_{i,j=1}^{\infty}
\text {str}(g|H({\frak L})_{(\alpha_i,a_j)})E^{(\alpha_i,a_j)}}.
\end{aligned}
\end{equation*}
Taking the logarithm and using the formal power series 
$\log(1-x)=-\sum_{k=1}^{\infty}\frac{1}{k}x^k$, we obtain
{\allowdisplaybreaks\begin{eqnarray*}
\lefteqn{
\sum_{\ag \in \GG}\sum_{r=1}^{\infty} 
\frac{1}{r}\text {str}(g^r|{\frak L}_{\ag})E^{r\ag} } \\
& &=-\log
(1-\sum_{i,j=1}^{\infty}\text {str}(g|H({\frak L})_{(\alpha_i,a_j)})E^{(\alpha_i,a_j)})^{-1} \\
& &=\sum_{m=1}^{\infty}\frac{1}{m}
\left(\sum_{i,j=1}^{\infty}\text {str}(g|H({\frak L})_{(\alpha_i,a_j)})
E^{(\alpha_i,a_j)} \right)^m \\
& &=\sum_{m=1}^{\infty}\frac{1}{m}
\sum \Sb s=(s_{ij}) \\ \sum s_{ij}=m  \endSb
\frac{(\sum s_{ij})!}{\prod s_{ij}!}
\prod\text {str}(g|H({\frak L})_{(\alpha_i,a_j)})^{s_{ij}} 
E^{\sum s_{ij}(\alpha_i,a_j)} \\
& &=\sum_{(\tau, b) \in \GG}
    \left(\sum_{s \in T(\tau, b)}
         \frac{(\sum s_{ij}-1)!}{\prod s_{ij}!}
         \prod\text {str}(g|H({\frak L})_{(\alpha_i,a_j)})^{s_{ij}}\right)
         E^{(\tau, b)} \\
& &=\sum_{(\tau, b) \in \GG}W_{g}(\tau, b) E^{(\tau, b)}.
\end{eqnarray*}}
It follows that 
\begin{equation*}
W_{g}(\tau, b)=\sum \Sb k> 0 \\ (\tau, b) =k \ag \endSb
\frac{1}{k}\text {str}(g^k|{\frak L}_{\ag}).
\end{equation*}
Hence, by M\"{o}bius inversion, we obtain
\begin{equation*}
\text {str}(g|{\frak L}_{\ag})=\sum \Sb  d> 0   \\ d(\tau,b)=\ag \endSb 
\frac{1}{d}\mu(d)W_{g^d}(\tau,b).
\end{equation*}
\qed 

\vskip 2mm
\noindent
{\it Remark.} \  Note that our supertrace formula is a generalization 
of the {\it trace formula} for graded Lie algebras obtained in \cite{KK}.
Also, when $g=1$, our supertrace formula yields the 
{\it superdimension formula} 
for homogeneous subspaces of graded Lie superalgebras (cf. \cite{Ka97}).

\vskip 2cm 
\section{Homology of Lie Superalgebras}

The homology group of a Lie algebra is defined
as the torsion group of its universal enveloping algebra and 
characterized by Koszul's complex.
Similarly, we define the homology group of a graded Lie superalgebra
as the torsion group of its universal enveloping algebra
viewed as a supplmented algebra (cf. \cite{CE}, \cite{Fu}).
In this section, we will show that there is a super-analogue of
Koszul's complex for a graded Lie superalgebra and characterize
its homology group. Our argument follows the framework given in 
\cite {CE} and \cite{Fu}. 

\vskip 2mm 
Let ${\frak L}={\frak L}_0\oplus {\frak L}_1$ be a 
Lie superalgebra and $U=U({\frak L})$ be its universal enveloping algebra.
For each $k\geq 0$, define
\begin{equation} \label {3.1}
C_{k}=C_k({\frak L})=\bigoplus_{p+q=k}\Lambda^p({\frak L}_0)
\otimes S^q({\frak L}_1).
\end{equation}
Consider the following chain complex 
$(M_{k}, \pr^{\pm}_{k})$ 
$(k\geq -1)$, where 
$$M_{k} = \cases U \otimes_{\mathbb{C}} C_{k} ({\frak L}) \ \ &
\text {if} \ k\ge 0, \\
\mathbb{C} \ \ & \text {if} \ k=-1,
\endcases$$ 
and the differentials $\pr_k:M_k \rightarrow 
M_{k-1}$ are given by
\begin{equation} \label {3.2}
\begin{aligned}
& \pr_{k}(u\otimes(x_1\wed\cdots\wed x_p)\otimes(y_1\cdots y_q)) \\
&=\sum_{1\leq s <t \leq p}(-1)^{s+t}
  u\otimes([x_s,x_t]\wed x_1\wed\cdots\wed\widehat{x_s}\wed\cdots\wed 
\widehat{x_t}\wed\cdots\wed x_p)
 \otimes(y_1\cdots y_q) \\
& +\sum_{s=1}^{p}\sum_{t=1}^{q}(-1)^{s}
    u\otimes(x_1\wed\cdots\wed\widehat{x_s}\wed\cdots\wed x_p)
    \otimes([x_s,y_t]y_1\cdots\widehat{y_t}\cdots y_q)  \\
&-\sum_{1\leq s < t\leq q }
    u\otimes([y_s,y_t]\wed x_1\wed\cdots\wed x_p)
    \otimes(y_1\cdots\widehat{y_s}\cdots\widehat{y_t}\cdots y_q) \\
&+\sum_{s=1}^{p}(-1)^{s+1}
    (u\cdot x_s)\otimes(x_1\wed\cdots\wed\widehat{x_s}\wed\cdots\wed x_p)
    \otimes(y_1\cdots y_q)  \\
& + (-1)^p\sum_{t=1}^{q}(u\cdot y_t)
   \otimes(x_1\wed\cdots\wed x_p)\otimes(y_1\cdots\widehat{y_t}\cdots y_q)
\end{aligned}
\end{equation}
for $k\ge 1$, $\pr_{0}$ is an augmentation map 
extracting the constant term, and $\pr_{-1}=0$. 
Then it is easy to verify that $\pr_{k-1} \circ \pr_{k}=0$.

\vskip 2mm
\begin{prop} \label{Prop 3.1} \
The chain complex 
$M=(M_{k}, \pr_{k})$ is a free resolution of the
trivial 1-dimensional left $U$-module  $\mathbb{C}$.
\end{prop}

\noindent
{\it Proof.} \ 
{\bf Case 1.} \ Suppose first that ${\frak L}$ is abelian. 
Then, by the Poincar\'e-Birkhoff-Witt Theorem, we have 
\begin{equation}
U({\frak L})\simeq S({\frak L}_0)\otimes \Lambda ({\frak L}_1)
\end{equation}
as a $\mathbb{C}$-vector space,
and the differential maps are simplified to 
\begin{equation*}
\begin{aligned}
& \pr_k(u\otimes(x_1\wed\cdots\wed x_p)\otimes(y_1\cdots y_q) \\
&=\sum_{s=1}^{p}(-1)^{s+1}
    (u\cdot x_s)\otimes(x_1\wed\cdots\wed\widehat{x_s}\wed\cdots\wed x_p)
    \otimes(y_1\cdots y_q) \\
& + (-1)^p\sum_{t=1}^{q}(u\cdot y_t)\otimes(x_1\wed\cdots\wed x_p)
    \otimes(y_1\cdots\widehat{y_t}\cdots y_q).
\end{aligned}
\end{equation*}

\vskip 2mm 
We define the homotopy map $D: M_{k} \rightarrow M_{k+1}$ by
\begin{equation*}
\begin{aligned}
& D((u_{1} \cdots u_{s} v_{1} \cdots v_{t})
\otimes (x_{1} \wedge \cdots \wedge x_{p})
\otimes (y_{1} \cdots y_{q})) \\
&= \sum_{i=1}^s (u_{1} \cdots \widehat {u_{i}} \cdots u_{s} v_{1} \cdots v_{t})
\otimes (u_{i} \wedge x_{1} \wedge \cdots \wedge x_{p})
\otimes (y_{1} \cdots y_{q}) \\
&= \sum_{j=1}^t (-1)^{p+t-j} (u_{1} \cdots u_{s} 
v_{1} \cdots \widehat {v_{j}} \cdots v_{t})
\otimes (x_{1} \wedge \cdots \wedge x_{p})
\otimes (v_{j} y_{1} \cdots y_{q}),
\end{aligned}
\end{equation*}
where $u_{i} \in {\frak L}_{0}$ and $v_{j} \in {\frak L}_{1}$.
Then it is easy to verify that 
$$D \pr + \pr D = (p+q+ s+t) \, 1_{M},$$
which implies that our complex $M=(M_{k}, \pr_{k})$ is exact
in this case. 

\vskip 2mm
\noindent
{\bf Case 2.} \ 
Let ${\frak L}$ be any Lie superalgebra. Recall the construction of
the universal enveloping algebra of ${\frak L}$.
Put $U^n({\frak L})=\pi(\bigoplus_{k=0}^{n} {\frak L}^{\otimes k})$, 
where $\pi:\mathcal{T}({\frak L})\arw U({\frak L})$ 
is the canonical projection map and 
$X^{(p)}_n=U^{p-n}({\frak L})\otimes C_n({\frak L})$ for $0\leq n \leq p$.

\vskip 2mm 
Consider the complex
\begin{equation} \label {3.8}
0\larw X^{(p)}_p\stackrel{\pr_p}{\larw}\cdots\stackrel{\pr_2}{\larw}X^{(p)}_1
\stackrel{\pr_1}{\larw}X^{(p)}_0\stackrel{\epsilon}{\larw}\CN\larw 0.
\end{equation}
We will use the induction on $p$ to prove the exactness of (\ref{3.8}), 
which would imply the exactness of our original sequence.

\vskip 2mm 
Note that the exactness of (\ref{3.8}) is equivalent to that of the 
following sequence:
\begin{equation} \label{3.9}
0\larw X^{(p)}_p\stackrel{\pr_p}{\larw}\cdots\stackrel{\pr_2}{\larw}X^{(p)}_1
\stackrel{\pr_1}{\larw}X^{(p)}_0/\CN\larw 0.
\end{equation}
For $p=0,1$, our assertion is clear.
Suppose (\ref{3.9}) is exact for $p-1$.
Consider the diagram 
\begin{equation*} 
\begin{array}{ccccccccc}
 &0 & & & &0 & &0 &   \\
 &\daw & & & &\daw & &\daw &   \\
0\arw&0 &\arw &\cdots &\stackrel{\pr}{\arw} &X^{(p-1)}_1 
&\stackrel{\pr}{\arw} &X^{(p-1)}_0/\CN &\arw 0   \\
 &\daw & & & &\daw & &\daw &   \\
0\arw&X^{(p)}_p & \stackrel{\pr} \arw &\cdots &\stackrel{\pr}{\arw} &X^{(p)}_1
&\stackrel{\pr}{\arw} &X^{(p)}_0/\CN &\arw 0   \\
 &\daw & & & &\daw & &\daw &   \\
0\arw&X^{(p)}_p &\stackrel{\overline {\pr}}\arw &\cdots &
\stackrel{\overline{\pr}}{\arw} &X^{(p)}_1/X^{(p-1)}_1
&\stackrel{\overline{\pr}}{\arw} &X^{(p)}_0/X^{(p-1)}_0&\arw 0.   \\
 &\daw & & & &\daw & &\daw &   \\
 &0 & & & &0 & &0 &  
\end{array}
\end{equation*}

\vskip 2mm 
Note that 
\begin{equation*}
X^{(p)}_n/X^{(p-1)}_n \simeq \left( U^{p-n}({\frak L})/ 
U^{p-n-1}({\frak L})\right)\otimes C_n({\frak L}) \ \ (0 \le n \le p-1)
\end{equation*}
as a $\CN$-vector space, and that, 
by the  Poincar\'{e}-Birkhoff-Witt Theorem,
\begin{equation*}
 U^{p-n}({\frak L})/U^{p-n-1}({\frak L})
\simeq\bigoplus_{k+l=p-n}S^k({\frak L}_0)\otimes \Lambda^l({\frak L}_1)
\end{equation*}
by the canonical homomorphism and $X^{(p)}_p\simeq\CN\otimes C_p({\frak L})$.
Finally, the differentials in bottom row are given by 
\begin{equation} \label {3.10}
\begin{aligned}
& \overline{\pr_n}
(u\otimes(x_1\wed\cdots\wed x_p)\otimes(y_1\cdots y_q)+ X^{(p-1)}_n) \\
& =\sum_{s=1}^{p}(-1)^{s+1}
    (u\cdot x_s)\otimes
    (x_1\wed\cdots\wed\widehat{x_s}\wed\cdots\wed x_p)
    \otimes(y_1\cdots y_q) \\
& + (-1)^p\sum_{t=1}^{q}
   (u\cdot y_t)\otimes(x_1\wed\cdots\wed x_p)
   \otimes(y_1\cdots\widehat{y_t}\cdots y_q) + X^{(p-1)}_{n-1}. 
\end{aligned}
\end{equation}
Thus the bottom row can be viewed as a subcomplex of 
the complex considered in Case 1, 
and hence it is exact.
Therefore, the middle row is exact, which completes our induction argument.

\vskip 2mm 
Since $C({\frak L})$ is a $\CN$-vector space, 
$U\otimes C({\frak L})$ is the direct sum of free $U$-modules.
Hence, $U\otimes C({\frak L})$ is a free 
resolution of the 1-dimensional  trivial left $U$-module $\CN$. 
\qed

\vskip 3mm 
Let $V$ be a right $U$-module. The homology module $H_{k}({\frak L},V)$
is defined to be
\begin{equation} \label {3.11}
H_{k}({\frak L},V) = \text{Tor}^{U}_{k}(V, \mathbb {C})
=H_{k}(V \otimes_{U} M),
\end{equation}
where $M=(M_k, \pr_k)$ is the free resolution
of the 1-dimensional trivial left $U$-module $\mathbb {C}$ given by 
Proposition \ref{Prop 3.1}.
Since we have 
$$V \otimes_{U} M \cong 
V \otimes_{U} U \otimes_{\mathbb{C}} C_{k}({\frak L})
\cong V \otimes_{\mathbb {C}} C_{k}({\frak L}) \ \ (k\ge 0),$$
and $$V \otimes_{U} \mathbb{C} \cong V,$$
as $\mathbb{C}$-vector spaces, 
we obtain the following {\it standard complex} for the
homology of Lie superalgebras:
\begin{equation} \label {3.12}
\cdots\larw V\otimes_{\mathbb {C}} C_k({\frak L})
\stackrel{d_k}{\larw}V \otimes_{\mathbb {C}}C_{k-1}({\frak L})
\stackrel{d_{k-1}}{\larw}
\cdots\stackrel{d_1}{\larw} V\otimes_{\mathbb {C}}C_{0}({\frak L})
\cong V \larw 0,
\end{equation}
where the differentials 
$d_{k}=1_{V} \otimes_{U} \pr_{k}:
V \otimes_{\mathbb {C}} C_{k} ({\frak L})
\longrightarrow V \otimes_{\mathbb {C}} C_{k-1} ({\frak L})$ are
given by
\begin{equation} \label {3.13}
\begin{aligned}
& d_{k}(v\otimes(x_1\wed\cdots\wed x_p)\otimes(y_1\cdots y_q)) \\
&=\sum_{1\leq s <t \leq p}(-1)^{s+t}
  v\otimes([x_s,x_t]\wed x_1\wed\cdots\wed\widehat{x_s}\wed\cdots\wed 
\widehat{x_t}\wed\cdots\wed x_p)
 \otimes(y_1\cdots y_q) \\
& +\sum_{s=1}^{p}\sum_{t=1}^{q}(-1)^{s}
    v\otimes(x_1\wed\cdots\wed\widehat{x_s}\wed\cdots\wed x_p)
    \otimes([x_s,y_t]y_1\cdots\widehat{y_t}\cdots y_q)  \\
&-\sum_{1\leq s < t\leq q }
    v\otimes([y_s,y_t]\wed x_1\wed\cdots\wed x_p)
    \otimes(y_1\cdots\widehat{y_s}\cdots\widehat{y_t}\cdots y_q) \\
&+\sum_{s=1}^{p}(-1)^{s+1}
    (v\cdot x_s)\otimes(x_1\wed\cdots\wed\widehat{x_s}\wed\cdots\wed x_p)
    \otimes(y_1\cdots y_q)  \\
& +(-1)^p\sum_{t=1}^{q}(v\cdot y_t)
   \otimes(x_1\wed\cdots\wed x_p)\otimes(y_1\cdots\widehat{y_t}\cdots y_q)
\end{aligned}
\end{equation}
for $v \in V$, $x_{i} \in {\frak L}_{0}$, $y_{j} \in {\frak L}_{1}$. 
In particular, when $V=\mathbb {C}$ is the trivial 1-dimensional 
right $U$-module, the standard complex (\ref{3.12}) is reduced 
to the complex (\ref{2.8}) given in Section 2. 

\vskip 2mm
As a corollary of Proposition \ref{Prop 3.1}, we obtain:

\vskip 2mm 
\begin{cor} \ 
Let $\A$ be a countable abelian group with a coloring map $\theta$
and $V=\oplus_{a\in\A}V_a$ be a $\theta$-colored superspace.
Let ${\frak L}$ be the free Lie superalgebra generated by $V$. 
Then we have 
\begin{equation}
H_{k}({\frak L}) = \cases V  \ \ & \text {if} \ k=1, \\
0 \ \ & \text {if} \ k\ge 2.
\endcases
\end{equation}
\end{cor}

\noindent
{\it Proof.} \ Note that the universal enveloping algebra of ${\frak L}$ is 
isomorphic to the tensor algebra $\mathcal{T}(V)$ generated by $V$.
Consider the following free resolution of the 1-dimensional 
trivial left $\mathcal{T}(V)$-module  $\CN$:
\[0\arw \mathcal{T}(V)\otimes V\stackrel{\varphi}{\arw}
\mathcal{T}(V)\stackrel{\psi}{\arw} 
\CN\simeq \mathcal{T}(V)/\mathcal{T}(V)\otimes V \arw 0,\]
where $\varphi$ is the natural inclusion map(=multiplication map) 
and $\psi$ is the augmentation map.

\vskip 2mm 
By tensoring with trivial right $\mathcal{T}(V)$-module $\CN$, we get
\[0\arw \CN\otimes_{\mathcal{T}(V)}\mathcal{T}(V)\otimes V
\stackrel{1\otimes\varphi}{\arw}\CN\otimes_{\mathcal{T}(V)}\mathcal{T}(V)
\stackrel{1\otimes\psi}{\arw}\CN\otimes\CN\arw 0. \]

Since $1\otimes\varphi$ is a zero map, 
$H_1(\mathcal{L})=\text{Ker}\,(1\otimes\varphi)=
\CN\otimes_{\mathcal{T}(V)}\mathcal{T}(V)\otimes V\simeq V$, 
and $H_k(\mathcal{L})=0$ for $k\geq 2$.
\qed

\vskip 2cm

\section{Free Lie Algebras and Free Lie Superalgebras}

%Free Lie algebras and free Lie superalgebras

Let $\A$ be a countable abelian group with a coloring map $\theta$
and $\G$ be a countable abelian semigroup satisfying the finiteness
condition given in Section 2. 
Let $V=\bigoplus_{\ag \in \GG} V_{\ag}$ be a $(\GG)$-graded 
$\theta$-colored superspace with finite dimensional homogeneous
subspaces, and let ${\frak L}=\bigoplus_{\ag \in \GG} 
{\frak L}_{\ag}$ be the free Lie superalgebra generated by $V$.
Suppose that a group $G$ acts on $V$ preserving
the $(\GG)$-gradation. Then the action of $G$ on $V$ can be
extended to ${\frak L}$ by Lie superalgebra automorphisms.

\vskip 2mm
As we have seen in Section 3, we have 
$H_{1}({\frak L})=V$ and $H_{k}({\frak L})=0$ for $k\ge 2$.
Hence, for any $g\in G$, the generalized denominator identity for
the free Lie superalgebra ${\frak L}$ is equal to
\begin{equation} \label {4.1}
\begin{aligned}
 \prod_{\ag \in \GG} & \exp\left( - \sum_{k=1}^{\infty} \frac{1}{k}
\text {str} (g^k | {\frak L}_{\ag}) E^{k\ag} \right)\\
&=1-\text {sch}_{g} V =1-\sum_{\ag \in \GG} 
\text {str} (g| V_{\ag}) E^{\ag}.
\end{aligned}
\end{equation}
For any $g\in G$, let $P(V, \GG)=\{ \ag \in \GG \ | \ 
\text {str}(g|V_{\ag})\neq 0\}$ and let $\{ (\alpha_i, a_{j}) \ 
|\ i,j=1,2,3,\cdots \}$ be an enumeration of the set 
$P(V, \GG)$. For each $(\tau, b) \in \GG$, we denote by $T(\tau, b)$
the set of all paritions of $(\tau, b)$ into a sum of $(\alpha_i, a_{j})$'s
as defined in (\ref{2.15}), and let $W_{g}(\tau, b)$ be the Witt partition 
function as defined in (\ref{2.16}). 
Then our supertrace formula (\ref{2.17}) yields 
\begin{equation} \label {4.2}
\text {str}(g|{\frak L}_{\ag})= \sum \Sb d>0 \\ \ag =d(\tau, b) \endSb
\frac{1}{d} \mu(d)  \sum_{s\in T(\tau, b)} \frac
{(|s|-1)!}{s!} \prod_{i,j} 
\text {str} (g^d|V_{(\alpha_i, a_j)})^{s_{ij}}.
\end{equation}

\vskip 2mm
\noindent
{\it Remark.} \ The supertrace formula (\ref{4.2}) 
is a generalization of the
superdimension formula for the free Lie superalgebras given in 
\cite{Ka96} and \cite{Ka97}.

\vskip 2mm
\noindent
\begin{ex} \rm {In this example, we will consider the simplest application
of our supertrace formula to free Lie superalgebras. 
Consider the superspace $V=V_{0} \oplus V_{1}$ with 
$\dim V_{0}=r$ and $\dim V_{1} =s$ for some $r,s \in \mathbb {Z}_{>0}$,
and let ${\frak L}$ be the free Lie superalgebra generated by $V$.
By setting $\text {deg} v=1$ for all $v\in V$, the free Lie superalgebra 
${\frak L}$ has a $\mathbb {Z}_{>0}$-gradation ${\frak L} 
=\bigoplus_{n=1}^{\infty} {\frak L}_{n}$ induced by $V$. 
Let $G=GL(V_{0}) \times GL(V_{1}) \cong GL(r) \times GL(s)$.
For $(g,h)\in G$, we denote by $t_g=\tr(g|V_{0})$ and $t_{h}=\tr(h|V_{1})$.
Since $P(V, \mathbb {Z}_{>0})=\{1\}$ and $\text {str}((g,h) |V)
=\tr(g|V_0)-\tr(h|V_{1}) = t_{g}-t_{h}$ for all $(g,h) \in G$, 
the Witt partition function $W_{(g,h)}(n)$ is given by
$$W_{(g,h)}(n)=\frac{(n-1)!}{n!} (t_{g} - t_{h})^{n}
=\frac{1}{n} (t_g -t_h)^n.$$
Hence, by the supertrace formula (\ref{4.2}), we have
\begin{equation} \label {4.3}
\text {str}((g,h)|{\frak L}_{n})=\frac{1}{n} \sum_{d|n} \mu(d)
(t_{g}-t_{h})^{\frac{n}{d}}. \hskip 1cm \qed
\end{equation}}
\end{ex}

\vskip 2mm
Next, we will discuss the relation of the free Lie algebra
and the free Lie superalgebra generated by the same vector space
$V=\bigoplus_{\ag \in \GG} V_{\ag}$ on which a group $G$ acts
preserving the $(\GG)$-gradation.
By neglecting the coloring map $\theta$ on $\A$,
consider $V$ as a $(\GG)$-graded vector space, and let
$L=\oplus_{\ag\in\GG}L_{\ag}$ be the free Lie algebra generated by $V$.
Then the action of $G$ on $V$ can be extended to $L$ by Lie algebra
automorphisms. 
Hence, for any $g\in G$, the generalized ednominator identity for $L$
is equal to 
\begin{equation} \label {4.4}
\begin{aligned}
 \prod_{\ag\in\GG} & \exp \left(-\sum_{k=1}^{\infty}\frac{1}{k}
\tr(g^k|L_{\ag})e^{k\ag} \right) \\
& =1-\ch_{g}V =1-\sum_{\ag \in \GG} \tr(g|V_{\ag}) e^{\ag}.
\end{aligned}
\end{equation} 

Since $E^{\ag}=\psi(a) e^{\ag}$ and $\text {str}(g|V_{\ag})
=\psi(a) \tr(g|V_{\ag})$ for all $\ag \in \GG$, 
the right-hand sides of (\ref{4.1}) and (\ref{4.4}) are the same.
Moreover, note that, for an $n\times n$ complex matrix $A$ with
eigenvalues $x_{1}, \cdots, x_{n}$, Lemma 2.1 implies
\begin{equation}
\begin{aligned}
& \exp\left( \sum_{k=1}^{\infty} \frac{1}{k} \tr(A^k) t^k \right)
 \exp\left( \sum_{k=1}^{\infty} \frac{1}{k} (-1)^k \tr(A^k) t^k \right)\\
&=\prod \frac{1}{1-x_{i} t} \prod \frac{1}{1+x_{i} t} 
=\prod \frac{1}{1-x_{i}^2 t^2} \\
&=\sum_{n=0}^{\infty} h_{n}(x_{1}^2, \cdots, x_{n}^2) t^{2n} 
=\exp \left( \sum_{k=1}^{\infty} \frac{1}{k} \tr(A^{2k}) t^{2k}\right).
\end{aligned}
\notag
\end{equation}
Therefore, the left-hand side of (\ref{4.4}) is equal to
\begin{equation} 
\begin{aligned}
 \prod_{\ag \in \GG_0}& \exp \left(-\sum_{k=1}^{\infty}\frac{1}{k}
\tr(g^k|L_{\ag})e^{k\ag} \right) 
\prod_{(\beta, b) \in \GG_1} \exp \left(-\sum_{k=1}^{\infty}\frac{1}{k}
\tr(g^k| L_{(\beta, b)})e^{k(\beta, b)} \right) \\
=& \prod_{\ag\in\GG_0} \exp \left(-\sum_{k=1}^{\infty}\frac{1}{k}
\tr(g^k|L_{\ag})E^{k\ag} \right)\\
& \times \prod_{(\beta, b) \in\GG_1} \exp \left(-\sum_{k=1}^{\infty}
\frac{(-1)^k}{k}\tr(g^k|L_{(\beta, b)})E^{k(\beta, b)} \right) \\
=&\prod_{\ag\in\GG_0} \exp \left(-\sum_{k=1}^{\infty}\frac{1}{k}
\tr(g^k|L_{\ag})E^{k\ag} \right) \\
& \times \prod_{ (\beta, b) \in\GG_1}\exp \left(\sum_{k=1}^{\infty}\frac{1}{k}
\tr(g^{k}|L_{(\beta, b)})E^{k(\beta, b)} \right) \\
& \times \prod_{ (\beta, b) \in\GG_1} \exp \left(-\sum_{k=1}^{\infty}
\frac{1}{k}\tr(g^{2k}|L_{(\beta, b)})E^{2k(\beta, b)}\right) \\
\end{aligned}
\notag
\end{equation}
\begin{equation}
\begin{aligned}
=&\prod_{\ag\in\GG} \exp \left(-\sum_{k=1}^{\infty}\frac{1}{k}
\psi(a) \tr(g^k|L_{\ag})E^{k\ag} \right) \\
& \times \prod \Sb (\alpha, a)=2 (\beta, b) \\ 
(\beta, b) \in \GG_1 \endSb 
\exp \left(-\sum_{k=1}^{\infty}\frac{1}{k} 
\tr(g^{2k}| L_{(\beta, b)})E^{k\ag} \right).
\end{aligned}
\notag
\end{equation}

Hence, for any $g\in G$, the generalized denominator identity 
for the free Lie superalgebra ${\frak L}$ is the same as
\begin{equation} \label {4.5}
\begin{aligned}
 \prod_{\ag \in \GG} & \exp \left( -\sum_{k=1}^{\infty} \frac{1}{k}
\text {str} (g^k| {\frak L}_{\ag}) E^{k\ag} \right)\\
=& 1-\sum_{\ag \in \GG} \text {str} (g|V_{\ag}) E^{\ag} \\
=& \prod_{\ag\in \GG} \exp \left(-\sum_{k=1}^{\infty}\frac{1}{k}
\psi(a) \tr(g^k|L_{\ag})E^{k\ag} \right) \\
& \times \prod \Sb (\alpha, a)=2 (\beta, b) \\ 
(\beta, b) \in \GG_1 \endSb 
\exp \left(-\sum_{k=1}^{\infty}\frac{1}{k} 
\tr(g^{2k}| L_{(\beta, b)})E^{k\ag} \right),
\end{aligned}
\end{equation}
which yields:

\begin{prop} \label {Prop 4.2}
Let $\A$ be a countable abelian group with a coloring map 
$\theta$ and $\G$ be a countable abelian semigroup satisfying
the finiteness condition given in Section 2. 
For a $(\GG)$-graded superspace $V=\bigoplus_{\ag \in \GG} V_{\ag}$
with finite dimensional homogeneous subspaces, let 
$L=\bigoplus_{\ag \in \GG} L_{\ag}$ be the free Lie algebra 
generated by $V$, and let 
${\frak L}=\bigoplus_{\ag \in \GG} {\frak L}_{\ag}$ be the
free Lie superalgebra generated by $V$. 
Suppose that a group $G$ acts on $V$ preserving the $(\GG)$-gradation.
Then, for any $g\in G$, we have
\begin{equation} \label {4.6}
\text {str}(g|{\frak L}_{\ag})=\psi(a) \tr(g|L_{\ag})
+\sum \Sb (\alpha, a)=2(\beta, b) \\ (\beta, b) \in \GG_1 \endSb 
\tr(g^2 | L_{(\beta, b)}). \
\end{equation}
\qed
\end{prop}

\noindent
{\it Remark.} \ When $g=1$, we recover the superdimension formula 
\begin{equation} \label {4.7}
\text {sdim} \, {\frak L}_{\ag} =\psi(a) \dim L_{\ag}
+\sum \Sb (\alpha, a)=2(\beta, b) \\ (\beta, b) \in \GG_1 \endSb 
\dim L_{(\beta, b)},
\end{equation}
which was obtained in  \cite {Ka97}.

\vskip 2cm

\section{Decomposition of Free Lie Superalgebras}

In this section, we investigate the structure of the free Lie
superalgebra generated by the natural representation of the
{\it general linear Lie superalgebra} $gl (k,l)$. 
For  $k,l\in\mathbb{Z}_{\geq 0}$, let  ${\frak L}=gl(k,l)$ 
be the space of all $(k+l)\times(k+l)$ matrices, and set 
\begin{equation} \label{5.1}
{\frak L}_0=\left\{ \left. 
\left(\begin{matrix}
       A & 0 \\  0 & D \end{matrix}
\right)
\right| 
\text{ $A$ is a $k\times k$ matrix and $D$ is an $l\times l$ matrix }
\right\},
\end{equation}     
\begin{equation} \label{5.2}
{\frak L}_1=\left\{\left.
\left(\begin{matrix}  0 & B \\  C & 0  \end{matrix}
\right)
\right| 
\text{ $B$ is a $k\times l$ matrix and $C$ is an $l\times k$ matrix }
\right\}.
\end{equation}
Then ${\frak L}={\frak L}_0\oplus {\frak L}_1$ and 
${\frak L}$ becomes a $\mathbb{Z}_2$-graded Lie superalgebra
called the \textit{general linear Lie superalgebra} with the superbracket
defined by
\begin{equation}
[X,Y]=XY-(-1)^{\alpha\beta}YX
\ \ \ \ \text {for} \ \ 
X\in {\frak L}_{\alpha},\,\,Y\in {\frak L}_{\beta},\,\,
\alpha,\beta\in\mathbb{Z}_2.
\end{equation}

\vskip 2mm

A $\mathbb{Z}_2$-graded superspace $V=V_0\oplus V_1$
is an {\it ${\frak L}$-module} if there is a
$\mathbb{Z}_2$-graded Lie superalgebra homomorphism
\begin{equation}
\phi:{\frak L} \longrightarrow gl(V) \cong gl(k',l'), 
\ \ \text{where }k'=\dim V_0,\,l'=\dim V_1.
\end{equation}
For ${\frak L}$-modules $V^1,\cdots,V^n$, the tensor product 
$V^1\otimes\cdots\otimes V^n$ becomes an ${\frak L}$-module 
with the action of ${\frak L}$ defined by 
\begin{equation*}
x\cdot (v_1\otimes\cdots\otimes v_n)
=\sum_{i=1}^{n}(-1)^{\deg x(\sum_{j < i}\deg v_j)}
v_1\otimes\cdots\otimes x \cdot v_i\otimes\cdots\otimes v_n
\end{equation*}
for $x\in {\frak L}$, $v_{i} \in V^{i}$. 

\vskip 2mm 

We recall some of the basic theory of symmetric functions. 
A \textit{partition} is any (finite or infinite) sequence
$\ld=(\ld_1,\ld_2,\cdots)$ of non-negative integers 
in decreasing orders $\ld_1\geq\ld_2\geq\cdots$ and containing
only finitely many non-zero terms.
The non-zero $\ld_i$'s are called the \textit{parts} of $\ld$.
The number of parts is the \textit{length} of $\ld$, denoted by $l(\ld)$, 
and the sum of parts is the \textit{weight} of $\ld$, denoted by $|\ld|$.
If $|\ld|=n$, then we say that $\ld$ is a {\it partition of $n$},
denote by $\ld\vdash n$.
We also use the notation which indicates the number of times
each integer occurs as a part. That is,  
$\ld=(1^{m_1},2^{m_2},\cdots )$ 
means that $i$ appears as a part of $\ld$ exactly $m_i$ times.
Each partition $\ld$ uniquely determines a \textit{Young diagram}.
The \textit{conjugate} of a partition $\ld$ is the partition ${\ld}'$
whose diagram is the transpose of that of  $\ld$
(i.e., the diagram obtained by reflection in the main diagonal). 
There is a partial ordering on the set of parititions of $n$,
called the \textit{natural (partial) ordering}
or the {\it dominance relation}, 
which is defined as follows:
\begin{equation*}
\ld \ge \mu \hspace{3mm}\text{if and only if}\hspace{3mm}
\sum_{i=1}^{k}\ld_i \geq \sum_{i=1}^{k}\mu_i\hspace{3mm}
\text{for all $k \ge 1$}.
\end{equation*}

\vskip 2mm

Let $x_1,\cdots,x_k$ be the indeterminates and $\ld$ a partition of $n$
with $l(\ld)\leq k$.
A \textit{Schur function} corresponding to $\ld$ is a symmetric function
in $x_1,\cdots,x_k$ defined by
\begin{equation}
S_{\ld}(x)=S_{\ld}(x_1\cdots,x_k)=
\frac{\left|x_{i}^{\ld_j+n-j}\right|}{\left|x_{i}^{n-j}\right|}.
\end{equation}
For a partition $\m=(\m_1,\m_2,\cdots)$, we define 
\begin{equation}
p_{\m}(x)=p_{\m}(x_1,\cdots,x_k)=\prod p_{\m_i}(x_1,\cdots,x_k),
\end{equation}
where $p_{\m_i}(x_1,\cdots,x_k)=\sum x_{j}^{\m_i}$ is the 
$\mu_{i}$-th power sum.
It is known that the set of Schur  functions
$\{\,S_{\ld}(x)\,|\,l(\ld)\leq k\,\}$ (resp. $\{\,p_{\ld}(x) |
\ l(\ld) \le k\}$)
forms a $\mathbb{Z}$-basis (resp.  $\mathbb{Q}$-basis)
of the ring of symmetric functions in $x_1,\cdots,x_k$,
and $\{\,S_{\ld}(x)\,|\,l(\ld)\leq k, \lambda \vdash n\,\}$ 
(resp. $\{\,p_{\ld}(x) | \ l(\ld) \le k, \ld \vdash n\}$)
forms a $\mathbb{Z}$-basis (resp.  $\mathbb{Q}$-basis)
of the subspace of symmetric functions in $x_1,\cdots,x_k$ with 
homogeneous degree $n$. 
The functions $S_{\ld}(x)$ and $p_{\ld}(x)$ are related as follows:

\vskip 2mm 

\begin{lem} \label {Lem 5.1} {\rm (\cite{Mac79})} \ 
For each partition $\rho$ of $n$, we have 
\begin{equation}
p_{\rho}(x)=\sum_{\ld\vdash n}\chi_{\ld}^{\rho}S_{\ld}(x),
\end{equation}
where $\chi_{\ld}^{\rho}$ is the character value of the irreducible 
representation of $S_n$ corresponding to $\ld$ at the conjugacy class
of type $\rho$.
\qed
\end{lem}

\vskip 2mm

A partition $\ld$ is called a \textit{$(k,l)$-hook partition} if
$\ld_{k+1}\leq l$.
Set $H(k,l;n)=\{\,\ld\vdash n\,|\,\ld_{k+1}\leq l\,\}$.
For convenience, we will often write $H$ for $H(k,l;n)$. 
For $\ld\in H(k,l;n)$, a \textit{$(k,l)$-hook Schur function} in 
$x_1,\cdots,x_k,y_1,\cdots,y_l$ corresponding to $\ld$ is a function
defined by
\begin{equation}
HS_{\ld}(x,y)=HS_{\ld}(x_1,\cdots,x_k,y_1,\cdots,y_l)= 
\sum_{\m < \ld}S_{\m}(x)S_{\ld'/\m'}(y),
\end{equation}
where $S_{\ld'/\mu'}(y)$ is the skew Schur function 
corresponding to $\ld',\mu'$ (cf. \cite{BR}, \cite{Mac79}, \cite{R}).
Equivalently,
\begin{equation}
HS_{\ld}(x,y)=\sum_{\m < \ld}S_{\m}(x)\left(\sum_{\n}N^{\ld'}_{\m'\n}S_{\n}(y)\right),
\end{equation}
where $N^{\ld'}_{\m'\n}$ is the {\it Littlewood-Richardson coefficient}
corresponding to $\ld',\m',\n$.
The combinatorial interpretation
of the Littlewood-Richardson coefficients is given as follows:
for given partitions $\ld,\mu,\nu$ satisfying $|\ld|=|\mu|+|\nu|$
and $\ld_i \geq \mu_i$ for all $i$,
$N^{\ld}_{\mu \nu}$ is the number of ways the Young diagram
for $\mu$ can be expanded to the Young diagram for $\ld$ 
by a {\it strict $\nu$-expansion} (if the above conditions do not 
hold, $N^{\ld}_{\mu \nu}=0$). 
For a partition 
$\nu=(\nu_1,\cdots,\nu_k)$, a \textit{$\nu$-expansion} of a Young diagram
for $\mu$ is obtained by first adding $\nu_1$ boxes 
to the Young diagram for $\mu$ with no two boxes in the same
column and putting the integer 1 in each of these $\nu_1$ boxes;
then adding similarly $\nu_2$ boxes with a 2, continuing until
finally $\nu_k$ boxes are added with the integer $k$. 
The expansion is called \textit{strict} if, when the integers in the 
boxes are listed from right to left, starting with the top row 
and working down, and one looks at the first $t$ entries in this list
(for any $t$ between 1 and $|\nu|$), each integer $p$
between 1 and $k-1$ occurs at least as many times as the next integer
$p+1$.

\vskip 2mm 

In \cite{BR}, Berele and Regev showed that, for a $\mathbb{Z}_2$-graded
superspace $V=V_0\oplus V_1$ viewed as a $gl(k,l)$-module
($k=\dim V_0,\,l=\dim V_1$), $V^{\otimes n}$ is completely
reducible as a $gl(k,l)$-module and its irreducible components are
parametrized by  ($k,l$)-hook partitions $\ld$, denoted by $V_{\ld}$.  

\vskip 2mm 

Let ${\frak L}(V)$ be the free Lie superalgebra generated by $V$.
Set ${\frak L}(V)_{(n,0)}={\frak L}(V)\cap (V^{\otimes n})_0$ ,
    ${\frak L}(V)_{(n,1)}={\frak L}(V)\cap (V^{\otimes n})_1$ and
        ${\frak L}_{n}(V)={\frak L}(V)_{(n,0)}\oplus 
{\frak L}(V)_{(n,1)}$.
Thus ${\frak L}(V)$ has an ($\mathbb{N}\times\mathbb{Z}_2$)-gradation:
${\frak L}(V)=\bigoplus_{\mathbb{N}\times\mathbb{Z}_2}
{\frak L}(V)_{(n,\alpha)}$.
The homogeneous subspace ${\frak L}_{n}(V)$ is 
a $gl(k,l)$-submodule of $V^{\otimes n}$, where 
the action of $x\in gl(k,l)$ induced from $V^{\otimes n}$ is given by 
\begin{equation*}
\begin{aligned}
& x\cdot[v_1,[v_2,[\cdots[v_{n-1},v_n]\cdots]]] \\
&=  \sum_{i=1}^{n}(-1)^{\deg x(\sum_{j < i}\deg v_j)}
   [v_1,[\cdots[xv_i[\cdots[v_{n-1},v_n]\cdots]]]],
\end{aligned}
\end{equation*}
and it is completely reducible 
with its irreducible components parametrized by $\ld's \in H$.

\vskip 2mm 

>From now on, we simply write 
${\frak L}_n$, ${\frak L}_{(n,0)}$ and ${\frak L}_{(n,1)}$ 
for ${\frak L}_{n}(V)$, ${\frak L}(V)_{(n,0)}$ and 
${\frak L}(V)_{(n,1)}$, respectively.
For each $\ld\in H(k,l;n)$, let $c_{\ld}$ be the multiplicity of
$V_{\ld}$ in ${\frak L}_n$.
To compute $c_{\ld}$, we would like to apply our supertrace formula.
By exponentiating the action of the Lie algebra 
$gl(k,l)_0=gl(k)\times gl(l)$ on $V$,
the group $GL(k) \times GL(l)$ acts on $V$, and hence 
$V^{\otimes n}$ and ${\frak L}_{n}(V)$ become
$GL(k)\times GL(l)$-modules. We first recall:

\vskip 2mm

\begin{prop} \label {Prop 5.2} {\rm (\cite{BR})} \
Let $\ld\in H(k,l;n)$ and $V_{\ld}$ be the corresponding 
irreducible representation of $gl(k,l)$.
Suppose that $g=(g_0,g_1)\in GL(k)\times GL(l)$ 
has eigenvalues $x_1,\cdots,x_k$ for $g_{0}$ and 
$y_1,\cdots,y_l$ for $g_1$.
Then we have 
\begin{equation}
\tr(g|V_{\ld})=HS_{\ld}(x,y).
\end{equation}
\qed
\end{prop}

\vskip 2mm 

Our supertrace formula (\ref {4.2}) for free Lie superalgebras yields
\begin{equation} \label {5.12}
\text {str}(g|{\frak L}_{(n,\alpha)})=
\sum_{d|(n,\alpha)}\frac{1}{d}\m(d)
\sum_{s\in T((n,\alpha)/d)}\frac{(|s|-1)!}{s!}
     \text {str}(g_0^d|V_0)^{s_0}\text {str}(g_1^d|V_1)^{s_1},
\end{equation}
where
{\allowdisplaybreaks\begin{eqnarray*}
T(m,0)&=&\{\,(s_0,s_1)\,|\,s_0+s_1=m,\ s_1\equiv 0\,\text{(mod 2)}\,\}, \\
T(m,1)&=&\{\,(s_0,s_1)\,|\,s_0+s_1=m,\ s_1\equiv 1\,\text{(mod 2)}\,\}. 
\end{eqnarray*}} 
Hence, for $g=(g_{0}, g_{1}) \in GL(k) \times GL(l)$, 
we have 
\begin{equation*}
\begin{aligned}
& \tr(g|{\frak L}_{(n,0)})=\text {str}(g|{\frak L}_{(n,0)}) \\
 &= \sum_{d|(n,0)}\frac{1}{d}\m(d)
    \sum_{s\in T((n,0)/d)}\frac{(|s|-1)!}{s!}
     \text {str}(g_0^d|V_0)^{s_0}\text {str}(g_1^d|V_1)^{s_1} \\
 &= \sum_{d|(n,0)}\frac{1}{d}\m(d)
    \sum_{s\in T((n,0)/d)}\frac{(|s|-1)!}{s!}
    (-1)^{s_1}p_{d}(x)^{s_0}p_{d}(y)^{s_1} \\                            
 &=     \sum_{d|(n,0)}\frac{1}{d}\m(d)
    \sum_{s\in T((n,0)/d)}\frac{(|s|-1)!}{s!}(-1)^{s_1}
        \left(\sum_{\ld\vdash ds_0}\chi^{(d^{s_0})}_{\ld}S_{\ld}(x)\right)
        \left(\sum_{\m\vdash ds_1}\chi^{(d^{s_1})}_{\m}S_{\m}(y)\right) \\
 &= \sum_{d|(n,0)}\frac{1}{d}\m(d)
    \sum_{s\in T((n,0)/d)}\frac{(|s|-1)!}{s!}(-1)^{s_1}
    \left(\sum \Sb \ld\vdash ds_0 \\ \mu \vdash ds_1 \endSb 
\chi^{(d^{s_0})}_{\ld}\chi^{(d^{s_1})}_{\m}  S_{\ld}(x)S_{\m}(y)\right),
\end{aligned}
\end{equation*}
and
\begin{equation*}
\begin{aligned}
& \tr(g|{\frak L}_{(n,1)})=-\text {str}(g|{\frak L}_{(n,1)}) \\
& = -\sum_{d|(n,1)}\frac{1}{d}\m(d)
     \sum_{s\in T((n,1)/d)}\frac{(|s|-1)!}{s!}
     \text {str}(g_0^d|V_0)^{s_0}\text {str}(g_1^d|V_1)^{s_1} \\
 & = -\sum \Sb d|n \\ \text{$d$:odd} \endSb 
     \frac{1}{d}\m(d)
     \sum_{s\in T((n,1)/d)}\frac{(|s|-1)!}{s!}
     (-1)^{s_1} \tr(g_0^d|V_0)^{s_0}\tr(g_1^d|V_1)^{s_1} \\      
 & = \sum \Sb  d|n \\ \text{$d$:odd} \endSb 
     \frac{1}{d}\m(d)
     \sum_{s\in T((n,1)/d)}\frac{(|s|-1)!}{s!}
    p_{d}(x)^{s_0}p_{d}(y)^{s_1} \\                              
 & =    \sum \Sb d|n \\ \text{$d$:odd} \endSb 
     \frac{1}{d}\m(d)
     \sum_{s\in T((n,1)/d)}\frac{(|s|-1)!}{s!}
        \left(\sum_{\ld\vdash ds_0}\chi^{(d^{s_0})}_{\ld}S_{\ld}(x)\right)
        \left(\sum_{\m\vdash ds_1}\chi^{(d^{s_1})}_{\m}S_{\m}(y)\right) \\
 & = \sum \Sb d|n \\ \text{$d$:odd}  \endSb 
     \frac{1}{d}\m(d)
     \sum_{s\in T((n,1)/d)}\frac{(|s|-1)!}{s!}
     \left(\sum \Sb  \ld\vdash ds_0 \\  \m\vdash ds_1 \endSb
  \chi^{(d^{s_0})}_{\ld}\chi^{(d^{s_1})}_{\m} S_{\ld}(x)S_{\m}(y)\right).
\end{aligned}
\end{equation*}
It follows that 
\begin{equation}
\begin{aligned}
& \  \tr(g|{\frak L}_n)=\tr(g|{\frak L}_{(n,0)})+\tr(g|{\frak L}_{(n,1)}) \\
 &=\sum \Sb d|n \\  \text{$d$:odd} \endSb 
    \frac{1}{d}\m(d)
    \sum_{s_0+s_1=n/d}\frac{(|s|-1)!}{s!}
    \left(\sum \Sb \ld\vdash ds_0 \\ \m\vdash ds_1 \endSb
  \chi^{(d^{s_0})}_{\ld}\chi^{(d^{s_1})}_{\m}  S_{\ld}(x)S_{\m}(y)
        \right)  \\
 &+\sum \Sb   d|n \\ \text{$d$:even} \endSb 
    \frac{1}{d}\m(d)
    \sum_{s_0+s_1=n/d}(-1)^{s_1}\frac{(|s|-1)!}{s!}
    \left(\sum \Sb  \ld\vdash ds_0 \\ \m\vdash ds_1 \endSb
  \chi^{(d^{s_0})}_{\ld}\chi^{(d^{s_1})}_{\m}
                  S_{\ld}(x)S_{\m}(y)
        \right)  \\
 &=\sum \Sb   l(\ld)\leq k \\  l(\m)\leq l \\  |\ld|+|\m|=n \endSb
     a_{\ld\m}S_{\ld}(x)S_{\m}(y),
\end{aligned}
\end{equation}
where 
\begin{equation*}
\begin{aligned}
a_{\ld\m}=&
\sum \Sb d|n \\ \text{$d$:odd} \endSb 
\frac{\m(d)}{d}
\sum \Sb s_0+s_1=n/d \\  ds_0=|\ld| \endSb 
\frac{(|s|-1)!}{s!}\chi^{(d^{s_0})}_{\ld}\chi^{(d^{s_1})}_{\m} \\
 &+\sum \Sb  d|n \\ \text{$d$:even} \endSb 
\frac{\m(d)}{d}
\sum \Sb  s_0+s_1=n/d \\  ds_0=|\ld| \endSb 
(-1)^{s_1}\frac{(|s|-1)!}{s!}\chi^{(d^{s_0})}_{\ld}\chi^{(d^{s_1})}_{\m} \\
=&\frac{1}{n}\sum_{d|\,|\ld|,|\m|}\m(d)
   \frac{(n/d)!}{(|\ld|/d)!(|\m|/d)!}(-1)^{(d-1)\frac{|\m|}{d}}
   \chi_{\ld}^{(d^{|\ld|/d})}\chi_{\m}^{(d^{|\m|/d})}.   
\end{aligned}
\end{equation*}           

\vskip 2mm 

Note that $a_{\ld\m}$ is the multiplicity of $W_{\ld}\otimes W_{\m}$ 
in ${\frak L}_n$, where $W_{\ld}$ (resp. $W_{\m}$) is the irreducible 
$GL(k)$-module (resp. $GL(l)$-module) corresponding to $\ld$ (resp. $\m$).  
Recall that $c_{\ld}$ is the multiplicity of $V_{\ld}$ in 
${\frak L}_n$.
By Proposition \ref {Prop 5.2}, we have 
\begin{equation} \label{5.14}
\tr(g|{\frak L}_n)=\sum_{\ld\in H}c_{\ld}HS_{\ld}(x,y).
\end{equation}
Therefore, we obtain:

\vskip 2mm
\begin{prop} \label {Prop 5.3}

If $\ld=(\tau_1, \tau_2, \cdots, \tau_{k}, \tau_{k+1}, \cdots)$, 
then we have 
\begin{equation} \label{5.15}
c_{\ld}=a_{\ld_0\ld_1}
       -\sum \Sb \m\in H-\{\ld\} \\  \m_0 > \ld_0  \endSb
        c_{\m}N^{\m'}_{\ld_0'\ld_1},
\end{equation}
where $\ld_0=(\tau_1,\cdots,\tau_k)$ and $\ld_1=(\tau_{k+1},\cdots)'$.
\end{prop}
\noindent
{\it Proof.} \ 
Comparing the equations (\ref{5.14}) and (\ref{5.15}), we have 
\begin{equation*}
\begin{aligned}
& \sum_{\m\in H}c_{\m}HS_{\m}(x,y)
= c_{\ld}HS_{\ld}(x,y)+\sum_{\m\neq\ld}c_{\m}HS_{\m}(x,y) \\
&= c_{\ld}\left(\sum_{\n<\ld}S_{\n}(x)S_{\ld'/\n'}(y)\right)
    +\sum_{\m\neq\ld}c_{\m}HS_{\m}(x,y) \\
&= c_{\ld}\left(S_{\ld_0}(x)S_{\ld_1}(y)+(*)\right)
    +\sum_{\m\neq\ld}c_{\m}
        \left(\sum_{\tau<\m}S_{\tau}(x)
        \sum_{\sigma}N^{\m'}_{\tau'\sigma}S_{\sigma}(y)\right) \\
&= \left(c_{\ld}+\sum \Sb  \m\in H \\ \m_0>\ld_0 \endSb 
                 c_{\m}N^{\m'}_{\ld_0'\ld_1}  \right)
S_{\ld_0}(x)S_{\ld_1}(y)+(*)' \\
&=      \sum \Sb  l(\ld)\leq k \\  l(\m)\leq l \\ |\ld|+|\m|=n  \endSb
     a_{\ld\m}S_{\ld}(x)S_{\m}(y).
\end{aligned}                                                   
\end{equation*}

\vskip 2mm 
Since $\{\,S_{\ld}(x)S_{\m}(y)\,|\,l(\ld)\leq k,\,l(\m)\leq l\,\}$ is
linearly independent (\cite{BR}), we obtain 
\begin{equation*} 
a_{\ld_0\ld_1}=c_{\ld}
       +\sum \Sb \m\in H-\{\ld\} \\  \m_0 > \ld_0  \endSb 
        c_{\m}N^{\m'}_{\ld_0'\ld_1}, 
\end{equation*}
which completes the proof. 
\qed

\vskip 3mm

If $l(\ld)\leq k$ (or $|\ld_0|=n$) for $\ld\in H(k,l;n)$, then 
(\ref{5.15}) yields
\begin{equation} \label{5.16}
c_{\ld}=a_{\ld_0\ld_1}=\frac{1}{n}\sum_{d|n}\m(d) \chi^{(d^{n/d})}_{\ld}.
\end{equation}
If all $c_{\ld}$'s are known for $m\leq |\ld_0|\leq n$, then
for any $\ld\in H(k,l;n)$ with $|\ld_0|=m-1$, we have 
\begin{equation} \label {5.17}
\begin{aligned}
c_{\ld}&=a_{\ld_0\ld_1}
          -\sum \Sb  \m\in H-\{\ld\} \\  \m_0 > \ld_0  \endSb
                            c_{\m}N^{\m'}_{\ld_0'\ld_1} \\
      &=a_{\ld_0\ld_1}
         -\sum \Sb   \m\in H \\ \m_0 > \ld_0    \\ |\m_0|\geq m \endSb
                            c_{\m}N^{\m'}_{\ld_0'\ld_1},                
\end{aligned}
\end{equation}
since $\m_0 > \ld_0$ and $|\m_0|=m-1$ imply $\m_0=\ld_0$, 
and in this case, 
$N^{\m'}_{\ld_0'\ld_1}=0$ unless $\ld=\m$ (cf. \cite{Mac79}).
Following the above inductive step,
 we can completely determine the value of all $c_{\ld}$ 
for $\ld\in H(k,l;n)$:

\vskip 2mm 

\begin{prop} \label {Prop 5.4}

Under the above hypothesis, we have 
\begin{equation}
{\frak L}_n(V)=\bigoplus_{\ld\in H(k,l;n)}V_{\ld}^{\oplus c_{\ld}},
\end{equation}
where $V_{\ld}$ is the irreducible $gl(k,l)$-module 
corresponding to $\ld$
and $c_{\ld}$ is determined by {\rm (\ref{5.15})}. 
\qed
\end{prop}

\vskip 2mm

\noindent
{\it Remark.} \ Note that the multiplicity $c_{\ld}$ is  given
explicitly in closed form when $l(\ld)\leq k$ 
and is given by a recursive formula otherwise.

\vskip 2mm 
\begin{cor} \label {Cor 5.5} {\rm (cf. \cite{Ba}, \cite{HK})} \ 
Under the above hypothesis, if $V_1=0$, then we have 
\begin{eqnarray}
{\frak L}_n(V)&=&\bigoplus \Sb 
 \ld\vdash n \\ l(\ld)\leq k \endSb 
V_{\ld}^{\oplus c_{\ld}}, \\
c_{\ld}&=&\frac{1}{n}\sum_{d|n}\mu(d)\chi^{(d^{n/d})}_{\ld},
\end{eqnarray}
where $V_{\ld}$ is the irreducible $gl(k)$-module 
corresponding to $\ld$.
\qed
\end{cor}

\vskip 2cm 

\section{Generalized Kac-Moody Superalgebras}

The {\it generalized Kac-Moody superalgebras} arise naturally
in the context of {\it Monstrous Moonshine} (\cite{B88}, \cite{B92}), 
automorphic forms with infinite product expansions (\cite{B95}, 
\cite{GN1}--\cite{GN3}), 
and the string theory (\cite{HM}). In this section, we give a conjectural
formula for the homology modules over (the negative part of) 
generalized Kac-Moody superalgebras, and derive a closed form 
supertrace formula for generalized Kac-Moody superalgebras with
group actions.

\vskip 2mm
Let $I$ be a countable (possibly infinite) index set.
A real square matrix $A=(a_{ij})_{i,j\in I}$
is called a \it Borcherds-Cartan matrix \rm if it satisfies:
(i) $a_{ii}=2$ or $a_{ii}\le 0$  for
all $i\in I$, \
(ii) $a_{ij}\le 0$ if $i\neq j $, and $a_{ij}\in \mathbb {Z}$
if $a_{ii}=2$, \
(iii) $a_{ij}=0$ implies $a_{ji}=0.$
We say that an index $i$ is {\it real} if $a_{ii}=2$
and {\it imaginary} if $a_{ii}\le 0$.
We denote by $I^{re}=\{i\in I|\ a_{ii}=2\}$,
$I^{im}=\{i\in I|\ a_{ii}\le 0\}$.
Let $\underline m=(m_{i} \in \mathbb {Z}_{>0} |\ i\in I)$
be a sequence of positive integers
such that $m_i=1$ for all $i\in I^{re}$.
We call $\underline {m}$ the {\it charge} of $A$.
In this paper, we assume that the Borcherds-Cartan matrix $A$ is
{\it symmetrizable}, i.e., there is a diagonal
matrix $D=diag(s_{i}|\ i\in I)$ with $s_i >0$ ($i\in I$)
such that $DA$ is symmetric. 

\vskip 2mm
Let $C=(\theta_{ij})_{i,j\in I}$ be a complex matrix satisfying
$\theta_{ij} \theta_{ji}=1$ for all $i,j\in I$. Thus we have $\theta_{ii}
=\pm 1$ for all $i\in I$. We call $i\in I$ an {\it even index}
if $\theta_{ii}=1$ and an {\it odd index} if $\theta_{ii}=-1$. 
We denote by $I^{even}$ (resp. $I^{odd}$) the set of all even
(resp. odd) indices. 
We say that a Borcherds-Cartan matrix $A=(a_{ij})_{i,j\in I}$
is {\it colored by $C$} if $a_{ii}=2$ and 
$\theta_{ii}=-1$ imply $a_{ij}$ are even integers
for all $j\in I$.
In this case, the matrix $C$ is called a {\it coloring matrix} of $A$.

\vskip 2mm
Let ${\frak h}=(\bigoplus_{i\in I} \mathbb {C} h_i)
\oplus (\bigoplus_{i\in I} \mathbb {C} d_i)$ be a complex vector
space with a basis $\{ h_i, d_i| \ i\in I \}$, and for each $i\in I$
define a linear functional $\alpha_i \in {\frak h}^*$ by
\begin{equation} \label {6.1}
\alpha_i(h_j)=a_{ji}, \ \ \alpha_i(d_j)=\delta_{ij} \ \
\text {for all} \ \ j\in I. 
\end{equation}
The free abelian group $Q=\bigoplus_{i\in I} \mathbb {Z} \alpha_i$
generated by $\alpha_i$'s $(i\in I)$ is called the {\it root lattice} 
associated with $A$. 
Let $\Pi=\{ \alpha_i | \ i \in I \}$ and $B$ be a basis of ${\frak h}^*$
extending $\Pi$. Set $B'=B \setminus \Pi$. 
Since $A$ is assumed to be symmetrizable, 
there is a symmetric bilinear form $(\ |\ )$ on ${\frak h}^*$ defined by
\begin{equation} \label {6.2}
\begin{aligned}
& (\alpha_i|\alpha_j)=s_i a_{ij} \ \  
\text {for} \ i,j \in I, \\
& (\lambda | \alpha_i)= \lambda(s_i h_i) \ \ \text {for} \ 
\lambda \in B', \\
& (\lambda | \mu) =0 \ \ \text {for} \ \lambda, \mu \in B'
\end{aligned}
\end{equation}
(cf. see, for example, \cite{Im}).
Let $Q^{+}=\sum_{i\in I} \mathbb {Z}_{\ge 0} \alpha_i$ 
and $Q^{-}=-Q^{+}$.
There is a partial ordering $\ge$ on $Q$ given by
$\lambda \ge \mu$ if and only if $\lambda-\mu \in Q^{+}$.
The coloring matrix $C=(\theta_{ij})_{i,j\in I}$ defines a bimultiplicative
map $\theta: Q \times Q \rightarrow \mathbb {C}^{\times}$ on $Q$ by
\begin{equation} \label {6.3}
\begin{aligned}
& \theta (\alpha_i, \alpha_j)=\theta_{ij} \ \ \text {for all} \ i,j \in I, \\
& \theta (\alpha+\beta, \gamma)=\theta(\alpha, \gamma)
\theta (\beta, \gamma),\\
& \theta (\alpha, \beta+\gamma)
=\theta(\alpha, \beta) \theta(\alpha, \gamma)
\end{aligned}
\end{equation}
for all $\alpha, \beta, \gamma \in Q$. 
Note that, since $\theta_{ij}\theta_{ji}=1$ for all $i,j\in I$,
$\theta$ satisfies 
\begin{equation} \label {6.4}
\theta(\alpha, \beta) \theta(\beta,\alpha)=1 \ \ \text {for all} \ \
\alpha, \beta \in Q. 
\end{equation}
That is, $\theta$ is a coloring map on $Q$. 
In particular, $\theta(\alpha,\alpha)=\pm 1$ for all $\alpha \in Q$.
We say $\alpha \in Q$ is {\it even} if $\psi(\alpha)=\theta(\alpha, \alpha)=1$
and {\it odd} if $\psi(\alpha)=\theta(\alpha, \alpha)=-1$.

\vskip 2mm
The {\it generalized Kac-Moody superalgebra}
${\frak g}={\frak g}(A,\underline m, C)$ associated with a symmetrizable 
Borcherds-Cartan matrix $A=(a_{ij})_{i,j \in I}$ of charge 
$\underline m=(m_i|\ i \in I)$ with a 
coloring matrix $C=(\theta_{ij})_{i,j\in I}$
is the $\theta$-colored Lie superalgebra generated by the
elements $h_i, d_i$ ($i\in I$), $e_{ik}$, $f_{ik}$ 
$(i\in I, \ k=1,2,\cdots,m_i)$ with the defining relations:
\begin{equation} \label {6.5}
\begin{aligned}
&[h_i,h_j]=[h_i, d_j]=[d_i,d_j] =0, \\
&[h_i,e_{jl}]=a_{ij} e_{jl},
\ \ [h_i,f_{jl}]=-a_{ij} f_{jl}, \\
&[d_i, e_{jl}]=\delta_{ij}e_{jl},
\ \ [d_i, f_{jl}]=-\delta_{ij}f_{jl},\\
&[e_{ik},f_{jl}]=\delta_{ij}\delta_{kl}h_{i}, \\
&(ad e_{ik})^{1-a_{ij}}(e_{jl})
=(ad f_{ik})^{1-a_{ij}}(f_{jl})=0 \ \  
\text {if}\ a_{ii}=2 \ \text {and} \  i\neq j,\\
&[e_{ik},e_{jl}]=[f_{ik},f_{jl}]=0 \ \ \text {if} \ a_{ij}=0
\end{aligned} 
\end{equation}
$\text {for}\  i,j \in I, \ k=1,\cdots, m_i, \ l=1,\cdots,m_j.$

\vskip 2mm 
The abelian subalgebra ${\frak h}=\left(\bigoplus_{i\in I} \mathbb {C}
h_i \right) \bigoplus \left (\bigoplus_{i\in I} \mathbb {C}d_i \right)$
is called
the {\it Cartan subalgebra} of ${\frak g}$,
and the linear functionals $\alpha_i \in {\frak h}^*$ $(i\in I)$
defined by (\ref{6.1}) are called the {\it simple roots} of ${\frak g}$.
For each $i\in I^{re}$, let $r_i \in GL({\frak h}^*)$ be the
{\it simple reflection} on ${\frak h}^*$ defined by
\begin{equation*}
r_i(\lambda)=\lambda-\lambda(h_i) \alpha_i \ \ \text {for} \ \ 
\lambda \in {\frak h}^*. 
\end{equation*}
The subgroup $W$ of $GL({\frak h}^*)$ generated by the $r_i$'s
$(i\in I^{re})$ is called the {\it Weyl group} of the
generalized Kac-Moody superalgebra ${\frak g}$.

\vskip 2mm
The generalized Kac-Moody superalgebra ${\frak g}=
{\frak g}(A, \underline m, C)$
has the {\it root space decomposition}
${\frak g}=\bigoplus_{\alpha\in Q} {\frak g}_{\alpha},$
where
\begin{equation*}
{\frak g}_{\alpha}=\{x\in {\frak g}|\ [h,x]=\alpha(h)x\ \ \text {for all}
\ h\in {\frak h}\}. 
\end{equation*}
Note that
${\frak g}_{\alpha_i}=\mathbb {C}e_{i,1}\oplus \cdots \oplus
\mathbb {C}e_{i,m_i}$ and
${\frak g}_{-\alpha_i}=\mathbb {C}f_{i,1}\oplus \cdots \oplus
\mathbb {C}f_{i,m_i}$.
We say that $\alpha\in Q^{\times}$ is a {\it root} of ${\frak g}$ if
${\frak g}_{\alpha}\neq 0$.
The subspace ${\frak g}_{\alpha}$ is called the {\it root space}
of ${\frak g}$ attached to $\alpha$.
A root $\alpha$ is called {\it real} if $(\alpha|\alpha)>0$
and {\it imaginary} if $(\alpha|\alpha) \le 0$.
In particular, a simple root $\alpha_i$ is real if $a_{ii}=2$
(i.e., $i\in I^{re}$) and imaginary if $a_{ii}\le 0$ (i.e., $i\in I^{im}$).
Note that the imaginary simple roots may have multiplicity $> 1$.
A root $\alpha>0$ (resp. $\alpha<0$) is called {\it positive}
(resp. {\it negative}).
One can show that all the roots are either positive or negative.
We denote by $\Phi$, $\Phi^+$, and $\Phi^-$
the set of all  roots, positive roots, and negative roots, respectively.
We also denote by $\Phi_{0}$ (resp. $\Phi_{1}$) the set of
all even (resp. odd) roots of ${\frak g}$. Hence, for example,
$\Phi^{+}_{0}$ will denote the set of all positive even roots
of ${\frak g}$.
Define the subalgebras ${\frak g}^{\pm}=\bigoplus_{\alpha \in \Phi^{\pm}}
{\frak g}_{\alpha}$. Then we have the {\it triangular decomposition}\,:
$${\frak g}={\frak g}^{-} \oplus {\frak h} \oplus {\frak g}^{+}.$$

\vskip 2mm
We can define a nondegenerate symmetric bilinear form on ${\frak h}$
by
\begin{equation} \label {6.5-1}
(h_i | h)=\frac{1}{s_i} \alpha_i (h) \ \ \text {and} \ \ 
(d_i| d_j)=0
\end{equation}
for all $h\in {\frak h}$, $i, j \in I$ (cf. \cite{K90}). 
The symmetric bilinear form $(\ |\ )$ on ${\frak h}$ defined
by (\ref{6.5-1}) can be
extended to a nondegenerate, {\it supersymmetric}, and {\it
invariant} bilinear form on the generalized Kac-Moody
superalgebra ${\frak g}$:

\begin{prop} \label{Prop 6.0} {\rm (cf. \cite{Im}, \cite{K78})} \ 
Let ${\frak g}={\frak g}(A, \underline {m}, C)$ be the
generalized Kac-Moody superalgebra associated with a
Borcherds-Cartan data $(A, \underline {m}, C)$. 
Then there exists a nondegenerate bilinear form $(\ | \ )$
on ${\frak g}$ satisfying the following conditions\,{\rm :}

{\rm (a)} The bilinear form $(\ |\ )$ is supersymmetric, i.e., we have
$$(x|y)=\theta (\beta, \alpha) (y|x) \ \ \text {for} \ 
x\in {\frak g}_{\alpha}, y\in {\frak g}_{\beta}.$$

{\rm (b)} The bilinear form $(\ |\ )$ is invariant, .i.e., 
we have 
$$([x,y]|z)=(x| [y, z]) \ \ \text {for} \ 
x, y, z \in {\frak g}. $$

{\rm (c)} The bilinear form  $(\ |\ )|_{{\frak h}}$ 
is given by {\rm (\ref{6.5-1})}.

\vskip 2mm

{\rm (d)} $({\frak g}_{\alpha} | {\frak g}_{\beta})=0$ if 
$\alpha + \beta \neq 0$.

\vskip 2mm

{\rm (e)} The root spaces ${\frak g}_{\alpha}$ and ${\frak g}_{-\alpha}$
are nondegenerately paired with respect to $(\ |\ )$, i.e., the 
bilinear form $(\ | \ )|_{{\frak g}_{\alpha}+{\frak g}_{-\alpha}}$is 
nondegenerate.

\vskip 2mm

{\rm (f)} If $\phi: Q \rightarrow {\frak h}$ is a linear map satisfying 
$\phi(\alpha_i)=s_i h_i$ for all $i\in I$, then we have 
$$[x, y]=\theta(\alpha, \alpha) (x|y) \phi(\alpha) 
\ \ \text {for} \ x\in {\frak g}_{\alpha}, 
y\in {\frak g}_{-\alpha}. \ \ \ \ \ \ \ \ \ \ \qed$$ 
\end{prop}

\vskip 2mm
A ${\frak g}$-module $V$ is
called {\it ${\frak h}$-diagonalizable}
if it admits a {\it weight space decomposition}
$V=\bigoplus_{\mu\in {\frak h}^*} V_{\mu}$, where
\begin{equation*}
V_{\mu}=\{ v\in V | \ h\cdot v=\mu(h) v \ \ \text {for all} \ h\in
{\frak h} \}. 
\end{equation*}
If $V_{\mu} \neq 0$, then $\mu$ is called a {\it weight} of $V$
and $V_{\mu}$ is called the {\it $\mu$-weight space}.
We denote by $P(V)$ the set of all weights of $V$.
When all the weight spaces are finite dimensional,
we define the {\it character} of $V$ to be
\begin{equation} \label {6.6}
\text {ch} V=\sum_{\mu\in {\frak h}^*} \left(\text {dim}
V_{\mu}\right) e^{\mu}, 
\end{equation}
where $e^{\mu}$ are the basis elements of the group algebra 
{\bf C}[${\frak h}^*$] with the multiplication given by
$e^{\mu}e^{\nu}=e^{\mu+\nu}$ for $\mu, \nu \in {\frak h}^*$.

\vskip 2mm 
We denote by ${\cal O}$ the category of ${\frak h}$-diagonalizable
${\frak g}$-modules with finite dimensional weight spaces such that
there exist a finite number of linear functionals $\lambda_1, \cdots,
\lambda_s$ satisfying $P(V) \subset \cup_{i=1}^s D(\lambda_i)$,
where $D(\lambda)=\{ \mu \in {\frak h}^*| \ \mu\le \lambda\}$.
The morphisms in ${\cal O}$ are the usual ${\frak g}$-module 
homomorphisms. 

\vskip 2mm
The most important example of the ${\frak g}$-modules in category
${\cal O}$ may be the class of highest weight modules. 
An ${\frak h}$-diagonalizable
${\frak g}$-module $V$ is called a {\it highest
weight module} with highest weight $\lambda\in {\frak h}^*$ if there is a
nonzero vector $v_{\lambda} \in V$ such that 
(i) $e_{ik}\cdot v_{\lambda}=0$ for all $i\in I$, $k=1,\cdots,m_i$, 
\ (ii) $h\cdot v_{\lambda}=\lambda(h) v_{\lambda}$ for all $h\in {\frak h}$,
\ (iii) $V=U({\frak g})\cdot v_{\lambda}$. The vector $v_{\lambda}$ is called
a {\it highest weight vector}. For a highest weight module $V$ with highest
weight $\lambda$, we have (i) $V=U({\frak g}^{-})\cdot v_{\lambda}$,
\ (ii) $V=\bigoplus_{\mu\le \lambda} V_{\mu}$,
$V_{\lambda}=\mathbb {C}v_{\lambda}$,
and \ (iii) $\text {dim} V_{\mu} < \infty$ for all $\mu\le \lambda$.

\vskip 2mm
Let ${\frak b}^{+}={\frak h} \oplus {\frak g}^{+}$, and let
$\mathbb {C}_{\lambda}$ be the 1-dimensional ${\frak b}^{+}$-module
defined by $h\cdot 1=\lambda(h) 1$ for all $h\in {\frak h}$ and
${\frak g}^{+}\cdot 1=0$. The induced module
$M(\lambda)=U({\frak g}) \otimes_{U({\frak b}^{+})} \mathbb {C}_{\lambda}$
is called the {\it Verma module} over ${\frak g}$ with highest weight
$\lambda$. Every highest weight ${\frak g}$-module with highest weight
$\lambda$ is a homomorphic image of $M(\lambda)$ and the Verma module
$M(\lambda)$ contains a unique maximal submodule $J(\lambda)$.
Hence the quotient $V(\lambda)=M(\lambda)/J(\lambda)$ is irreducible.

\vskip 2mm
Let $V$ be a ${\frak g}$-module in category ${\cal O}$. We define 
the {\it Casimir operator} on $V$ as follows (cf. \cite{Im}, 
\cite{K78}). 
Take a linear functional $\rho \in {\frak h}^*$ satisfying
$\rho(h_i)= \frac {1}{2} a_{ii}$ for all $i\in I$. 
Such a linear functional is called a {\it Weyl vector} of ${\frak g}$.
Note that $(\rho| \alpha_i)= \rho(s_i h_i)
=\frac{1}{2} (\alpha_i| \alpha_i)$ $(i\in I)$. 
For each positive root $\alpha \in \Phi^{+}$, choose a basis 
$\{e_{\alpha}^{(i)}|\ i=1, \cdots, \text {dim} {\frak g}_{\alpha}\}$ 
of the root space ${\frak g}_{\alpha}$ and let $\{e_{-\alpha}^{(i)}|
\ i=1, \cdots, \text {dim} {\frak g}_{\alpha} \}$ be the {\it dual
basis} of ${\frak g}_{-\alpha}$ in the sense that 
$(e_{\alpha}^{(i)} | e_{-\alpha}^{(j)}) = \delta_{ij} 
\theta (\alpha, \alpha)$. 
Then for a weight vector $v \in V_{\lambda}$ of weight
$\lambda$, the Casimir operator $\Omega$ on $V$ is defined by
\begin{equation} \label {6.7}
\Omega(v)=(\lambda+ 2 \rho | \lambda) v_{\lambda} + 2 \sum_{\alpha \in \Phi^{+}}
\sum_{i=1}^{\text {dim} {\frak g}_{\alpha}} 
e_{-\alpha}^{(i)} e_{\alpha}^{(i)} v_{\lambda}.
\end{equation}
The Casimir operator commutes with the action of ${\frak g}$ on $V$.
Recall that a weight vector $v\in V$ is called a {\it primitive vector}
if there is a ${\frak g}$-submodule $W$ of $V$ such that $v\notin W$
and $e_{ik} \cdot v \in W$ for all $i\in I$, $k=1,2,\cdots, m_i$ 
(cf. \cite{K90}).
The Casimir operator satisfies the following properties:

\begin{prop} \label {Prop:Casimir}

\ \ {\rm (a)} If $V$ is a highest weight ${\frak g}$-module in 
Category ${\cal O}$ with highest weight $\Lambda$, then we have 
\begin{equation} \label {6.8}
\Omega = (\Lambda + 2\rho | \Lambda) I_{V}.
\end{equation}

\ \ {\rm (b)} If $v$ is a primitive vector with weight $\lambda$, 
then there exists a ${\frak g}$-submodule $W \subset V$ such that
$v \notin W$ and
\begin{equation} \label {primitive}
\Omega(v)=(\lambda + 2\rho | \lambda) v \ \ \text {mod} \ W.
\ \ \ \ \ \ \ \ \ \ \qed
\end{equation}

\end{prop}

\vskip 2mm 
Let $P^{+}$ be the set of all 
linear functionals $\lambda\in {\frak h}^*$ satisfying
\begin{equation} \label {6.9}
\cases
\lambda(h_i) \in \mathbb {Z}_{\ge 0} \ \ \text {for all} \
i\in I^{re}, \\
\lambda(h_i) \in 2 \mathbb {Z}_{\ge 0} \ \ \text {for all} \
i\in I^{re} \cap I^{odd}, \\
\lambda(h_i) \ge 0 \ \ \text {for all} \ i\in I^{im}.
\endcases
\end{equation}
The elements of $P^{+}$ are called the {\it dominant integral weights}.
For a dominant integral weight $\Lambda \in P^{+}$, 
let
\begin{equation} \label {6.10}
\begin{aligned}
& \Phi^{+}(\Lambda) =\{ \beta =  \sum_{i\in I^{im}} k_i \alpha_i 
\in Q^{+}| \ (\Lambda | \alpha_i)=0 \ \ \text {for} \ k_i \ge 1,  \\
& \ \ (\alpha_i | \alpha_j)=0 \ \ \text {for} \ k_i, k_j \ge 1, \
i\neq j, 
\ \ (\alpha_i | \alpha_i)=0 \ \ \text {for} \ k_i \ge 2 \}.
\end{aligned}
\end{equation}
For such an element $\beta \in \Phi^{+}(\Lambda)$, we denote 
$|\beta|=\sum_{i\in I^{im}}  k_i$ and 
\begin{equation} \label {6.11}
\varepsilon (\beta)=\prod_{i\in I^{im} \cap I^{even}} \binom {m_i} {k_i}
\prod_{j \in I^{im} \cap I^{odd}} \binom{m_j + k_j -1} {k_j}.
\end{equation}
Then the character of the irreducible highest weight
module $V(\Lambda)$ with highest weight $\Lambda\in P^+$ is determined
by the {\it Weyl-Kac-Borcherds formula}\,:

\vskip 2mm
\begin{prop} \label {Prop 6.1} {\rm (\cite{Mi}, \cite{Ray})} 
\begin{equation} \label {6.12}
\text {ch}V(\Lambda)=\frac
{\prod_{\alpha \in \Phi^{-}_{1}} (1+e^{\alpha})^{\text {dim}
{\frak g}_{\alpha}}}
{\prod_{\alpha \in \Phi^{-}_{0}} (1-e^{\alpha})^{\text {dim}
{\frak g}_{\alpha}}}
\sum \Sb w\in W \\ \beta \in \Phi^{+}(\Lambda) \endSb
(-1)^{l(w)+|\beta|} \varepsilon(\beta)
e^{w(\Lambda+\rho-\beta)-\rho}.
\end{equation}

Letting $\Lambda=0$, we obtain the {\it denominator identity}:

\begin{equation} \label {6.13}
\frac {\prod_{\alpha\in \Phi^{-}_{0}}
(1-e^{\alpha})^{\text {dim}{\frak g}_{\alpha}}}
{\prod_{\alpha\in \Phi^{-}_{1}}
 (1+e^{\alpha})^{\text {dim}{\frak g}_{\alpha}}}
=\sum \Sb w\in W \\ \beta \in \Phi^{+}(0) \endSb
(-1)^{l(w)+|\beta|} \varepsilon(\beta) e^{w(\rho-\beta)-\rho}. 
\ \ \ \ \ \ \qed
\end{equation}
  
\end{prop}

\vskip 2mm
\begin{cor} \label {Cor 6.2}
A highest weight ${\frak g}$-module $V$ with highest weight $\Lambda
\in P^{+}$ and highest weight vector $v_{\Lambda}$ is irreducible 
if and only if it satisfies\,{\rm :}
\begin{equation} \label {6.14}
\begin{aligned}
& f_{ik}^{\Lambda(h_i)+1} \cdot v_{\Lambda}=0 \ \ \text {\rm if}
\ i\in I^{re}, \\
& f_{ik} \cdot v_{\Lambda} =0 \ (k=1, \cdots, m_i) \ \ \text {\rm if}
\ \Lambda (h_i)=0. 
\ \ \ \ \ \ \ \qed
\end{aligned}
\end{equation}
  
\end{cor} 
\vskip 2mm

Suppose a group $G$ acts on the generalized Kac-Moody superalgebra
${\frak g}={\frak g}(A, \underline {m}, C)$ by Lie superalgebra 
automorphisms of ${\frak g}_{\pm}$ 
preserving the root space decomposition.
Assume further that $\text {str}(g|{\frak g}_{\alpha})
=\text {str}(g|{\frak g}_{-\alpha})$ for all $g\in G$ and $\alpha \in Q^{+}$. 
We would like to apply our supertrace formula (\ref{2.17}) to derive a closed
form formula for $\text {str}(g| {\frak g}_{\alpha})$ $(g\in G, \alpha \in Q)$.

\vskip 2mm 
Let $J$ be a finite subset of $I^{re}$, and let 
$\Phi_{J}=\Phi \cap (\sum_{j\in J} \mathbb {Z} \alpha_j)$,
$\Phi^{\pm}_{J}=\Phi^{\pm} \cap \Phi_{J}$,
and
$\Phi^{\pm}(J)=\Phi^{\pm} \setminus \Phi^{\pm}_{J}$.
We also denote 
$Q_{J}=Q \cap (\sum_{j\in J} \mathbb {Z} \alpha_j)$,
$Q^{\pm}_{J}=Q^{\pm} \cap Q_{J}$,
and
$Q^{\pm}(J)=Q^{\pm} \setminus Q^{\pm}_{J}$.
Let ${\frak g}_{0}^{(J)}={\frak h} \oplus \left(\bigoplus_{\alpha
\in \Phi_{J}} {\frak g}_{\alpha} \right),$
and ${\frak g}_{\pm}^{(J)}=\bigoplus_{\alpha \in \Phi^{\pm}(J)}
{\frak g}_{\alpha}$.
Then we have the {\it triangular decomposition}\,:
\begin{equation} \label {6.15}
{\frak g}={\frak g}_{-}^{(J)} \oplus {\frak g}_{0}^{(J)}
\oplus {\frak g}_{+}^{(J)}, 
\end{equation}
where ${\frak g}_{0}^{(J)}$ is the 
Kac-Moody superalgebra (with an
extended Cartan subalgebra) associated with the generalized 
Cartan matrix $A_{J}=(a_{ij})_{i,j\in J}$ and the set of odd indices
$J^{odd}=J\cap I^{odd}=\{j\in J| \ \theta_{jj}=-1 \}$,
and ${\frak g}_{-}^{(J)}$ (resp. ${\frak g}_{+}^{(J)}$) is
a direct sum of irreducible highest weight (resp. lowest
weight) modules over ${\frak g}_{0}^{(J)}$ (cf. \cite{K78}).

\vskip 2mm 
To apply our supertrace formula (\ref{2.17}) to the Lie superalgebra
${\frak g}_{-}^{(J)}$, we need to compute the generalized 
supercharacters of the homology modules $H_{k}({\frak g}_{-}^{(J)})$.
We conjecture that the ${\frak g}_{0}^{(J)}$-module structure
of $H_{k}({\frak g}_{-}^{(J)})$ is determined by {\it Kostant's
formula}. More precisely, let $W_{J}$ be the subgroup of $W$ 
generated by the simple reflections
$r_j$ with $j\in J$, and let
$W(J)=\{ w\in W|\ \Phi_{w} \subset \Phi^{+}(J) \}$,
where $\Phi_{w}=\{ \alpha \in \Phi^{+}|\ w^{-1} \alpha <0 \}$.
Thus $W_{J}$ is the Weyl group of the Kac-Moody superalgebra
${\frak g}_{0}^{(J)}$ and $W(J)$ is the set of right coset
representatives of $W_{J}$ in $W$. That is, $W=W_{J} W(J)$.
Let us denote by
$\Phi^{\pm}_{J, i}= \Phi_{J} \cap \Phi^{\pm}_{i}$ $(i=0,1)$
and
$\Phi^{\pm}_{i}(J)=\Phi^{\pm}_{i} \setminus \Phi^{\pm}_{J, i}$
$(i=0,1)$.
Then we have the following conjecture for the structure of $H_{k}
({\frak g}_{-}^{(J)})$:

\vskip 3mm
\noindent
{\bf Conjecture:} 
\begin{equation} \label{6.16}
H_{k}({\frak g}_{-}^{(J)}) \cong 
\sum \Sb w\in W(J) \\ \beta \in \Phi^{+}(0)
\\ l(w)+|\beta| =k \endSb
V_{J} (w(\rho-\beta)-\rho)^{\oplus \varepsilon (\beta)},  
\end{equation}
where $V_{J}(\mu)$ denotes the irreducible highest weight module
over the Kac-Moody superalgebra ${\frak g}_{0}^{(J)}$ with
highest weight $\mu$. \hskip 1cm \qed

\vskip 2mm
\noindent
{\it Remark.} \ The formula (\ref{6.16}) was first introduced by Kostant 
for finite dimensional simple Lie algebras (\cite{Ko}). It was generalized
to symmetrizable Kac-Moody algebras (\cite{GL}, \cite{Li}) and
generalized Kac-Moody algebras (\cite{N}). But, in general, 
it is not true for
Lie superalgebras. Still, we conjecture that Kostant's formula 
holds for generalized Kac-Moody superalgebras associated with the
Borcherds-Cartan matrix which is colored by a coloring matrix.

\vskip 2mm
For each $k\ge 1$, let
\begin{equation}  \label {6.17}
H_k^{(J)}=\sum \Sb w\in W(J) \\ \beta \in \Phi^{+}(0) 
\\ l(w)+|\beta|=k \endSb
V_{J} (w(\rho-\beta)-\rho)^{\oplus \varepsilon (\beta)}, 
\end{equation}
and define the {\it homology superspace} 
\begin{equation}  \label {6.18}
H^{(J)}=\sum_{k=1}^{\infty} (-1)^{k+1} H_k^{(J)}
=H_1^{(J)} \ominus H_2^{(J)} \oplus H_3^{(J)} \ominus \cdots, 
\end{equation}
an alternating direct sum of superspaces.
Then, for all $g\in G$ and $\alpha \in Q_{-}$, we have
\begin{equation} \label {6.19}
\text {str}(g|H^{(J)}_{\alpha})=\sum \Sb w\in W(J) \\ \beta \in \Phi^{+}(0)
\\ l(w)+|\beta| \ge 1 \endSb
(-1)^{l(w)+|\beta|+1} \varepsilon (\beta) \,
\text {str}(g|V_{J} (w(\rho-\beta)-\rho)_{\alpha}),
\end{equation}
and the generalized denominator identity can be written as
\begin{equation} \label {6.20}
\begin{aligned}
 \prod_{\alpha \in \Phi^{-}(J)} & \exp \left(- \sum_{k=1}^{\infty}
\frac{1}{k} \text {str}(g^k|{\frak g}_{\alpha}) E^{k\alpha} \right)
=1-\text {sch}_{g} H^{(J)} \\
& =\sum \Sb w\in W(J) \\ \beta \in \Phi^{+}(0) \endSb 
(-1)^{l(w)+|\beta|} \varepsilon (\beta) \,
\text {sch}_{g} V_{J} (w(\rho-\beta)-\rho). 
\end{aligned}
\end{equation}

\vskip 2mm
Let $P(H^{(J)})=\{ \tau \in Q_{-}|\ 
\text {str} (g| H^{(J)}_{\tau}) \neq 0 \}
= \{ \tau_1, \tau_2, \tau_3, \cdots \}$, and 
$t_{g}(i)=\text {str} (g|H^{(J)}_{\tau_i})$.
For $\tau \in Q_{-}$, define the set 
$T^{(J)}(\tau)$ of all partitions of $\tau$ into a sum of $\tau_i$'s
as in (\ref{2.15}) and define the Witt partition function
$W_{g}^{(J)}(\tau)$ as in (\ref{2.16}). 
Then, our supertrace formula (\ref{2.17}) yields:

\vskip 2mm
\noindent
\begin{thm} \label {Thm 6.5}
Suppose that a group $G$ acts on the generalized Kac-Moody superalgebra
${\frak g}={\frak g}(A, \underline {m}, C)$ by Lie superalgebra 
automorphisms of ${\frak g}_{\pm}$ preserving the
root space decomposition. 
Assume further that $\text {str}(g|{\frak g}_{\alpha})
=\text {str}(g|{\frak g}_{-\alpha})$ for all $g\in G$ and $\alpha \in Q^{+}$. 
Let $J$ be a finite subset of $I^{re}$.
Then, for all $g\in G$ and $\alpha \in \Phi^{-}(J)$, we have
\begin{equation} \label {6.21}
\begin{aligned}
& \text {str} (g| {\frak g}_{\alpha}) =\sum \Sb d | \alpha \endSb
\frac {1}{d} \mu(d) W_{g^d}^{(J)} \left( \frac {\alpha}{d} \right) \\
& =\sum_{d|\alpha} \frac {1}{d} \mu(d)
\sum_{s \in T^{(J)}(\alpha /d)}  \frac
{(|s|-1)!}{s!} \prod t_{g^d}(i)^{s_i},
\end{aligned}
\end{equation}
where $\mu$ is the classical M\"obius function, and, for a
positive integer $d$, $d|\alpha$ denotes $\alpha=d\tau$ for some
$\tau \in Q_{-}$, in which case $\alpha /d =\tau$. \qed
\end{thm}

\vskip 2mm
We would like to remind the readers that our
supertrace formula (\ref{6.21}) has been derived under the assumption 
that Kostant's formula holds for generalized Kac-Moody superalgebras.
However, in some special cases, we can derive the supertrace 
formula for generalized Kac-Moody superalgebras wihtout assuming 
Kostant's formula.

\vskip 2mm
For example, suppose that the Borcherds-Cartan matrix 
$A$ satisfies: (i) the set $I^{re}$ is finite,
(ii) $a_{ij} \neq 0$ for all $i,j \in I^{im}$. 
Then it is known that the Lie superalgebra ${\frak g}_{-}^{(J)}
=\bigoplus_{\alpha \in \Delta^{-}(J)} {\frak g}_{\alpha}$ 
(resp. ${\frak g}_{+}^{(J)}
=\bigoplus_{\alpha \in \Delta^{+}(J)} {\frak g}_{\alpha}$)
is isomorphic to the free Lie superalgebra generated by the superspace 
$V=\bigoplus_{j \in I^{im}} V_{J}(-\alpha_j)^{\oplus m_j}$
(resp. $V^*=\bigoplus_{j \in I^{im}} V_{J}^*(-\alpha_j)^{\oplus m_j}$),
where $V_{J}(\mu)$ (resp. $V_{J}^*(\mu)$) denotes the irreducible
highest weight (resp. lowest weight) module over the 
Kac-Moody superalgebra ${\frak g}_{0}^{(J)}$ with highest weight $\mu$
(resp. lowest weight $-\mu$) (cf. \cite{Ka97}). 
Hence, for any $g\in G$, the generalized denominator identity 
of the Lie superalgebra 
${\frak g}_{-}^{(J)}$ is the same as
\begin{equation} \label {6.22}
\begin{aligned}
\prod_{\alpha \in \Phi^{-}(J)} & \exp \left(- \sum_{k=1}^{\infty}
\frac{1}{k} \text {str}(g^k|{\frak g}_{\alpha}) E^{k\alpha} \right)
=1- \text {sch}_{g} V \\
&=1- \sum_{j\in I^{im}} m_{j} \text {sch}_{g} V_{J} (-\alpha_{j}). 
\end{aligned}
\end{equation}
Hence, by our supertrace formula (\ref{2.17}), we obtain 
\begin{equation} \label {6.23}
\text {str} (g| {\frak g}_{\alpha})  =\sum_{d|\alpha} \frac {1}{d} \mu(d)
\sum_{s \in T^{(J)}(\alpha /d)} \frac
{(|s|-1)!}{s!} \prod t_{g^d}(i)^{s_i},
\end{equation}
where $t_{g}(i)=\text {str} (g| V_{\tau_i})
=\sum_{j\in I^{im}} m_{j} \text {str} (g| V_{J} (-\alpha_{j})_{\tau_i})$
for $i=1,2,3,\cdots$.

\vskip 2cm
\section {Monstrous Lie Superalgebras}

Let $F(q)=\sum_{n=-1}^{\infty} f(n) q^n = q^{-1}+f(1)q+f(2) q^2 +
\cdots $ be a normalized $q$-series such that $f(-1)=1$, $f(0)=0$,
and $f(n) \in \mathbb {Z}$ for all $n \ge 1$. 
We say that $F(q)$ is {\it convergent} if it converges for all
$q=e^{2\pi i \tau}$, $\text {Im} \tau >0$.
Observe that, for
each $m\ge 1$, there exists a unique polynomial 
$P_m(t) \in \mathbb {Z} [t]$
such that
$$P_{m}(F) \equiv q^{-m} \ \ \text {mod} \ q \mathbb {Z}[[q]].$$
A convergent normalized $q$-series
$F(q)=\sum_{n=-1}^{\infty} f(n) q^n = q^{-1}+f(1)q+f(2) q^2 +
\cdots $ is called {\it replicable} if for all $m>0$ and 
for all $a|m$, there exist normalized $q$-series
$F^{(a)}(q)=\sum_{n=-1}^{\infty} f^{(a)}(n) q^n 
= q^{-1}+f^{(a)}(1)q+f^{(a)}(2) q^2 + \cdots $
with $f^{(a)}(-1)=1$, $f^{(a)}(0)=0$,
and $f^{(a)} (n) \in \mathbb {Z}$ for all $n \ge 1$
such that 
\begin{equation} \label {7.1}
F^{(1)}=F, \ \ \text {and} \ \ 
\sum \Sb ad=m \\ 0 \le b < d \endSb
F^{(a)} \left( \frac {a\tau+b}{d} \right)
=P_{m}(F), 
\end{equation}
where $q=e^{2\pi i \tau}$, $\text {Im} \tau >0$.

\vskip 2mm
In \cite{Ka97}, the replicable functions were characterized in terms of
product identities. More precisely, a normalized $q$-series
$F(q)=q^{-1}+f(1)q+f(2) q^2 +\cdots =\sum_{n=-1}^{\infty} f(n) q^n$ 
is replicable if and only if 
for all $k\ge 1$, there exist normalized $q$-series
$F^{(k)}(q)=\sum_{n=-1}^{\infty} f^{(k)}(n) q^n$ satisfying
the product identity
\begin{equation} \label {7.2}
\prod_{m,n=1}^{\infty}
\exp \left(-\sum_{k=1}^{\infty} \frac{1}{k} 
f^{(k)}(mn) p^{km} q^{kn} \right)
=1-\sum_{i,j=1}^{\infty} f(i+j-1) p^i q^j. 
\end{equation}

\vskip 2mm
Now, we will consider the generalized Kac-Moody superalgebras 
associated with the normalized $q$-series.
We take $I=\{-1\} \cup \{1,2,3, \cdots\}$ to be
the index set, and let $A=(-(i+j))_{i,j\in I}$ be the
Borcherds-Cartan matrix of the Monster Lie algebra (\cite{B92}). 
Consider a normalized $q$-series $F(q)=\sum_{n=-1}^{\infty}
f(n)q^n$ such that $f(-1)=1$, $f(0)=0$,
and $f(n) \in \mathbb {Z}$ for all $n\ge 1$.
We define the charge of the matrix $A$ to be
$\underline {m}=(|f(i)| : \ i\in I)$, and choose a coloring
matrix $C=(\theta_{ij})_{i,j \in I}$ such that
$\theta_{ii}=1$ if $f(i)>0$ and $\theta_{ii}=-1$ if $f(i)<0$.
That is, an index $i\in I$ is even if $f(i)>0$ and
odd if $f(i) <0$. 

\vskip 2mm

The generalized Kac-Moody superalgebra 
${\frak L}(F)={\frak g}(A, \underline {m}, C)$ 
associated with the above data is called the 
{\it Monstrous Lie superalgebra} associated with the normalized
$q$-series $F(q)=\sum_{n=-1}^{\infty} f(n)q^n$. 
For example, the Monstrous Lie superalgebra associated with
the elliptic modular function $J(q)=j(q)-744$ is the Monster Lie algebra,
and the Monstrous Lie superalgebras associated with the Thompson series
$T_g(q)=\sum_{n=-1}^{\infty} c_g(n) q^n$ are the Monstrous Lie 
superalgebras given in \cite{B92}. 

\vskip 2mm

Let ${\frak L}(F)$ be the Monstrous Lie superalgebra
associated with a normalized $q$-series $F(q)=\sum_{n=-1}^{\infty}
f(n) q^n$. 
By our choice of the charge and the coloring matrix, 
we see that $\alpha_{-1}$ is the only real even simple root,
$\alpha_{i}$ $(i\ge 1)$ is an imaginary even 
simple root with multiplicity $f(i)$  if $f(i)>0$, and 
an imaginary odd simple root with multiplicity $-f(i)$ if $f(i)<0$. 
If we identify the simple roots $\alpha_{-1}$ with $(1,-1) 
\in II_{1,1}$ and $\alpha_i$ with $(1,i) \in II_{1,1}$ $(i\ge 1)$, 
then the Monstrous Lie superalgebra ${\frak L}(F)$ 
becomes a $II_{1,1}$-graded Lie superalgebra, 
and we have
$${\frak L}(F)_{+}=\bigoplus \Sb m>0 \\ n\in \text {\bf Z} \endSb 
{\frak L}(F)_{(m,n)}, \ \ \ \ 
{\frak L}(F)_{-}=\bigoplus \Sb m>0 \\ n\in \text {\bf Z} \endSb 
{\frak L}(F)_{(-m,n)}.$$
Suppose that a group $G$ acts on 
${\frak L}(F)$ by Lie superalgebra automorphisms of ${\frak L}(F)_{\pm} $
preserving the $II_{1,1}$-gradation.
Assume further that $\text {str} (g| {\frak L}(F)_{(m,n)})=
\text {str}(g|{\frak L}(F)_{(-m,-n)})$ for all $g\in G$, $m,n>0$, and  
$$g\cdot f_{-1}=f_{-1}, \ \ g\cdot e_{-1}=e_{-1} \ \ \text {for all}
\ g\in G. $$
We will apply our supertrace formula (\ref{6.21}) (actually (\ref{6.23}))
to this setting. 

\vskip 2mm
Let $V_{-1}={\frak L}_{(-1,1)}=\mathbb {C} f_{-1}$, 
$V_{i}={\frak L}(F)_{(-1,-i)}
=\mathbb{C} f_{i,1} \oplus \cdots \oplus \mathbb {C} f_{i, |f(i)|}$
for $i\ge 1$, and $V=\bigoplus_{n=-1}^{\infty} V_{n}$.
For $g\in G$ and $i\in I=\{-1\} \cup \{1,2,3,\cdots \}$,
define $f_{g}(i)=\text {str}(g| V_{i})$. 
Then $V=\bigoplus_{n=-1}^{\infty} V_{n}$ is a 
$\mathbb {Z}$-graded representation of the group $G$ such that
$\text {sch} (V)=F(q)=\sum_{n=-1}^{\infty} f(n) q^n$ and 
$\text {sch}_{g}(V)= F_{g}(q)=\sum_{n=-1}^{\infty}f_{g}(n) q^n$.

\vskip 2mm
Consider a basis of $V_{i}$ consisting of the eigenvectors 
$v_{i,k}$ of $g\in G$ with eigenvalues $\lambda_{i,k}$
$(k=1,\cdots , |f(i)|)$. Then we have
$$f_{g}(i)=\text {str} (g|V_{i}) 
=\psi(1,i) \sum_{k=1}^{|f(i)|} \lambda_{i,k}.$$
To apply our supertrace formula (\ref{6.21}), take 
$J=\{-1\}$ and consider the corresponding 
triangular decomposition
$${\frak L}(F)={\frak L}(F)_{-}^{(J)} \oplus {\frak L}(F)_{0}^{(J)}
\oplus {\frak L}(F)_{+}^{(J)}.$$
Then ${\frak L}(F)_{0}^{(J)}=\langle e_{-1}, f_{-1},
h_{-1} \rangle + {\frak h} \cong sl(2, \mathbb {C}) +{\frak h}$,
and the subalgebra 
${\frak L}(F)_{-}^{(J)}$ (resp. ${\frak L}(F)_{+}^{(J)}$)
is isomorphic to the free Lie superalgebra generated by the
superspace $H^{(J)}=\bigoplus_{i=1}^{\infty} 
V_{J}(-\alpha_i)^{\oplus |f(i)|}$ (resp. 
$H^{(J)*}=\bigoplus_{i=1}^{\infty} V_{J}^{*}(-\alpha_{i})^{\oplus |f(i)|}$),
where $V_{J}(-\alpha_i)$ and $V_{J}^{*}(-\alpha_{i})$
are the $i$-dimensional irreducible
representations of $sl(2, \mathbb {C})$ 
(since $-\alpha_{i}(h_{-1})=i-1$).
We will concentrate on the Lie superalgebra ${\frak L}(F)_{-}^{(J)}$.

\vskip 2mm
Since the weights of $V_{J}(-\alpha_{i})$ are
$-\alpha_{i}=(-1, -i)$, $-\alpha_i -\alpha_{-1}=(-2, -i+1)$,
$\cdots$, $-\alpha_{i}-(i-1) \alpha_{-1} =(-i, -1)$, 
we have
$$P(H^{(J)})=\{(-i,-j)| \ i, j \in \mathbb {Z}_{>0}\}.$$
Note that each $f_{i,k}$ generates an $i$-dimensional irreducible 
representation of the Lie algebra $sl(2, \mathbb {C})$ generated
by $e_{-1}, f_{-1}, h_{-1}$, and hence so does each
$v_{i,k}$ $(i\ge 1, k=1,2,\cdots, |f(i)|)$. 
Thus we have
\begin{equation*}
\begin{aligned}
H^{(J)}_{(-i,-j)} &=\bigoplus_{k=1}^{|f(i+j-1)|}
\mathbb{C} (adf_{-1})^{i-1} (f_{i+j-1, k}) \\
&=\bigoplus_{k=1}^{|f(i+j-1)|}
\mathbb{C} (adf_{-1})^{i-1} (v_{i+j-1, k}).
\end{aligned}
\end{equation*}
It follows that
$$
\begin{aligned}
g\cdot (ad f_{-1})^{i-1} (v_{i+j-1,k}) &=(ad(g\cdot f_{-1}))^{i-1}
(g\cdot v_{i+j-1,k})\\
&=\lambda_{i+j-1,k} (ad f_{-1})^{i-1} (v_{i+j-1,k}),
\end{aligned}
$$
which yields
$$\text {str}(g|H^{(J)}_{(-i, -j)})
=\psi(1, i+j-1) \sum_{k=1}^{|f(i+j-1)|} \lambda_{i+j-1, k}
=f_{g}(i+j-1)$$
for all $g\in G$, $i,j\in \mathbb{Z}_{>0}$.
Hence the generalized denominator identity for $g\in G$ 
of the Lie superalgebra ${\frak L}(F)_{-}^{(J)}$ is equal to
\begin{equation} \label {7.3}
\begin{aligned}
\prod_{m,n=1} \exp & \left( -\sum_{k=1}^{\infty}
\frac{1}{k} \text {str}(g^k |{\frak L}(F)_{(-m,-n)}^{(J)}) 
p^{km} q^{kn} \right) \\
&=1-\sum_{i,j=1}^{\infty} f_{g}(i+j-1) p^i q^j.
\end{aligned}
\end{equation}

\vskip 2mm
For $k, l>0$, we have
$$T^{(J)}(k,l)=T(k,l)=\{s=(s_{ij})_{i,j\ge 1} | \ 
s_{ij} \in \text {\bf Z}_{\ge 0}, \ 
\sum s_{ij} (i,j)=(k,l) \}, $$
the set of all partitions of $(k,l)$ into a sum of ordered 
pairs of positive integers, and the Witt partition function
$W_{g}^{(J)}(k,l)$ is given by 
$$
\begin{aligned}
W_{g}^{(J)}(k,l)&=\sum_{s\in T(k,l)} \frac 
{(|s|-1)!}{s!} \prod \text {str}(g|H^{(J)}_{(-i,-j)})^{s_{ij}} \\
&=\sum_{s\in T(k,l)} \frac 
{(|s|-1)!}{s!} \prod f_{g}(i+j-1)^{s_{ij}}.
\end{aligned}
$$
Therefore, our supertrace formula (\ref{6.21}) yields:

\begin{prop} \label {Prop 7.1}
Let ${\frak L}(F)=\bigoplus_{(m,n) \in II_{1,1}} 
{\frak L}(F)_{(m,n)}$ be the Monstrous Lie superalgebra 
associated with a normalized $q$-series $F(q)=\sum_{n=-1}^{\infty}
f(n) q^n$ and let $V=\bigoplus_{n=-1}^{\infty} V_{n}
=\bigoplus_{n=-1}^{\infty} {\frak L}(F)_{(-1,-n)}$.
Suppose that a group $G$ acts on 
${\frak L}(F)$ by Lie superalgebra automorphisms of ${\frak L}(F)_{\pm} $
preserving the $II_{1,1}$-gradation. Assume 
further that $\text {str} (g| {\frak L}(F)_{(m,n)})=
\text {str}(g|{\frak L}(F)_{(-m,-n)})$ for all $g\in G$, $m,n>0$, and 
$$g\cdot f_{-1}=f_{-1}, \ \ g\cdot e_{-1}=e_{-1} \ \ \text {for all}
\ g\in G. $$
Then, for all $g\in G$ and $m,n>0$, we have 
\begin{equation}\label {7.4}
\begin{aligned}
& \text {str} (g| {\frak L}(F)_{(m,n)})=
\text {str}(g|{\frak L}(F)_{(-m,-n)}) \\
&=\sum_{d|(m,n)} \frac {1}{d} \mu(d)
\sum_{s\in T(\frac {m}{d},\frac {n}{d})} \frac 
{(|s|-1)!}{s!} \prod f_{g^d}(i+j-1)^{s_{ij}},
\end{aligned}
\end{equation}
where $f_{g}(i)=\text {str}(g|V_{i})=\text {str}(g|{\frak L}_{(-1, -i)})$ 
for all $g\in G$ and $i=1,2, 3, \cdots$.  \qed 
\end{prop}

\vskip 2mm
\begin{cor} \label {Cor 7.2}
Suppose that the generalized supercharacters
$F_{g}(q)=\sum_{n=-1}^{\infty} f_{g}(n) q^n$ are replicable
with the replicates $F_{g}^{(a)}$ $(a\ge 1)$ for all $g\in G$.
Then we have
\begin{equation} \label {7.5}
\text {str}(g|{\frak L}(F)_{(m,n)}) =\text {str}(g|{\frak L}(F)_{(-m,-n)}) 
 = \sum \Sb ad | (m,n) \endSb 
\frac{1}{ad} \mu(d) f_{g^d}^{(a)} \left(\frac{mn}{a^2 d^2}\right).
\end{equation}
\end{cor}

\noindent
{\it Proof.} \ Since the generalized supercharacter $F_{g}(q)$
is replicable, by (\ref{7.1}), the generalized denominator
identity for $g\in G$ of the Lie superalgebra ${\frak L}(F)_{-}^{(J)}$
is equal to 
\begin{equation*}
\begin{aligned}
\prod_{m,n=1}^{\infty} \exp & \left(-\sum_{k=1}^{\infty} \frac{1}{k}
\text {str}(g^k| {\frak L}(F)_{(-m,-n)})p^{km} q^{kn} \right) \\
&= 1- \sum_{i,j=1}^{\infty} f_{g}(i+j-1) p^i q^j \\
&= \prod_{m,n=1}^{\infty} \exp \left(-\sum_{a=1}^{\infty} 
\frac{1}{a} f_{g}^{(a)}(mn) p^{am} q^{an} \right).
\end{aligned}
\end{equation*}
By taking the logarithm, we have
\begin{equation*}
\begin{aligned}
\text {log} \prod_{m,n=1}^{\infty} 
\exp & \left(-\sum_{k=1}^{\infty} \frac{1}{k}
\text {str}(g^k| {\frak L}(F)_{(-m,-n)})p^{km} q^{kn} \right)\\
& =- \sum_{m,n=1}^{\infty} \sum_{k=1}^{\infty} \frac{1}{k}
\text {str}(g^k| {\frak L}(F)_{(-m,-n)}) p^{km} q^{kn} \\
&= - \sum_{m,n=1}^{\infty} \sum \Sb k>0 \\ k|(m,n) \endSb 
\frac{1}{k} \text {str} (g^k| {\frak L}(F)_{(-\frac{m}{k}, -\frac{n}{k})})
p^m q^n,
\end{aligned}
\end{equation*}
and
\begin{equation*}
\begin{aligned}
\text {log} \prod_{m,n=1}^{\infty} 
\exp & \left(-\sum_{a=1}^{\infty} \frac{1}{a}
f_{g}^{(a)}(mn) p^{am} q^{an} \right)\\
& =- \sum_{m,n=1}^{\infty} \sum_{a=1}^{\infty} \frac{1}{a}
f_{g}^{(a)} (mn) p^{am} q^{an} \\
&= - \sum_{m,n=1}^{\infty} \sum \Sb a>0 \\ a|(m,n) \endSb 
\frac{1}{a} f_{g}^{(a)} \left( \frac{mn}{a^2} \right) p^m q^n.
\end{aligned}
\end{equation*}
It follows that 
$$\sum \Sb k>0 \\ k|(m,n) \endSb \frac {1}{k} \text {str} (g^k |
{\frak L}(F)_{(-\frac{m}{k}, -\frac{n}{k})})
=\sum \Sb a>0 \\ a|(m,n) \endSb \frac{1}{a} 
f_{g}^{(a)}\left( \frac{mn}{a^2} \right).$$
Therefore, by M\"obius inversion, we obtain
$$\text {str}(g|{\frak L}(F)_{(-m,-n)}) = \sum \Sb ad | (m,n) \endSb 
\frac{1}{ad} \mu(d) f_{g^d}^{(a)} \left(\frac{mn}{a^2 d^2}\right).
\ \ \ \ \ \ \ \qed$$

\vskip 2cm 
\section {Orbit Lie Superalgebras for the Dynkin Diagram Automorphisms
of Borcherds-Cartan Data}

In the following two sections, we will discuss the applications of 
the generalized denominator identities for Dynkin diagram automorphisms
of generalized Kac-Moody superalgebras. 
The notion of {\it orbit Lie superalgebras} will play an important role 
in our discussion. Most of our result is a generalization 
of the corresonding result
on Kac-Moody algebras and generalized Kac-Moody algebras 
obtained in \cite{FSS} and \cite{FRS}. 
Let $A=(a_{ij})_{i,j \in I}$ be a Borcherds-Cartan matrix of charge
$\underline {m} = (m_i| \ i\in I)$ with a coloring matrix 
$C=(\theta_{ij})_{i,j\in I}$, and let ${\frak g}={\frak g}(A, 
\underline {m}, C)$ be the corresponding generalized Kac-Moody superalgebra. 

\vskip 2mm
\noindent
\begin{df} \label{Def 8.0} {\rm A bijection $\sigma: I \rightarrow I$ 
is called a {\it Dynkin diagram automorphism} of the Borcherds-Cartan data 
$(A, \underline {m}, C)$ if it satisfies: \ 
(i) $a_{ij}=a_{\sigma(i), \sigma(j)}$ for all $i, j\in I$, \ 
(ii) $\sigma(I^{even}) \subset I^{even}$, \ 
$\sigma(I^{odd}) \subset I^{odd}$, \ 
(iii) $m_{\sigma(i)}=m_i$ for all $i\in I$. }
\end{df}

\vskip 2mm
\noindent 
In this paper, we will assume that the Dynkin diagram automorphism
has  finite order $N$. 

\vskip 2mm
Given a Dynkin diagram automorphism $\sigma$ of a Borcherds-Cartan data
$(A, \underline {m}, C)$, we extend it to a Lie superalgebra automorphism
of ${\frak g}={\frak g}(A, \underline {m}, C)$ by defining
\begin{equation} \label {8.1}
\begin{aligned}
& \sigma(e_{ik})=e_{\sigma(i), k} \ \ \text {for} \ i\in I, 1 \le k \le m_i, \\
& \sigma(f_{ik})=f_{\sigma(i), k} \ \ \text {for} \ i\in I, 1 \le k \le m_i, \\
& \sigma(h_i)=h_{\sigma(i)}, \ \ \sigma(d_i)=d_{\sigma(i)} \ \ \text {for} \ i\in I.
\end{aligned}
\end{equation}
It is straightforward to check that $\sigma$ is indeed a Lie superalgebra
automorphism of ${\frak g}$. 
Let $U=U({\frak g})$ be the universal enveloping algebra of ${\frak g}$.
Then the Dynkin diagram automorphism $\sigma$ of ${\frak g}$ can be 
extended to an algebra automorphism of $U({\frak g})$ by
$$\sigma(x_1 \cdots x_k)=\sigma(x_1) \cdots \sigma(x_k) \ \ \ \text {for} 
\ \  x_i \in {\frak g},$$
which we will also denote by $\sigma$. 

\vskip 2mm
Moreover, the Dynkin diagram automorphism $\sigma$ induces an 
automorphism $\sigma^*$ of ${\frak h}^*$ defined by
\begin{equation} \label {8.2}
\sigma^*(\lambda) (h)= \lambda(\sigma^{-1}(h)) \ \ \text {for} \ 
\lambda \in {\frak h}^*, \ h\in {\frak h}.
\end{equation}
In particular, for all $i\in I$, we have 
$\sigma^*(\alpha_i)=\alpha_{\sigma(i)}$. 
Since $\sigma$ has order $N$, $\sigma^*$ also has order $N$, and
we have an eigenspace decomposition 
$${\frak h}^*=\bigoplus_{k=0}^{N-1} ({\frak h}^*)^{(k)},$$
where $({\frak h}^*)^{(k)}=\{\lambda \in {\frak h}^*| \ 
\sigma^*(\lambda)=e^{\frac{2 \pi i k}{N}} \lambda \}$.
The elements of $({\frak h}^*)^{(0)}=\{\lambda \in {\frak h}^*| \ 
\sigma^*(\lambda)= \lambda \}$ are called {\it symmetric}. 

\vskip 2mm
Now, we will construct a generalized Kac-Moody superalgebra 
$\breve{\frak g} = {\frak g}(\sigma)$, called the {\it orbit 
Lie superalgebra}, 
associated with the Dynkin diagram 
automorphism $\sigma$ of ${\frak g}$.
Let $\sigma$ be a Dynkin diagram automorphism  
of the Borcherds-Cartan data $(A,\underline {m}, C)$.
Since $\sigma$ permutes the elements of $I$, $I$ is decomposed into a 
disjoint union of $\sigma$-orbits: $I= \bigsqcup [i]$,
where $[i]=\{\sigma^k(i)| \ k=0,1, \cdots, N_{i}-1\}$. 
Let $N_i=|[i]|$ denote the number of elements in the $\sigma$-orbit
$[i]$. 
For each $i\in I$, let $\bar i$ be 
the smallest element of the $\sigma$-orbit of $i$, 
and let $\widehat {I} =\{ \bar i | \ i\in I \}$ 
be the set of such representatives.
Set $\widehat {I}^{re} = \widehat {I} \cap I^{re}$, 
$\widehat {I}^{im} = \widehat {I} \cap I^{im}$,
$\widehat {I}^{even} = \widehat {I} \cap I^{even}$, and
$\widehat {I}^{odd} = \widehat {I} \cap I^{odd}$.

\vskip 2mm 
Let $\widehat{A} =\widehat{A}(\sigma)=(\widehat {a}_{\bar i, \bar j})_{\bar i,
\bar j \in \widehat{I}}$  
be the square matrix
indexed by $\widehat {I}$ whose entries are defined by
\begin{equation} \label {8.3}
\widehat {a}_{\bar i, \bar j}=\epsilon_j 
\sum_{k=0}^{N_j -1} a_{i, \sigma^k(j)} 
\ \ \text {for} \ \bar i, \bar j \in \widehat {I}, 
\end{equation}
where $\epsilon_j = 1-\sum_{k=1}^{N_j -1} a_{j, \sigma^k(j)}$.
Note that $\epsilon_j= \epsilon_{\sigma(j)}$ for all $j\in I$.

\vskip 2mm
%In this paper, we will assume that $\epsilon_i$'s satisfy the following
%condition: 
%\begin{equation} \label {8.4}
%\begin{aligned}
%& \epsilon_i \le 2 \ \ \text {for} \ i\in I^{re}, \\
%& \epsilon_i=1 \ \ \text {for} \ i\in I^{im}.
%\end{aligned}
%\end{equation}
Let 
\begin{equation} \label {8.4}
\breve{I}= 
\{\,i \in \widehat{I}\,|\,\epsilon_i \le 2 \ \ \text {for} \ i\in I^{re}, 
\,\,\epsilon_i=1 \ \ \text {for} \ i\in I^{im}\,\}.
\end{equation}
Since $\widehat {a}_{\bar i, \bar i}=\epsilon_i (a_{ii}+1-\epsilon_i)$,
the condition on $\breve{I}$ is equivalent to 
$\widehat {a}_{\bar i, \bar i}=a_{ii}$ for all $\bar i \in \widehat {I}$.
(The first half of the condition (\ref{8.4}) was introduced in \cite{FSS}
and was referred to as the {\it linking condition}.)
Also, we set 
$\breve{A} =\breve{A}(\sigma)=(\widehat {a}_{\bar i, \bar j})_{\bar i,
\bar j \in \breve{I}}$,  
$\breve{I}^{re} = \breve {I} \cap I^{re}$, 
$\breve {I}^{im} = \breve {I} \cap I^{im}$,
$\breve{I}^{even} = \breve{I} \cap I^{even}$, and
$\breve {I}^{odd} = \breve{I} \cap I^{odd}$.

One can easily verify that 
$\widehat {A}$ is a symmetrizable Borcherds-Cartan matrix. 
For instance, let $\widehat {s}_{\bar i}=N_i\epsilon_i  s_i$
$(\bar i \in \widehat {I})$ and let $\widehat {D}
=diag(\widehat {s}_{\bar i}| \ \bar i \in \widehat {I})$,
then $\widehat D \widehat A$ is symmetric. 

\vskip 2mm
We define the charge $\underline {\widehat m}=
(\widehat {m}_{\bar i} | \ \bar i \in \widehat I)$
and the coloring matrix 
$\widehat C=(\widehat {\theta}_{\bar i, \bar j})_{\bar i,
\bar j \in \widehat I}$ of $\widehat {A}$ by
\begin{equation} \label {8.5}
\widehat {m}_{\bar i}=m_i, \ \ 
\widehat {\theta}_{\bar i, \bar j}=\theta_{ij} \ \ 
\text {for} \ \bar i, \bar j \in \widehat {I}.
\end{equation}
Consider the generalized Kac-Moody superalgebra
$\widehat {\frak g}={\frak g}(\widehat A, \underline {\widehat m},
\widehat C)$ associated with the Borcherds-Cartan data
$(\widehat A, \underline {\widehat m}, \widehat C)$. 
We denote by
$\widehat {e}_{\bar i, k}$,
$\widehat {f}_{\bar i, k}$,
$\widehat {h}_{\bar i}$, $\widehat {d}_{\bar i}$ $(\bar i \in \widehat {I},
k=1, 2, \cdots, m_i)$ the generators of 
Lie superalgebra $\widehat {\frak g}$,  and by 
$\widehat {\frak h}=\left(\bigoplus_{\bar i \in \widehat I} 
\mathbb {C} \widehat h_{\bar i} \right) \oplus 
\left( \bigoplus_{\bar i \in \widehat I}
\mathbb {C} \widehat d_{\bar i} \right)$ the 
Cartan subalgebra of $\widehat {\frak g}$.
We also denote by $\widehat {\Pi}=\{\widehat {\alpha}_{\bar i} |
\ \bar i \in \widehat {I} \}$ the set of simple roots of $\widehat 
{\frak g}$ and by $\widehat {\Phi}$ the set of roots of 
$\widehat {\frak g}$. 
The {\it orbit Lie superalgebra} 
$\breve{\frak g}={\frak g}(\sigma)$ is defined to be 
the subalgebra of $\widehat {\frak g}$ generated by 
$\widehat {e}_{\bar i, k}$,
$\widehat {f}_{\bar i, k}$,$(\bar i \in \breve {I},k=1, 2, \cdots, m_i)$ 
$\widehat {h}_{\bar i}$, $\widehat {d}_{\bar i}$ $(\bar i \in \widehat {I},
k=1, 2, \cdots, m_i)$. 
We denote by $\breve {\Pi}=\{\widehat {\alpha}_{\bar i} |
\ \bar i \in \widehat {I} \}$ the set of simple roots of $\breve 
{\frak g}$ and by $\breve {\Phi}$ the set of roots of 
$\breve {\frak g}$.

\vskip 2mm
As in the case of ${\frak g}={\frak g}(A, \underline{m}, C)$,
there is a symmetric bilinear form on $\widehat {\frak h}^*$
satisfying  
$(\widehat {\lambda} | \widehat {\alpha_{\bar i}})
=\widehat {\lambda} (\widehat {s_{\bar i}} \widehat {h_{\bar i}})$
for $\widehat {\lambda} \in \widehat {\frak h}^*, 
\ \bar i \in \widehat {I}$ (see (\ref{6.2})). 
Define a $\mathbb {C}$-linear isomorphism 
$\phi: ({\frak h}^*)^{(0)}
\rightarrow \widehat {\frak h}^*$  by
\begin{equation} \label {8.7}
\begin{aligned}
& \phi(\lambda) (\widehat {h_{\bar i}}) = \lambda(h_i),\\
& \phi(\lambda) (\widehat {d_{\bar i}}) = \epsilon_{i}^{-1}\lambda(d_i) 
\ \ \text {for} \ \lambda \in ({\frak h}^*)^{(0)}, \ 
\bar i\in \widehat {I}.
\end{aligned}
\end{equation}
Then we obtain the following lemma which will be used in the next section.

\begin{lem} \label {Lemma 8.1}

For $\lambda$, $\mu$ $\in ({\frak h}^*)^{(0)}$, we have
$$(\lambda|\mu)=(\phi(\lambda)|\phi(\mu)).$$
Equivalently, for $\widehat{\lambda}$, $\widehat{\mu}$
$\in\widehat {\frak h}^*$, we have
$$(\widehat {\lambda}|\widehat{\mu})=
(\phi^{-1}(\widehat{\lambda})|\phi^{-1}(\widehat{\mu})). $$

\noindent
{\it Proof.} \ {\rm First, observe that, 
%if we define 
%$\widehat {\alpha}_{\bar i} \in \widehat {\frak h}^*$ by
%$$\widehat {\alpha}_{\bar i} (\widehat {h}_{\bar j})
%=\widehat {a}_{\bar j, \bar i}, \ \ \text {and} 
%\ \ \widehat {\alpha}_{\bar i} (\widehat {d}_{\bar j})
%=\epsilon_i \delta_{ij} \ \ \text {for} \ \ \bar i, \bar j 
%\in \widehat {I},$$  
%we have 
$\phi^{-1}(\widehat {\alpha}_{\bar i})
= \epsilon_i \sum_{k=0}^{N_i -1}  \alpha_{\sigma^k(i)}$ for $\bar i \in
\widehat {I}$.  Hence for $\widehat {\lambda} \in 
\widehat {\frak h}^*$, we obtain 
$$(\widehat {\lambda}| \widehat {\alpha_{\bar i}})
=\widehat {\lambda}(\widehat {s_{\bar i}} \widehat {h_{\bar i}})
=\widehat {s_{\bar i}} \widehat {\lambda} (\widehat {h_{\bar i}}),$$
and 
\begin{equation*}
\begin{aligned}
(\phi^{-1}(\widehat {\lambda})|\phi^{-1}(\widehat {\alpha_{\bar i}}))
& =(\phi^{-1}(\widehat {\lambda})| 
\epsilon_i \sum_{k=0}^{N_i -1} \alpha_{\sigma^k(i)})
=\epsilon_i \sum_{k=0}^{N_i -1} (\phi^{-1} (\widehat {\lambda})|
\alpha_{\sigma^k(i)}) \\
&=\epsilon_{i} \sum_{k=0}^{N_i-1} 
\phi^{-1}(\widehat {\lambda}) (s_{\sigma^k(i)} h_{\sigma^k(i)})
=\epsilon_i N_i \phi^{-1} (\widehat {\lambda}) (s_i h_i) \\
& =\epsilon_i N_i s_i \phi^{-1}(\widehat {\lambda})(h_i) 
=s_i \epsilon_i N_i \widehat {\lambda} (\widehat {h_{\bar i}}),
\end{aligned}
\end{equation*}
which proves our claim.} \qed
\end{lem} 

\vskip 2mm
Observe that if $\epsilon_{i}=1$, 
the $\sigma$-orbit of $i$ in the Dynkin diagram
of $\frak{g}$ is the disjoint union of 
$N_i$ copies of $A_{1}$ diagram or $B(0,1)$ diagram.
If $\epsilon_{i}=2$, then $N_i$ is even and 
the $\sigma$-orbit of $i$ in the Dynkin diagram of $\frak{g}$ is 
the disjoint union of $N_{i}/2$ copies of $A_{2}$ diagram.

\vskip 2mm
Next, we will show that the Weyl group $\breve{W}$ of
$\breve{\frak g}$ is isomorphic to a subgroup of $W$.
We may assume that $\breve{W}$ acts on $\widehat{\frak h}^*$
since $\breve{\frak g}$ is a subalgebra of 
$\widehat{\frak g}$. 
For each $ \bar i\in \breve{I}^{re}$, 
we denote by $\breve{w}_{\bar{i}}$ 
the simple reflection in $\breve{W}$
and define 
\begin{equation} \label {8.8}
\bar w_{\bar i}=\cases \prod_{k=0}^{N_i-1} w_{\sigma^k(i)}
\ \ & \text {if} \ \epsilon_i=1, \\
%w_i w_{\sigma(i)} w_i \ \ & \text {if} \ \epsilon_i=2
\prod_{l=0}^{N_{i}/2 -1}
w_{\sigma^{l}(i)}w_{\sigma^{l+N_{i}/2}(i)}w_{\sigma^{l}(i)}
\ \ & \text {if} \ \epsilon_i=2
\endcases
\end{equation}
where $w_{i}$ is simple reflection in $W$ $(i\in I^{re})$.
Since $a_{i, \sigma^k(i)}=0$ for all $k\not\equiv 0$ $(\text {mod} \ N_i)$
if $\epsilon_i=1$ and $a_{i, \sigma^{N_{i}/2}(i)}= -1$ if 
$\epsilon_{i}=2$, $\bar w_{\bar i}$ is well-defined. 
Let $\bar {W}$ be the subgroup of $W$ generated by
$\bar w_{\bar i}$'s $(\bar i \in \breve {I}^{re})$.

\begin{prop} \label {Prop 8.2} \ For all $\bar i \in \breve {I}^{re}$, 
we have\,{\rm :}

{\rm (a)} $\bar w_{\bar i}^2=1$.

\vskip 2mm

{\rm (b)} $\sigma^* \bar w_{\bar i} = \bar w_{\bar i} \sigma^*$.
In particular, $({\frak h}^*)^{(0)}$ is invariant under $\bar w_{\bar i}$.

\vskip 2mm

{\rm (c)} For each $\lambda \in ({\frak h}^*)^{(0)}$, we have
$$\bar w_{\bar i} (\lambda)=\lambda - \epsilon_i \lambda(h_i) 
\sum_{k=0}^{N_i-1} \alpha_{\sigma^k(i)}.$$

{\rm (d)}  For each simple reflection $\breve {w}_{\bar i} \in \breve {W}$,
we have $\bar w_{\bar i} = \phi^{-1} \breve {w}_{\bar i} \phi$. 

\vskip 2mm
\noindent
{\it Proof.} \ {\rm The proofs of (a), (b) and (c) are straightforward.

\vskip 2mm
For (d), let $\lambda \in ({\frak h}^*)^{(0)}$. Then  we have 
\begin{equation*}
\begin{aligned}
\phi^{-1} \breve {w}_{\bar i} \ \phi(\lambda)
&=\phi^{-1} (\phi(\lambda) - \phi(\lambda) (\widehat h_{\bar i})
\widehat \alpha_{\bar i}) \\
&=\phi^{-1} (\phi(\lambda) -\lambda(h_i) \widehat {\alpha}_{\bar i}) \\
&=\lambda -\lambda(h_i) \phi^{-1} (\widehat {\alpha}_{\bar i})\\
&=\lambda -\lambda(h_i) \epsilon_i \sum_{k=0}^{N_i-1} \alpha_{\sigma^k(i)}\\
&=\bar w_{\bar i} (\lambda),
\end{aligned}
\end{equation*}
which proves our assertion.  \qed}

\end{prop}

\begin{prop} \label {Prop 8.3} \ 
The Weyl group $\breve {W}$ of the orbit Lie superalgebra
$\breve {\frak g}$ is isomorphic to the subgroup 
$\bar W$ of $W$ under the map $\breve{w}_{\bar i} \mapsto 
\bar {w}_{\bar i}$ $(\bar {i} \in \breve{I}^{re})$. 

\noindent
{\it Proof.} \ {\rm 
See \cite{FSS} or \cite{FRS}.}
\qed
\end{prop}

{\it Remark.} \ The original proof of \ref{Prop 8.3} in \cite{FSS} consists of
a somewhat lengthy calculations. In \cite{FRS}, the proof was simplified
and it was also proved that 
$\bar{W}\cong \{\,w\in W\,|\,w\sigma^{*}=\sigma^{*}w\,\}$.

\vskip 2mm Under the above identification of $\bar W$ with $\breve {W}$, 
we define the {\it sign function} $\bar {\varepsilon} : \bar W
\rightarrow \{\pm 1 \}$ of $\bar W$  by
$$\bar {\varepsilon} (\bar w)=\varepsilon (\phi \bar w \phi^{-1})
=\varepsilon(\psi^{-1} (\bar w))
\ \ \text {for} \ \bar w \in \bar W,$$
where $\varepsilon$ denotes the sign function of the Weyl group 
$\breve {W}$ of the orbit Lie superalgebra $\breve {\frak g}
={\frak g}(\sigma)$.
Note that, even though $\bar W$ is a subgroup of $W$, we use the
sign function of $\breve{W}$, not the sign function of $W$. 

\vskip 2cm

\section {Generalized Characters of Highest Weight Modules 
for Dynkin Diagram Automorphisms}

Let $\Lambda \in P^{+}$ be a dominant integral weight and let 
$M(\Lambda)$ be the Verma module with highest weight $\Lambda$. 
By construction of $M(\Lambda)$, we have
$$M(\Lambda) \cong U({\frak g}) \big/ R(\Lambda),$$
where $R(\Lambda)$ is the left ideal of $U({\frak g})$
generated by the elements $e_{ik}$ $(i\in I, k=1, \cdots, m_i)$
and $h-\Lambda(h) 1$ $(h\in {\frak h})$. 
By (\ref{8.1}) and (\ref{8.2}), we see that
$\sigma(R(\Lambda))=R(\sigma^*(\Lambda))$. Hence the 
Dynkin diagram automorphism $\sigma$ induces a 
$\mathbb {C}$-linear isomorphism
\begin{equation}
\sigma: M(\Lambda) \rightarrow M(\sigma^*(\Lambda))
\end{equation}
between Verma modules. 
Note that we have $\sigma(M(\Lambda)_{\lambda})
=M(\sigma^*(\Lambda))_{\sigma^*(\lambda)}$ for all $\lambda 
\le \Lambda$. 

\vskip 2mm
Similarly, let $V(\Lambda)$ be the irreducible highest weight 
module with highest weight $\Lambda$. Then by Corollary \ref 
{Cor 6.2}, the Dynkin diagram automorphism induces a $\mathbb {C}$-linear
isomorphism 
\begin{equation}
\sigma: V(\Lambda) \rightarrow V(\sigma^*(\Lambda))
\end{equation}
between irreducible highest weight modules, and we have
$\sigma(V(\Lambda)_{\lambda})
=V(\sigma^*(\Lambda))_{\sigma^*(\lambda)}$ for all $\lambda 
\le \Lambda$.

\vskip 2mm
In particular, if $\Lambda \in ({\frak h}^*)^{(0)}$, i.e., if
$\sigma^*(\Lambda)=\Lambda$, then 
the Dynkin diagram automorphism $\sigma$ induces 
$\mathbb{C}$-linear automorphisms
on $M(\Lambda)$ and $V(\Lambda)$, and we can consider the generalized
characters for $\sigma$:
\begin{equation} \label {9.3}
\begin{aligned}
& \text {ch}_{\sigma} M(\Lambda)=\sum_{\lambda \le \Lambda}
\tr (\sigma| M(\Lambda)_{\lambda}) e^{\lambda}, \\
& \text {ch}_{\sigma} V(\Lambda)=\sum_{\lambda \le \Lambda}
\tr (\sigma| V(\Lambda)_{\lambda}) e^{\lambda}. 
\end{aligned}
\end{equation}
Note that we have only to take the sum over the weights 
$\lambda$ such that $\sigma^*(\lambda)=\lambda$. 

\vskip 2mm
Recall that the automorphism $\sigma^*$ permutes the roots of ${\frak g}$ in 
such a way that we have
$$\sigma^*(\Phi^{\pm})=\Phi^{\pm}, \ \ 
\sigma^*(\Phi_{0}^{\pm})=\Phi_{0}^{\pm}, \ \ 
\sigma^*(\Phi_{1}^{\pm})=\Phi_{1}^{\pm}.$$
Thus $\Phi^{\pm}$ is a disjoint union of $\sigma^*$-orbits of the elements
in $\Phi^{\pm}$. 
Let $S^{\pm}$ denote the set of representatives of $\sigma^*$-orbits
in $\Phi^{\pm}$ and set
$$S_{0}^{\pm}=S^{\pm} \cap \Phi_{0}, \ \ 
S_{1}^{\pm}=S^{\pm} \cap \Phi_{1}.$$
For each root $\alpha \in \Phi$, let $N_{\alpha}$ denote 
the number of elements in $\sigma^*$-orbit of $\alpha$, and define
\begin{equation}
(\alpha)=\sum_{k=0}^{N_{\alpha}-1} {\sigma^*}^k (\alpha)
\ \ \text {and} \ \ 
{\frak g}_{(\alpha)}=\bigoplus_{k=0}^{N_{\alpha} -1}
{\frak g}_{{\sigma^*}^k(\alpha)}.
\end{equation}
Then, by the Poincar\'e-Birkhoff-Witt Theorem, we have
\begin{equation*}
U({\frak g}^{\pm}) \cong S({\frak g}_{0}^{\pm}) \otimes 
\Lambda({\frak g}_{1}^{\pm}), 
\end{equation*}
where 
$$S({\frak g}_{0}^{\pm}) \cong \bigotimes_{\alpha \in S_{0}^{\pm}}
S({\frak g}_{(\alpha)}), \ \ \text {and} \ \ 
\Lambda({\frak g}_{1}^{\pm}) \cong \bigotimes_{\beta \in S_{1}^{\pm}}
\Lambda({\frak g}_{(\beta)})$$
as $\mathbb {C}$-vector spaces. 
Hence the generalized (super)character of the Verma module
$M(0)$ for the Dynkin diagram automorphism $\sigma$ is given by
\begin{equation*}
\begin{aligned}
 \ch_{\sigma} M(0) &=\text {sch}_{\sigma} M(0) = 
\sum_{\lambda \le 0} \tr(\sigma| M(0)_{\lambda})e^{\lambda}\\
&= \prod_{\alpha \in S_{0}^{+}} \left(\sum_{m\ge 0}
\tr(\sigma| S^{m}({\frak g}_{(-\alpha)})) e^{- m(\alpha) / N_{\alpha}} \right) \\
& \hskip 1cm \times 
 \prod_{\beta \in S_{1}^{+}} \left(\sum_{n\ge 0}
\tr(\sigma| \Lambda^n({\frak g}_{(-\beta)})) e^{-n(\beta)/ N_{\beta}}\right)\\
&= \prod_{\alpha \in S_{0}^{+}} \exp \left(\sum_{m\ge 1} \frac{1}{m} 
\tr(\sigma^m | {\frak g}_{(-\alpha)}) e^{- m(\alpha) / N_{\alpha}} \right)\\
& \hskip 1cm  \times 
\prod_{\beta \in S_{1}^{+}} \exp  \left(- \sum_{n\ge 1} \frac{(-1)^n}{n} 
\tr(\sigma^n | {\frak g}_{(-\beta)}) e^{-n(\beta)/ N_{\beta}}\right)\\
&= \prod_{\alpha \in S^{+}} \exp \left(\sum_{m\ge 1} \frac{1}{m} 
\text {str} (\sigma^m | {\frak g}_{(-\alpha)}) 
E^{- m(\alpha) / N_{\alpha}} \right) \\
&= \prod_{\alpha \in \Phi^{+}} \exp  \left(\sum_{m\ge 1} 
\frac{1}{m N_{\alpha}} \text {str}(\sigma^m | {\frak g}_{(-\alpha)}) 
E^{-m(\alpha)/ N_{\alpha}}\right).
\end{aligned}
\end{equation*}

\vskip 2mm 
Note that if $\rho$ is a Weyl vector of the generalized Kac-Moody 
superalgebra ${\frak g}$, then 
$\phi(\rho)$ is a Weyl vector of the orbit Lie superalgebra
${\breve {\frak g}}$. 

\vskip 2mm
\begin{prop} \label {Prop 9.1} \ 
For all $\bar w \in \bar W$, we have 
\begin{equation}
\bar w \left( e^{-\rho} \ch_{\sigma} M(0) \right)
=\bar \varepsilon (\bar w) e^{-\rho} \ch_{\sigma} M(0).
\end{equation}

\noindent
{\it Proof.} \ {\rm It suffices to prove our assertion for
the simple reflections $\bar w_{\bar i} \in \bar W$
$(\bar i \in \breve {I}^{re})$. 
For each $i\in I^{re}$, let 
$r_{i}=\exp (e_{i}') \exp (-f_{i}') \exp (e_{i}')$, where
the elements $e_{i}'$, $f_{i}'$ and $h_{i}'$ are defined  by
\begin{equation*}
\cases  e_{i}' =   e_{i}, \ \  f_{i}'=f_i, \ \ 
h_{i}'=h_{i} \ \ & \text {if} \  i \in I^{even}, \\ 
 e_{i}'=\frac {1}{4} [e_{i}, e_{i}], \ \ 
f_{i}'=\frac {1}{4} [f_{i}, f_{i}], \ \  
\ \ h_{i}'=\frac{1}{2} h_i \ \ & \text {if} \  i \in I^{odd}.
\endcases 
\end{equation*}
Then for each $i\in I$, 
the elements $e_{i}'$, $f_i'$, and $h_i'$ generate the
subalgebra of ${\frak g}$ isomorphic to $sl(2, \mathbb {C})$,
and by the same argument in \cite{K78}, we can verify that 
$r_{i} \in Aut({\frak g})$ and $r_{i}({\frak g}_{\alpha})
={\frak g}_{w_i(\alpha)}$. 

\vskip 2mm
We first consider the case when $\epsilon_i=1$. (Recall that, in
this case, $a_{i, \sigma(i)}=0$.)
Suppose $\alpha \neq \alpha_i$ and $\alpha \neq 2\alpha_i$. If we put 
$\bar r_{\bar i} =\prod_{k=0}^{N_i-1} r_{\sigma^k(i)}$,
then since $a_{i, \sigma^{k}(i)}=0$, we get 
$r_{\sigma^{k'}(i)} r_{\sigma^k(i)} = r_{\sigma^k(i)} r_{\sigma^{k'}(i)}$. 
It follows that 
$$\text {str} (\sigma^m| {\frak g}_{(-\alpha)})
=\text {str} (\bar r_{\bar i} \ \sigma^m \bar r_{\bar i}^{-1}|
{\frak g}_{(-\bar{w}_{\bar{i}} \alpha)})
=\text {str} (\sigma^m| {\frak g}_{(-\bar{w}_{\bar{i}} \alpha)}).$$
If $\alpha=\alpha_i$ and $i\in I^{even}$, then, 
since 
$$\text {str} (\sigma^m | {\frak g}_{(-\alpha_i)})
=\cases N_i \ \ & \text {if} \ m \equiv 0 \ (\text {mod} \, N_i), \\
0 \ \ & \text {otherwise},
\endcases
$$
we have 
$$\exp \left(\sum_{m=1}^{\infty} \frac{1}{m} 
\text {str} (\sigma^m| {\frak g}_{(-\alpha_i)}) 
E^{-m(\alpha_i)/ N_i} \right)
=(1-E^{-(\alpha_i)})^{-1},$$
and 
$$\bar w_{\bar i} (e^{-\rho}(1-E^{-(\alpha_i)})^{-1})
=e^{-\rho +(\alpha_i)} (1-E^{(\alpha_i)})^{-1}
=-e^{-\rho}(1-E^{-(\alpha_i)})^{-1}.$$
Similarly, if $\alpha=\alpha_i$ 
or $\alpha=2\alpha_i$ for $i\in I^{odd}$, then
\begin{equation*}
\begin{aligned}
& \exp \left(\sum_{m=1}^{\infty} \frac{1}{m} 
\text {str} (\sigma^m| {\frak g}_{(-\alpha_i)}) 
E^{-m(\alpha_i)/ N_i} \right) \\
& \times \exp \left(\sum_{m=1}^{\infty} \frac{1}{m} 
\text {str} (\sigma^m| {\frak g}_{(-2\alpha_i)}) 
E^{-m(2\alpha_i)/ N_i} \right) \\
&=(1-E^{-(\alpha_i)}) (1-E^{-(2\alpha_i)})^{-1},
\end{aligned}
\end{equation*}
and 
\begin{equation*}
\begin{aligned}
\bar w_{\bar i} &  (e^{-\rho}(1-E^{-(\alpha_i)})
(1-E^{-(2\alpha_i)})^{-1}) \\
& =e^{-\rho +(\alpha_i)} (1-E^{(\alpha_i)})
(1-E^{(2\alpha_i)})^{-1}\\
&=-e^{-\rho+(\alpha_i)}(1+e^{(\alpha_i)}) (1-e^{(2\alpha_i)})^{-1} \\
&=-e^{-\rho} (1-E^{-(\alpha_i)}) (1-E^{-(2\alpha_i)})^{-1}.
\end{aligned}
\end{equation*}

\vskip 2mm
Next, we consider the case when $\epsilon_i=2$. 
In this case, we have $a_{i, \sigma^{N_{i}/2}(i)}=-1$,
which implies that $i\in I^{even}$. 
Suppose $\alpha \neq \alpha_k$ and 
$\alpha \neq \alpha_{k} + \alpha_{\sigma^{N_{i}/2}(k)}$
for any $k=\sigma^{l}(i)$ $(1 \leq l\leq N_{i}/2)$.
Set $\bar r_{\bar i}=\prod_{l=0}^{N_{i}/2 -1} 
r_{\sigma^{l}(i)} r_{\sigma^{l+N_{i}/2}(i)} r_{\sigma^{l}(i)}$.
Since $a_{k,\sigma^{N_{i}/2}(k)}=-1$ implies 
$r_{k} r_{\sigma^{N_{i}/2}(k)} r_{k} 
= r_{\sigma^{N_{i}/2}(k)} r_{k} r_{\sigma^{N_{i}/2}(k)}$
(see, for example, \cite{KP}), we have $\sigma \bar r_{\bar i}
=\bar r_{\bar i} \sigma$ and 
$$\text {str} (\sigma^m| {\frak g}_{(-\alpha)})
=\text {str} (\bar r_{\bar i} \sigma^m \bar r_{\bar i}^{-1}|
{\frak g}_{(-\bar{w}_{\bar{i}} \alpha)})
=\text {str} (\sigma^m |{\frak g}_{(-\bar{w}_{\bar{i}} \alpha)}).$$
If $\alpha=\alpha_k$ or $\alpha=\alpha_k+ \alpha_{\sigma^{N_{i}/2}(k)}$
for some $k=\sigma^{l}(i)$ $(1 \leq l\leq N_{i}/2)$ 
(i.e. $(\alpha)=(\alpha_{i})$), then 
\begin{equation*}
\text {exp} \left( \sum_{m=1}^{\infty} \frac{1}{m} 
\text {str} (\sigma^m | {\frak g}_{(-\alpha_i)})
E^{-m(\alpha_i) / N_{i}} \right) 
=(1-E^{-(\alpha_i)})^{-1}
%=(1-E^{-\alpha_i -\alpha_{\sigma(i)}})^{-1},
\end{equation*}
\begin{equation*}
\text {exp} \left( \sum_{m=1}^{\infty} \frac{1}{m} 
\text {str} (\sigma^m | {\frak g}_{(-\alpha_i - \alpha_{\sigma^{N_{i}/2}(i)})})
E^{-2m(\alpha_i + \alpha_{\sigma^{N_{i}/2}(i)})/N_{i}} \right) 
=(1 + E^{-(\alpha_i)})^{-1},
\end{equation*}
and 
\begin{equation*}
\begin{aligned}
\bar w_{\bar i} & (e^{-\rho} (1-E^{(-\alpha_i)})^{-1}
(1+E^{(-\alpha_i)})^{-1})\\
&= - e^{-\rho} (1-E^{(-\alpha_i)})^{-1}
(1+E^{(-\alpha_i)})^{-1},
\end{aligned}
\end{equation*}
which completes the proof of our assertion. 
\qed }
\end{prop}

\vskip 2mm
\begin{prop} \label {Prop 9.2}
For $\Lambda \in P^{+}$ and $\bar w \in \bar W \subset W$, we have
\begin{equation} \label {9.6}
\bar w (\ch_{\sigma} V(\Lambda))=\ch_{\sigma} V(\Lambda).
\end{equation}

\noindent
{\it Proof.} \ {\rm  It suffices to prove (\ref {9.6})
for $\bar w=\bar w_{\bar i}$ $(\bar i \in \breve {I}^{re})$. 
Let $\pi: {\frak g} \rightarrow gl(V(\Lambda))$ denote
the representation of ${\frak g}$ on $V(\Lambda)$ and define
\begin{equation*}
\bar r^{\pi}_{\bar i}= \cases \prod_{k=0}^{N_i-1} r^{\pi}_{\sigma^k(i)} 
\ \ & \text {if} \ \epsilon_i=1, \\
%r^{\pi}_i r^{\pi}_{\sigma(i)} r^{\pi}_i 
\prod_{l=0}^{N_{i}/2 -1}
r_{\sigma^{l}(i)} r_{\sigma^{l+N_{i}/2}(i)} r_{\sigma^{l}(i)}
\ \ & \text {if} \ 
\epsilon_i=2.
\endcases
\end{equation*}
By the same argument in the preceding proposition, we have
$\sigma \bar r^{\pi}_{\bar i}=\bar r^{\pi}_{\bar i} \sigma$.
It follows that
$$\tr(\sigma| V(\Lambda)_{\lambda})
=\tr(\sigma|V(\Lambda)_{\bar w_{\bar i}\lambda}),$$
which proves our assertion. \qed
}
\end{prop}

\vskip 2mm

Now, we are in a position to state and prove the main result of this 
section. 

\begin{thm} \label {Thm 9.3} 
Let $\sigma$ be a Dynkin diagram automorphism of a generalized Kac-Moody
superalgebra ${\frak g}$ and let $\breve {\frak g}={\frak g}(\sigma)$ 
be the corresponding orbit Lie superalgebra.
Let $\phi: ({\frak h}^*)^{(0)} \rightarrow \widehat {\frak h}^*$ 
be the $\mathbb{C}$-linear isomorphism given by 
{\rm (\ref{8.7})}, and extend it to 
a $\mathbb {C}$-linear map $\phi: \mathbb {C}[[({\frak h}^*)^{(0)}]]
\rightarrow \mathbb{C}[[\widehat {\frak h}^*]]$ by 
$\phi(e^{\lambda})=e^{\phi(\lambda)}$ for $\lambda \in 
({\frak h}^*)^{(0)}$. 
Then, for a symmetric dominant integral weight $\Lambda \in P^+ \cap 
({\frak h}^*)^{(0)}$, we have
\begin{equation}
\phi(\ch_{\sigma} V(\Lambda))=\ch \breve{V}(\phi(\Lambda)),
\end{equation}
where $\breve{V}(\phi(\Lambda))$ is the irreducible highest weight
$\breve{\frak g}$-module with highest weight $\phi(\Lambda)$. 

\vskip 2mm

\noindent
{\it Proof.} \ {\rm Using the properties of Casimir operator given in
Proposition \ref{Prop:Casimir}, the same argument in \cite{K90} yields
\begin{equation} \label {9.8}
\ch_{\sigma} V(\Lambda)=\sum \Sb \lambda \le \Lambda \\
|\lambda+\rho|^2=|\Lambda+\rho|^2 \endSb c_{\lambda} 
\ch_{\sigma} M(\lambda).
\end{equation}
Since $\ch_{\sigma} M(\lambda)=e^{\lambda} \ch_{\sigma} M(0)$, 
by (\ref{9.8}), we have
\begin{equation} \label {9.9}
e^{\rho} R_{\sigma} \ch_{\sigma} V(\Lambda) 
=\sum \Sb \lambda \le \Lambda \\
|\lambda+\rho|^2=|\Lambda+\rho|^2 \endSb c_{\lambda} e^{\lambda+\rho},
\end{equation}
where 
$$R_{\sigma}=(\ch_{\sigma} M(0))^{-1} 
=\prod_{\alpha\in \Phi^{+}} \exp \left( -\sum_{m\ge 1} 
\frac{1}{m N_{\alpha}} \text {str}(\sigma^m | {\frak g}_{(-\alpha)}) 
E^{-m(\alpha)/ N_{\alpha}}\right).$$
We may assume that the sum in the right-hand side of (\ref{9.9})
is taken over the weights $\lambda$ which are invariant under $\sigma^*$.

\vskip 2mm
If $\bar w (\lambda + \rho)=\tau+\rho$, then
$c_{\tau}=\bar {\varepsilon}(\bar w) c_{\lambda}$.
Choose a weight $\tau$ such that $c_{\tau} \neq 0$. 
Then $\bar w (\tau +\rho) \le \Lambda + \rho$ for all $\bar w \in
\bar W$. 
Let $\lambda$ be an element of the set $\{\bar w(\tau + \rho)-\rho| \
\bar w \in \bar W \}$ such that $\text {ht} (\Lambda -\lambda)$ is
minimal.
Note that $(\lambda+\rho|\alpha_{i})\geq 0$ for all 
$\bar{i}\in \widehat{I}^{re} \setminus \breve{I}^{re}$ 
by definition of $\breve{I}$.
If $(\lambda + \rho | \alpha_i) < 0$ for some $\bar{i}\in \breve{I}^{re}$,
then
$\bar{w}_{\bar{i}}(\lambda+\rho)
=\lambda +\rho-\epsilon_{i}(\lambda+\rho)(h_{i})(\alpha_{i})$ and
$\text{ht}(\Lambda-\bar{w}_{\bar{i}}(\lambda)) < \text{ht}(\Lambda-\lambda)$,
which is a contradiction.
Hence, $(\lambda + \rho | \alpha_i) \geq 0$ for all $i\in I^{re}$.
Moreover, since $|\Lambda+\rho|^2=|\lambda+\rho|^2$, we have
\begin{equation} \label {9.10}
%\begin{aligned}
0= |\Lambda+\rho|^2 -|\lambda+\rho|^2 = (\Lambda|\Lambda-\lambda)
+(\lambda+2\rho|\Lambda-\lambda),
%&= (\Lambda | \sum n_i \alpha_i)+(\lambda+2\rho|\sum n_i \alpha_i),
%\end{aligned}
\end{equation}
where $\Lambda - \lambda =\sum_{i\in I} n_i \alpha_i$ lies in $Q^{+}$. 

\vskip 2mm 
If $i\in I^{re}$, then 
$$(\lambda+2\rho|\alpha_i) 
=(\lambda + \rho | \alpha_i) +(\rho | \alpha_i) >0.$$ 
If $i\in I^{im}$ and $n_i \ge 1$, then 
\begin{equation*}
\begin{aligned}
& (\lambda + 2\rho|\alpha_i)
=(\lambda + \alpha_i|\alpha_i)\\  
&=(\Lambda-\sum_{j\in I \setminus \{i\}} 
n_j \alpha_j -(n_i -1)\alpha_i|\alpha_i) \\
&=(\Lambda|\alpha_i)-\sum_{j\in I \setminus \{i\}} 
n_j (\alpha_j|\alpha_i) -(n_i -1)(\alpha_i|\alpha_i)\geq 0.
\end{aligned}
\end{equation*}                        
It follows that $n_i =0$ for all $i\in I^{re}$. 
That is, $\lambda=\Lambda-\sum_{i\in I^{im}}n_i \alpha_i$. 
>From (\ref{9.10}), for $k\in I^{im}$ with $n_k\geq 1$, we have
\begin{equation*}
(\lambda+2\rho|\alpha_k)=(\lambda+\alpha_k|\alpha_k)
=(\Lambda|\alpha_k)-\sum_{i\in I^{im}}n'_i (\alpha_i|\alpha_k)=0,
\end{equation*} 
where $n'_i=n_i$ for $i\neq k$ and $n'_k=n_k-1$.
Hence we have
\begin{eqnarray*}
(\Lambda|\alpha_k)&=&0,  \\
(\alpha_i|\alpha_k)&=&0 \ \mbox{  if $n_i\geq 1$ and $i\neq k$}, \\
(\alpha_k|\alpha_k)&=&0 \ \mbox{  if $n_k\geq 2$},
\end{eqnarray*}
and $k$ lies in $\breve{I}^{im}$.
Combining the above conditions and the fact that $n_i=n_{\sigma(i)}$ 
for all $i\in I$,we conclude that $\lambda$ has the form 
$$\lambda=\Lambda-\sum_{\bar i \in \breve{I}^{im}} n_{\bar i}(\alpha_i),$$ 
where 
\begin{eqnarray*}
((\alpha_i)|(\alpha_j))&=&0 \ \mbox{  if $\bar i\neq \bar j$ and 
$n_{\bar i}, n_{\bar j} \geq 1$},\\
(\Lambda|(\alpha_i))&=&0 \ \mbox{  if $n_{\bar i}\geq 1$}, \\
((\alpha_i)|(\alpha_i))&=&0 \ \mbox{  if $n_{\bar i} \geq 2$}. 
\end{eqnarray*} 
Let $\lambda,\lambda'$ be two such elements of the above form.
Suppose that there exists a $\bar{w}\in \bar{W}$ such that 
$\bar{w}(\lambda+\rho)=\lambda'+\rho$. 
Then, if $\lambda \in P^+$, we have $\bar{w}(\lambda)\leq \lambda$, 
since $\bar{w}(\lambda)$ is also 
a weight of the $\frak{g}$-module $V(\lambda)$.
Furthermore, by induction on the length of $\bar{w}$, 
we can verify  $\bar{w}(\rho)\leq \rho$.
Hence we obtain $\bar{w}(\lambda+\rho)\leq \lambda +\rho$ and
$\bar{w}(\lambda+\rho)=\lambda+\rho
-\sum_{\alpha \in \Phi^{+, re}} n_{\alpha}\alpha$
with $n_{\alpha}\geq 0$. On the other hand, if 
$\bar{w}(\lambda+\rho)=\lambda'+\rho$, we must have $n_{\alpha}=0$ for all
$\alpha \in \Phi^{+, re}$.  Hence $\bar{w}(\lambda+\rho)=\lambda+\rho$,
which implies $\lambda=\lambda'$.

\vskip 2mm
Let $\Phi^{+}(\Lambda)$ be the set defined in (\ref{6.10})
and let $\Phi^{+}(\Lambda)^{(0)}=\Phi^{+}(\Lambda)\cap (\frak{h}^*)^{(0)}$ 
be the set of elements
$\sum_{\bar i \in \breve{I}^{im}}n_{\bar i} (\alpha_i)$ 
such that $((\alpha_i)|(\alpha_j))=0$ if $n_{\bar i}, n_{\bar j} \geq 1$ 
for $\bar i \neq \bar j$,
$(\Lambda|(\alpha_i))=0$  if $n_{\bar i} \geq 1$, and
$((\alpha_i)|(\alpha_i))=0$ if $n_{\bar i} \geq 2$.
So (\ref {9.9}) can be written as follows\,:
\begin{equation*}
e^{\rho} R_{\sigma} \ch_{\sigma} V(\Lambda) 
=\sum_{\bar{w}\in \bar{W}}\bar{\epsilon}(\bar{w})
 \sum_{\beta\in \Phi^{+}(\Lambda)^{(0)}}
 c_{\Lambda-\beta}e^{\bar{w}(\Lambda-\beta+\rho)}.
\end{equation*}
Now, it remains to determine $c_{\Lambda-\beta}$ 
for $\beta\in S_{\Lambda}^{(0)}$.
First, we observe that $\ch_{\sigma} V(\Lambda)$ contains no term of the form
$e^{\Lambda-\sum_{i\in I^{im}}n_i\alpha_i}$, 
where $(\Lambda|\alpha_i)=0$ 
for all $i$ such that $n_i\geq 1$.
So, in order to determine $c_{\Lambda-\beta}$, we have only to compute
the coefficient of $e^{\Lambda-\beta}$ in $R_{\sigma}$.
Note that $((\alpha_i)|(\alpha_j))=0$ implies that 
$\alpha_i+\alpha_j \not\in \Phi$.

\vskip 2mm
For each $\bar i \in \breve{I}^{im}$, we have
\begin{equation*}
\exp \left( -\sum_{m\ge 1} 
\frac{1}{m} \text {str}(\sigma^m | {\frak g}_{(-\alpha_i)}) 
E^{-m(\alpha_i)/ N_i}\right)=(1-E^{(-\alpha_i)})^{m_i}.
\end{equation*}
For $\beta=\sum_{\bar i \in  \breve{I}^{im}} n_{\bar i}(\alpha_i)$, 
using the formal power series expansions
$$(1-x)^n=\sum_{k=0}^{n}\binom {n}{k}(-1)^k x^k$$
and 
$$(1+x)^{-n}=\sum_{k=0}^{\infty}\binom {n+k-1}{k}(-1)^k x^k,$$
we obtain  
\begin{equation*}
c_{\Lambda-\beta}=
\prod_{\bar i \in \breve{I}^{even}}\binom {m_{\bar i}}{n_{\bar i}} 
(-1)^{n_{\bar i}} \prod_{\bar i \in \breve{I}^{odd}}
\binom {m_{\bar i}+n_{\bar i}-1}{n_{\bar i}} (-1)^{n_{\bar i}}
=(-1)^{|\phi(\beta)|}\varepsilon(\phi(\beta)).
\end{equation*}
When $\Lambda=0$, we get 
\begin{equation*}
R_{\sigma}=\sum_{\bar{w}\in \bar{W}}\bar{\varepsilon}(\bar{w})
\sum_{\lambda\in \Phi^{+}(0)^{(0)}}c_{-\lambda}
e^{\bar{w}(\lambda+\rho)-\rho}.
\end{equation*}
Since $\phi((\alpha_i))=\widehat{\alpha}_{\bar i}$ for 
$\bar i\in \breve{I}^{im}$,  by Lemma \ref{Lemma 8.1}, we have  
$((\alpha_i)|(\alpha_j))=0$ 
if and only if $(\widehat{\alpha}_{\bar i}|\widehat{\alpha_{\bar j}})=0$ and
$(\Lambda|(\alpha_i))=0$ 
if and only if $(\phi(\Lambda)|\widehat{\alpha}_{\bar i})=0$.
It follows that $\phi(\Phi^{+}(\Lambda)^{(0)})
=\breve {\Phi}^{+}(\phi(\Lambda)) 
\subset \widehat{\frak{h}}^*$.
Therefore we obtain
\begin{equation*}
\begin{aligned}
\phi(\ch_{\sigma} V(\Lambda))&= \frac
{\sum_{\bar w \in \bar W} \bar {\varepsilon}(\bar w)
\sum_{\beta \in \Phi^{+}(\Lambda)^{(0)}} c_{\Lambda-\beta}
e^{\phi(\bar w (\Lambda -\beta + \rho)-\rho)}}
{\sum_{\bar w \in \bar W} \bar {\varepsilon}(\bar w)
\sum_{\gamma \in \Phi^{+}(0)^{(0)}} c_{-\gamma}
e^{\phi(\bar w (-\gamma + \rho)-\rho)}} \\
&=\frac{ \sum_{\breve{w} \in \widehat {W}} \varepsilon (\widehat {w}) 
\sum_{\widehat {\beta} \in \breve {\Phi}^{+}(\phi(\Lambda))}
c_{\Lambda - \phi^{-1}(\widehat {\beta})} 
e^{\breve {w}(\phi(\Lambda)-\widehat {\beta} +\rho)-\rho}}
{\sum_{\breve {w} \in \breve {W}} \varepsilon (\breve {w}) 
\sum_{\widehat {\gamma} \in \breve {\Phi}^{+}(0)}
c_{ -\phi^{-1}(\widehat {\gamma})} 
e^{\breve {w}(-\widehat {\gamma} +\rho)-\rho}} \\
&=\frac{ \sum_{\breve {w} \in \breve {W}} \varepsilon (\breve {w}) 
\sum_{\widehat {\beta} \in \breve{\Phi}^{+}(\phi(\Lambda))}
(-1)^{|\widehat {\beta}|} \varepsilon(\widehat {\beta}) 
e^{\breve {w}(\phi(\Lambda)-\widehat {\beta} +\rho)-\rho}}
{\sum_{\breve {w} \in \breve {W}} \varepsilon (\breve {w}) 
\sum_{\widehat {\gamma} \in \breve{\Phi}^{+}(0)}
(-1)^{|\widehat {\gamma}|} \varepsilon(\widehat {\gamma}) 
e^{\breve {w}(-\widehat {\gamma} +\rho)-\rho}} \\
&=\ch \breve{V}(\phi(\Lambda)).
\end{aligned}
\end{equation*}
\qed }
\end{thm}

\vskip 2mm

\begin{cor} \label {Cor 9.4}
\begin{equation} \label {9.11}
\begin{aligned}
& \phi(R_{\sigma})  = \prod_{\alpha \in \breve{\Phi}^{+}} 
(1-E^{-\alpha})^{\text {sdim} \, \breve {\frak g}_{\alpha}}, \\
& \phi(\ch_{\sigma} M(\Lambda))  =\ch \breve {M} (\phi(\Lambda)). \\
\end{aligned}
\end{equation}

\noindent
{\it Proof.} \ {\rm In the proof of Theorem 9.3, we have shown that
$$R_{\sigma} =\sum_{\bar w \in \bar W} \bar {\varepsilon}(\bar w)
\sum_{\gamma \in \Phi^{+}(0)^{(0)}} c_{-\gamma} 
e^{\bar w (-\gamma + \rho)-\rho}.$$
Hence we have 
\begin{equation*}
\begin{aligned}
\phi(R_{\sigma})&= \sum_{\bar w \in \bar W} \bar {\varepsilon}(\bar w)
\sum_{\gamma \in \Phi^{+}(0)^{(0)}} c_{-\gamma} 
e^{\phi(\bar w (-\gamma + \rho)-\rho)} 
= \sum_{\breve {w} \in \breve {W}} {\varepsilon}(\breve {w})
\sum_{\widehat{\gamma} \in \breve{\Phi}^{+}(0)} 
c_{-\phi^{-1}(\widehat{\gamma})} 
e^{\breve {w} (-\widehat {\gamma} + \rho)-\rho} \\
&= \sum_{\breve {w} \in \breve {W}} {\varepsilon}(\breve {w})
\sum_{\widehat{\gamma} \in \breve{\Phi}^{+}(0)} (-1)^{|\widehat {\gamma}|} 
\varepsilon(\widehat {\gamma}) 
e^{\breve{w} (-\widehat {\gamma} + \rho)-\rho} 
=\prod_{\alpha \in \breve {\Phi}^{+}} 
(1-E^{-\alpha})^{\text {sdim} \, \breve {\frak g}_{\alpha}}.
\end{aligned}
\end{equation*}
Since $\ch_{\sigma} M(\Lambda)=e^{\Lambda} R_{\sigma}^{-1}$, we obtain
\begin{equation*}
\phi(\ch_{\sigma} M(\Lambda)) 
=e^{\phi(\Lambda)} \prod_{\alpha \in \breve {\Phi}^{+}} 
(1-E^{-\alpha})^{-\text {sdim} \, \breve {\frak g}_{\alpha}}
=\ch \breve {M}(\phi(\Lambda)).
\end{equation*}
}
\qed
\end{cor}

\vskip 2mm

Set $\widehat {Q}=\bigoplus_{\bar i \in \widehat {I}} \mathbb {Z}
(\widehat {\alpha}_{\bar i} \big/ \epsilon_i N_i)$, and 
define a $\mathbb {Z}$-linear map  
$\widehat {\phi}: Q \rightarrow \widehat {Q}$ by
$$\widehat {\phi}(\alpha_{i}) 
=\widehat {\alpha}_{\bar i} \big/ \epsilon_i N_i \ \ 
\text {for} \ i\in I, \ k=0,1, \cdots, N_i-1.$$
Then $\widehat {\phi}$ is an extension of 
$\phi|_{Q \cap ({\frak h}^*)^{(0)}}$ to $Q$. From now on, 
we will also use $\phi$ for $\widehat {\phi}$. 

\vskip 2mm
For each $\beta \in \widehat {Q}$, define the
subspace ${\frak g}_{[\alpha]}$ of ${\frak g}$ to be 
$${\frak g}_{[\beta]}=\bigoplus_{\alpha \in \phi^{-1}(\beta)}
{\frak g}_{\alpha}.$$
Clearly, $[{\frak g}_{[\alpha]}, {\frak g}_{[\beta]}]
\subset {\frak g}_{[\alpha+\beta]}$, and hence ${\frak g}$
becomes a $\widehat {Q}$-graded Lie superalgebra with the coloring 
of $\widehat {Q}$ induced by that of $Q$. 
Note that the Dynkin diagram automorphism $\sigma$ preserves the
$\widehat {Q}$-gradation of ${\frak g}$. 
With respect to this gradation, the generalized denominator identity 
for the Dynkin diagram automorphism $\sigma$ of the generalized Kac-Moody 
superalgebra ${\frak g}$ is the same as
\begin{equation} \label {9.12}
\prod_{\beta \in \widehat {Q}^{+}} \exp \left( -\sum_{m=1}^{\infty}
\frac{1}{m} \text {str} (\sigma^m|{\frak g}_{[-\beta]}) E^{-m\beta} \right)
=1- \text {sch}_{\sigma} H({\frak g}_{-}).
\end{equation}
On the other hand, we have
\begin{equation} \label {9.13}
\begin{aligned} 
\phi(R_{\sigma})&= \prod_{\alpha \in \Phi^{+}} \exp \left(
-\sum_{m=1}^{\infty} \frac{1}{m N_{\alpha}} \text {str}
(\sigma^m | {\frak g}_{(-\alpha)}) 
E^{- m \phi((\alpha)) / N_{\alpha}} \right) \\
&=\prod_{\alpha \in \Phi^{+}} \exp \left( 
-\sum_{m=1}^{\infty} \frac{1}{m N_{\alpha}} \text {str}
(\sigma^m | {\frak g}_{(-\alpha)}) E^{- m \phi(\alpha)} \right) \\
&=\prod_{\beta \in \widehat {Q}^{+}} \exp \left( 
-\sum_{m=1}^{\infty} \frac{1}{m} \text {str}
(\sigma^m | {\frak g}_{[-\beta]}) E^{- m \beta } \right).
\end{aligned}
\end{equation}
Therefore, by combining (\ref{9.12}) and (\ref{9.13})
with Corollary \ref{Cor 9.4}, we obtain:

\begin{thm} \label {Thm 9.5}
The generalized denominator identity of the $\widehat {Q}$-graded 
generalized Kac-Moody superalgebra ${\frak g}$ for the Dynkin diagram 
automorphism $\sigma$ is the same as the denominator identity of the
orbit Lie superalgebra $\breve {\frak g}$. 
\qed
\end{thm}

\vskip 2mm
\noindent
{\it Remark.} \ (a) \ For the Kac-Moody algebras and generalized 
Kac-Moody algebras, 
Theorem \ref{Thm 9.3} was proved in \cite{FSS} and \cite {FRS}.
They also gave a list of Dynkin diagram automorphisms of 
Kac-Moody algebras ${\frak g}$ of finite growth and corresponding  
orbit Lie algebras $\breve {\frak g}$ (\cite{FSS}). 

\vskip 2mm
(b) \ Let ${\frak D}_{\sigma}({\frak g})$ denote the
generalized denominator identity of the $\widehat {Q}$-graded 
generalized Kac-Moody superalgebra ${\frak g}$ for the Dynkin 
diagram automorphism $\sigma$ and ${\frak D}(\breve {\frak g})$
the denominator identity of the orbit Lie superalgebra 
$\breve {\frak g}$. Then, Theorem \ref{Thm 9.5} implies
\begin{equation} \label {9.14}
{\frak D}_{\sigma}({\frak g})
={\frak D}(\breve {\frak g}).
\end{equation}

\vskip 2mm
(c) For the affine Kac-Moody superalgebras
introduced in \cite{K78}, we have
\begin{equation*}
\begin{aligned}
{\frak D}_{\sigma}(A^{(2)}(0, 2n-1))
& = {\frak D}(B^{(1)}(0, n-1)) \ \ \text {when $n\ge 2$ and 
$\sigma$ has order 2},\\
{\frak D}_{\sigma} (C^{(2)}(3))
& ={\frak D}(A^{(4)}(0,2)) \ \  \text {when $\sigma$ has order 2}, \\
{\frak D}_{\sigma} (C^{(2)}(0, 2n+1))
&= {\frak D}(A^{(4)}(0,2n)) \ \  \text {when $n\ge 2$ and $\sigma$ 
has order 2}, \\
{\frak D}_{\sigma} (C^{(2)}(0, 2n+2))
&={\frak D}(B^{(1)}(0,n)) \ \ \text {when $n\ge 2$ and $\sigma$ 
has order 2}.
\end{aligned}
\end{equation*}

\vskip 2mm
%(c) \ The identity (\ref{9.14}) can be written as 
%$$\prod_{\beta \in \widehat {Q}^{+}} \exp \left( -\sum_{m=1}^{\infty}
%\frac{1}{m} \text {str} (\sigma^m| {\frak g}_{[-\beta]}) E^{-m \beta} \right)
%=\prod_{\alpha \in \widehat {Q}^{+}} 
%(1-E^{-\alpha})^{\text {sdim} \, {\frak g}(\sigma)_{\alpha}},$$
%where ${\frak g}(\sigma)=\widehat {\frak g}$ is the orbit Lie superalgebra 
%corresponding to the Dynkin diagram automorphism $\sigma$ of ${\frak g}$. 
%By taking the logarithm of both sides, we have
%$$\sum_{\beta \in \widehat {Q}^{+}} \sum_{m=1}^{\infty}
%\frac{1}{m} \text {str} (\sigma^m| {\frak g}_{[-\beta]}) E^{-m \beta}
%=  \sum_{\alpha \in \widehat {Q}^{+}} \sum_{n=1}^{\infty} 
%\frac {1}{n} \text {sdim} \, {\frak g}(\sigma)_{\alpha} E^{-n \alpha}.$$
%If follows that 
%$$\sum \Sb k>0 \\ \beta = k \tau \endSb
%\frac{1}{k} \text {str} (\sigma^k| {\frak g}_{[-\tau]})
%=\sum \Sb a>0 \\ \alpha = a \gamma \endSb 
%\frac{1}{a} \text {sdim} \, {\frak g}(\sigma)_{\gamma}.$$
%Hence, by using M\"obius inversion, we obtain 
%\begin{equation} \label {9.15}
%\text {sdim} \, {\frak g}(\sigma)_{\alpha}
%=\sum_{ad | \alpha } \frac {1} {ad} \mu(d) \, \text {str}
%(\sigma^{ad} | {\frak g}_{[-\alpha /ad]})
%\end{equation}
%and
%\begin{equation} \label {9.16}
%\text {str}(\sigma | {\frak g}_{[-\alpha]})
%=\sum_{ad|\alpha} \frac{1}{ad} \mu(d) \, \text {sdim} \, 
%{\frak g}(\sigma^d)_{\alpha / ad},
%\end{equation}
%where $\mu$ denotes the classical M\"obius function. 
The identity (\ref{9.14}) can be rewritten  as 
$$\prod_{\beta \in \widehat {Q}^{+}} \exp \left( -\sum_{m=1}^{\infty}
\frac{1}{m} \text {str} (\sigma^m| {\frak g}_{[-\beta]}) E^{-m \beta} \right)
=\prod_{\alpha \in \widehat {Q}^{+}} 
(1-E^{-\alpha})^{\text {sdim} \, {\frak g}(\sigma)_{\alpha}},$$
where ${\frak g}(\sigma)=\breve {\frak g}$ is the orbit Lie superalgebra 
corresponding to the Dynkin diagram automorphism $\sigma$ of ${\frak g}$.
Let us consider ${\frak g}(\sigma^d)$ for a nonnegative integer $d$. 
Then we have the following objects as in the case of ${\frak g}(\sigma)$:

(i) $\widehat{I}(d)$,$\breve{I}(d)$ $N_i (d)$, $\epsilon_i (d)$, 
$\widehat{A}(d)$, $\breve{A}(d)$, 
$\widehat {Q}(d)=\bigoplus_{\bar i \in \widehat {I}(d)} \mathbb {Z}
(\widehat {\alpha}_{\bar i} \big/ \epsilon_i (d)  N_i (d))$,

(ii) a $\mathbb {Z}$-linear map  
$\phi_d : Q \rightarrow \widehat {Q}(d)$, where
$\phi_d (\alpha_{i}) 
=\widehat {\alpha}_{\bar i} \big/ \epsilon_i (d) N_i  (d)\ \ 
\text {for} \ i\in I, \ k=1, \cdots, N_i (d).$

Applying the above equation to ${\frak g}(\sigma^d)$ , we get
\begin{equation}\label{9.15}
\prod_{\beta \in \widehat {Q}^{+}(d)} \exp \left( -\sum_{m=1}^{\infty}
\frac{1}{m} \text {str} (\sigma^{md}|{\frak g}_{[-\beta]_{d}})
 E^{-m \beta} \right)
=\prod_{\alpha \in \widehat {Q}^{+}(d)} 
(1-E^{-\alpha})^{\text {sdim} \, {\frak g}(\sigma^d)_{\alpha}},
\end{equation}
where ${\frak g}_{[\beta]_d}=\bigoplus_{\alpha \in \phi_{d}^{-1}(\beta)}
{\frak g}_{\alpha}.$
Also, if we define a map $\tau$ of $\widehat{I}(d)$ by 
$\tau(\bar{i})=\overline{\sigma(i)}$,
then it induces a $\mathbb{Z}-$ linear map 
$\phi^{d} : \widehat{Q}(d) \rightarrow \widehat{Q}$
such that $\phi^{d} \phi_{d}(\beta)=\phi(\beta)$ for all $\beta \in Q$.
By taking $\phi^{d}$, (\ref{9.15}) can be written as 
\begin{equation}\label{9.16}
\prod_{\beta \in \widehat {Q}^{+}} \exp \left( -\sum_{m=1}^{\infty}
\frac{1}{m} \text {str} (\sigma^{md}|{\frak g}_{[-\beta]})
 E^{-m \beta} \right)
=\prod_{\alpha \in \widehat {Q}^{+}} 
(1-E^{-\alpha})^{\text {sdim} \, {\frak g}(\sigma^d)_{[\alpha]^{d}}},
\end{equation}
where ${\frak g}(\sigma^d)_{[\alpha]^{d}}
=\bigoplus_{\beta \in (\phi^{d})^{-1}(\alpha)}
{\frak g}(\sigma^d)_{\beta}.$
By taking the logarithm on both sides, we have
$$\sum_{\beta \in \widehat {Q}^{+}} \sum_{m=1}^{\infty}
\frac{1}{m} \text {str} (\sigma^{md}| {\frak g}_{[-\beta]}) E^{-m \beta}
=  \sum_{\alpha \in \widehat {Q}^{+}} \sum_{n=1}^{\infty} 
\frac {1}{n} \text {sdim} \, {\frak g}(\sigma^d)_{[\alpha]^d} E^{-n \alpha}.$$
It follows that 
$$\sum \Sb k>0 \\ \beta = k \tau \endSb
\frac{1}{k} \text {str} (\sigma^{dk}|{\frak g}_{[-\tau]})
=\sum \Sb a>0 \\ \beta = a \gamma \endSb 
\frac{1}{a} \text {sdim} \, {\frak g}(\sigma^d)_{[\gamma]^d}.$$
If we put $T_{\sigma^d}(\beta) 
=\sum_{a|\beta}\frac{1}{a}\text {sdim} \, {\frak g}(\sigma^d)_{[\beta/a]^d}$
for a nonnegative integer $d$ and $\beta\in \widehat{Q}$,
then by M\"obius inversion 
\begin{equation} \label {Mobius}
\begin{aligned}
\text {str}(\sigma^{k} | {\frak g}_{[-\alpha]})
&=\sum_{d|\alpha}\frac{1}{d}\mu(d)T_{\sigma^{dk}}(\alpha/d) \\
&=\sum_{ad|\alpha} \frac{1}{ad} \mu(d) \, \text {sdim} \, 
{\mathfrak g}(\sigma^{dk})_{[\alpha / ad]^{dk}},
\end{aligned}
\end{equation}
where $\mu$ denotes the classical M\"obius function.

Let ${\frak g}^{\sigma}$ be the fixed point 
subalgebra of ${\frak g}$ for a Dynkin diagram automorphism
$\sigma$. Then it is also graded by $\widehat{Q}$ and 
by (\ref{Mobius}) we get

\vskip 2mm
\begin{prop} \label{Prop 9.6} 
For each $\beta \in \widehat{Q}$,
\begin{equation}\label{9.17}
\text {sdim}{\frak g}^{\sigma}_{\beta}
=\frac{1}{|\sigma|}\sum_{k=1}^{|\sigma|}
\sum_{ad|\beta} \frac{1}{ad} \mu(d) \, 
\text {sdim}{\frak g}(\sigma^{dk})_{[\beta / ad]^{dk}},
\end{equation}
\end{prop}

\noindent
{\it Proof.} \ 
Note that ${\frak g}^{\sigma}_{\beta}$ is a representation of a 
finite cyclic group $\langle \sigma \rangle$. 
Hence $\text {dim}{\frak g}^{\sigma}_{\beta}$ 
is the multiplicity of the trivial representation of 
$\langle \sigma \rangle$ 
in ${\frak g}^{\sigma}_{\beta}$ and it is the same as the 
average of traces of
$\tau\in \langle \sigma \rangle$ on  ${\frak g}^{\sigma}_{\beta}$.
\qed

\vskip 5mm
\noindent
{\it Remark.} \ By the similar argument in \cite{B88}, we can see that
$\frak{g}^{\sigma}$ is a central extension of 
a generalized Kac-Moody superalgebra. However, 
the corresponding Borcherds-Cartan matrix  is not known in general. 
Also, the decomposition of $\frak{g}^{\sigma}$ with respect to
the $\widehat{Q}-$gradation is coarser than its root space decomposition.

\vskip 2cm

\small{

}
\end{document}